\documentclass[12 pt]{khbook}
\usepackage{amssymb}
\usepackage{amsmath}
\usepackage{eepic}
\usepackage{epsfig}
\usepackage{graphicx}
\usepackage[
english]{babel}
\author{Raf Cluckers}
\title{Cell decomposition and $p$-adic integration}
\date{}
\usepackage{fancyheadings}
\renewcommand{\chaptermark}[1]{\markboth{Chapter \thechapter.\ #1}{}}
\renewcommand{\sectionmark}[1]{\markright{\thesection\ #1}}
\newcommand{\clearemptydoublepage}{\clearpage{\pagestyle{empty}\cleardoublepage}}

\makeindex
\theoremstyle{plain}
\newtheorem{theorem}{Theorem}[section]
\newtheorem{proposition}[theorem]{Proposition}
\newtheorem{prop}[theorem]{Proposition}
\newtheorem{lemma}[theorem]{Lemma}
\newtheorem{cor}[theorem]{Corollary}

\theoremstyle{definition}
\newtheorem{definition}[theorem]{Definition}

\newtheorem{remark}[theorem]{Remark}
\newtheorem{remarks}[theorem]{Remarks}

\newtheorem{theorema}{Theorema}[chapter]
\newtheorem{definitie}[theorema]{Definitie}

\numberwithin{equation}{section}
\newcommand{\qed}{\rule{0mm}{0mm}\hfill$\square$\par}
\newenvironment{proof}{{\par\it Proof.}}{\qed}
\newcommand{\prooftitle}[1]{\par\textit{#1}\\}
\newenvironment{abstract}{{\textbf{Abstract.}}}{}

\newcommand{\Z}{\mathbb{Z}}
\newcommand{\R}{\mathbb{R}}
\newcommand{\C}{\mathbb{C}}
\newcommand{\N}{\mathbb{N}}
\newcommand{\Q}{\mathbb{Q}}
\newcommand{\F}{\mathbb{F}}
\newcommand{\PPP}{\mathbb{P}}
\newcommand{\PP}{\mathcal{P}}
\newcommand{\BK}{\mathcal{B}_K}
\newcommand{\TK}{K\langle X\rangle}

\newcommand{\OF}{\mathcal{O}(F)}

\newcommand{\SSS}{\mathcal{S}}
\newcommand{\aaa}{\alpha}
\newcommand{\bbb}{\beta}

\newcommand{\ddd}{\delta}
\newcommand{\lP}{\lambda P_n}

\newcommand{\Lm}{\mathcal{L}}

\newcommand{\Ldan}{\mathcal{L}_K}

\newcommand{\Lp}{\mathcal{L}_{Pres}}
\newcommand{\Lmac}{\mathcal{L}_{Mac}}

\DeclareMathOperator{\sq}{\square}
\newcommand{\Pm}{\mathcal{P}}

\newcommand{\M}{\ensuremath{\mathcal{M}}}
\newcommand{\Def}{\ensuremath{\mathcal{D}ef}}
\newcommand{\OOs}{\Omega}
\newcommand{\OOL}{\Omega_{\rm an}}
\newcommand{\OOK}{\Omega_{\rm an}}
\newcommand{\OOLp}{\Omega_{\rm simple}}
\newcommand{\cP}{\mathcal{P}}
\newcommand{\abs}[1]{\lvert#1\rvert}
\newcommand{\Lv}{\ensuremath{\mathcal{L}_{\rm v}}}
\newcommand{\cL}{\ensuremath{\mathcal{L}}}
\newcommand{\Lr}{\ensuremath{\mathcal{L_{\rm ring}}}}
\newcommand{\Lpas}{\cL_{\rm Pas}}
\newcommand{\ac}{\mathrm{ac}}
\renewcommand{\char}{{\rm char}}

\begin{document}

\frontmatter \clearemptydoublepage

\begin{titlepage}
\begin{center}
\large{KATHOLIEKE UNIVERSITEIT LEUVEN}\\
\normalsize{FACULTEIT WETENSCHAPPEN\\
DEPARTEMENT WISKUNDE}\\[15mm]
\includegraphics[height=32mm]{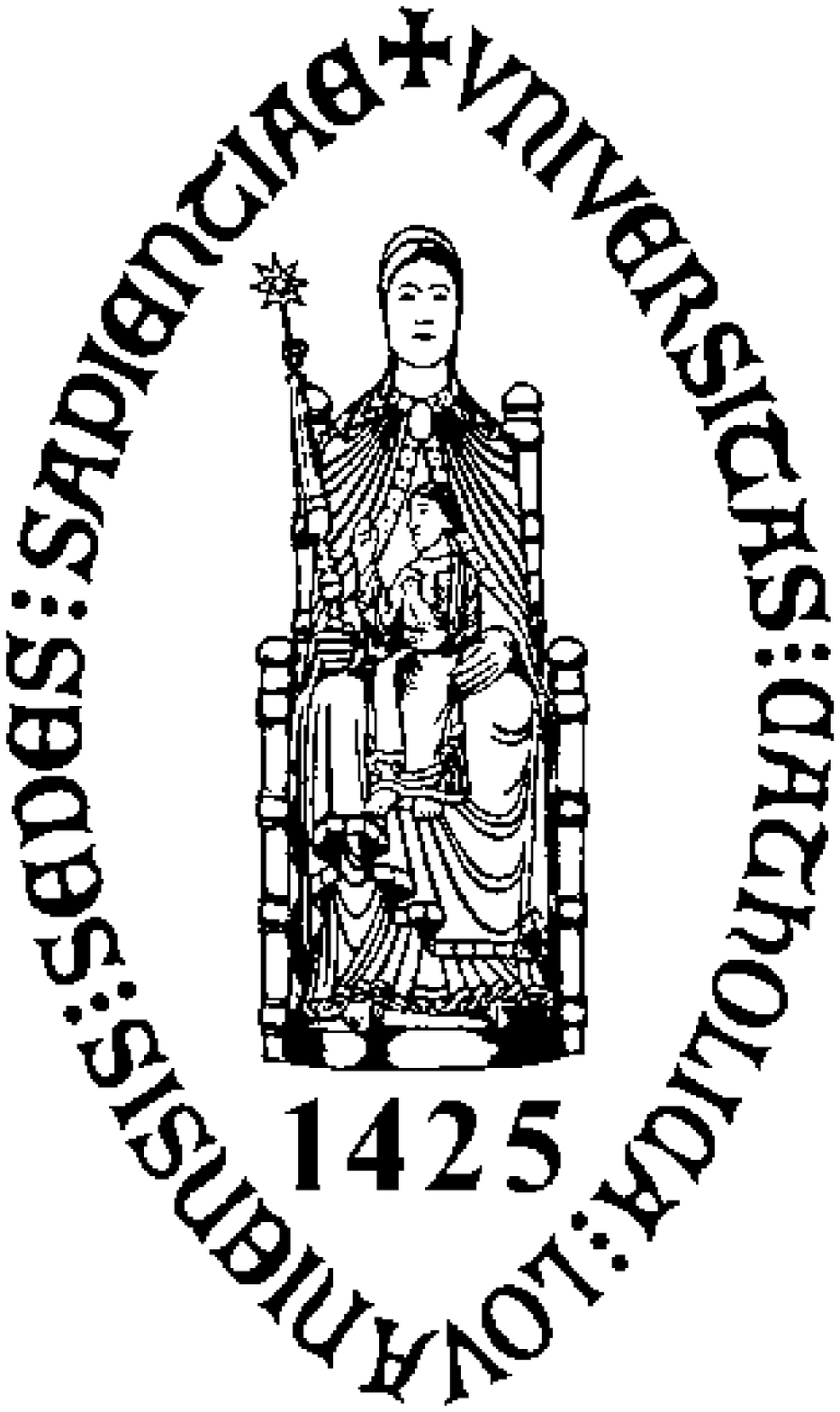}\\[15mm]
\textbf{\Large{
Cell decomposition and $p$-adic integration }\\[15mm]
\large{
Raf Cluckers\\[40mm]}}
\end{center}
Promotor: \textsc{Jan Denef}\\
18 december 2002
\end{titlepage}


\hyphenation{se-mi-al-ge-bra-isch se-mi-al-ge-bra-i-sche
ver-za-me-lingen in-te-gra-tie Gro-then-dieck}

\setcounter{page}{1}

\chapter{Voorwoord}
\begin{tabbing}
This is to position quotation tab's \= and \= author \kill
 \> Voorbij een bepaald punt\\
 \> kun je niet meer terug.\\
 \> Probeer dit punt te bereiken!\\
 \> \> \textsc{Franz Kafka}  \\
\end{tabbing}
De start van het doctoraat was erg onzeker,  zou er een beurs voor
mij zijn of niet, ik moest afwachten. Dit lijkt al erg lang
geleden, maar het is goed te beseffen dat ik geluk gehad heb:
enkele weken na mijn afstuderen in juni 1999 kreeg ik een beurs
toegewezen en kon ik er tegenaan gaan.
\par
Bij het be\"eindigen van mijn doctoraat besef ik dat ik enorm heb
genoten van het uitoefenenen van het beroep van wiskundige, en van
mijn steentje te kunnen bijdragen tot de wetenschap wiskunde. Het
doen van wiskundig onderzoek is een activiteit waarvan ik erg ben
gaan houden en is een grote drijfveer voor mij geweest. Laat ik
hier kort aan toevoegen dat ik hoop nog vele steentjes te kunnen
bijdragen, en nog menige jaren te kunnen meedraaien in het
wiskundig onderzoek. Natuurlijk hoop ik ook steeds veel tijd en
energie te kunnen vrijmaken voor vrienden, familie en allerlei
interessants  en leuks buiten het onderzoek, want dat is en blijft
voor mij het belangrijkste. Academisch onderzoek kan echter
voorgenoemde soms be\-moei\-lij\-ken: onzekerheid over toekomst,
onderzoek doen aan buitenlandse universiteiten, de druk om te
bouwen aan een goed curriculum vitae, enzovoort. Ik kijk
optimistisch naar de toekomst, en met tevredenheid naar de
voorbije jaren van onderzoek.
\par
 Samenvattend: het zijn fijne jaren
geweest te werken aan een doctoraat.
\par
\noindent \textit{Leuven, december  2002} \hfill \textsc{Raf
Cluckers}

\clearemptydoublepage

\chapter{Dankwoord}
\begin{tabbing}
This is to position quotation \= and \= author \kill
 \> Niet de kracht maar de volharding\\
 \> staat aan de oorsprong\\
 \> van een grote verwezenlijking.\\
 \> \> \textsc{Samuel Johnson}\\
\end{tabbing}
Ik bedank mijn promotor Prof.~dr.~Jan Denef voor de
enthousiasmerende gesprekken, voor de mooie problematieken waarin
hij mij ingeleid heeft en voor de vrijheid die hij me gaf, steeds
ondersteund door sterk advies.
\par
Graag wil ik ook mijn familie bedanken, mijn moeder en vader in
het bijzonder, die me steeds enorm hebben gesteund in mijn keuzes
en me gestimuleerd hebben om nieuwe horizonten te verkennen. Ook
alle vrienden, die me naar het einde van het doctoraat toe
misschien al eens minder zagen, wil ik bedanken voor het fijne en
ontspannende samenzijn. Vooral wil ik de persoon bedanken die mij
het dichtst bij staat: Erika, van wie ik hou, en wiens steun dag
na dag onschatbaar groot is geweest.
\par
Mijn dank gaat ook naar de leden van mijn afdeling voor de
interessante discussies, het leuke samenzijn, het vele gedeelde
werk binnen en buiten het onderzoek en voor de fijne sfeer op de
afdeling. In dezelfde optiek dank ik iedereen van het departement
wiskunde. Dank ook aan de KULeuven en het F.W.O.-Vlaanderen die
mij de mogelijkheid gaven in goede omstandigheden aan een
doctoraat te werken. Ik ben blij in een maatschappij te leven waar
men zoveel middelen vrijmaakt voor wetenschappelijk onderzoek.
\par
Veel dank ben ik verschuldigd aan jurylid Prof.~dr.~W.~Veys en
bu\-reau\-ge\-noot Wannes voor de vele wiskundige interacties,
exploraties en improvisaties.
\par
Ik dank ook hartelijk de verdere leden van de jury
Prof.~dr.~A.~Macintyre, Prof.~dr.~L.~Lipshitz,
Prof.~dr.~K.~Dekimpe en Prof.~dr.~L.~Verstraelen voor het
beoordelen van de thesis.
\par
Tenslotte, wat kan ik beter zeggen dan  het eindproduct voor te
leggen en iedereen  veel en intens leesplezier toe te wensen.

\clearemptydoublepage

\selectlanguage{english}

\tableofcontents
\clearemptydoublepage

\chapter{Introduction}
\begin{tabbing}
This is to position quotation tabbings \= and \= author \kill
 \> Beyond a certain point,\\
 \> one cannot return.\\
 \> Try to reach this point!\\
 \> \> \textsc{Franz Kafka}  \\
\end{tabbing}
The thesis consists of a collection of published papers and
preprints (each chapter representing one paper), together with an
extra introduction and some extra explanations. Since each chapter
contains an introduction where the general context is sketched and
the relevant definitions are given, we only give a short general
introduction. In this general introduction, we give an overview of
the topics of the thesis, we sketch some general context, and give
some references where more background can be found.
\par
To describe the thesis in one sentence: the research is mostly
inspired by the work of my thesis advisor, J.~Denef; it is a
combination of model theory and number theory in the lines of
J.~Igusa.
\section{General context}
It is typical for model theory, to study geometrical properties of
several algebraic objects. The most well-known topics coincide of
course with other domains of mathematics, like the study of
algebraic, semialgebraic, analytic and subanalytic sets in $\R^n$.
Structures that are studied in the thesis are:
\begin{quote}
I) The integers as an additive ordered group. This structure gives
rise to an algebra of sets closed under taking finite unions,
intersections, projections and Cartesian products. These sets  are
called Presburger sets. More precisely, a Presburger set is a
subset of $\Z^m$ for some $m$, which can be obtained by taking
finite unions, finite intersections, complements, projections, and
Cartesian products of $\N\subset\Z$ and of cosets of submodules of
$\Z^{m+e}$, where $e\geq 0$. A function between Presburger sets is
called a Presburger function if its graph is a Presburger set. A
Presburger set $X\subset \Z^m$ can always be described by a finite
disjunction of conjunctions of conditions of the forms
\begin{eqnarray*}
 f(x)\leq 0, &&
 g(x) \equiv 0 \bmod n,
\end{eqnarray*}
where $f$ and $g$ are  polynomials  of degree $\leq 1$ over $\Z$,
$n$ an integer, as proven by M.~Presburger in \cite{Presburger}.
 \par
II) The $p$-adic numbers as a field give rise to an algebra of
sets, also closed under the operations mentioned in I). These sets
are called semi-algebraic sets. More precisely, a semialgebraic
set is a subset of $\Q_p^m$ for some $m$ which can be obtained by
taking finite unions, finite intersections, complements and
projections of zerosets $X\subset \Q_p^{m+e}$ of polynomials over
$\Q_p$, where $e\geq0$. A function between semialgebraic sets is
called a semialgebraic function if its graph is a semialgebraic
set. A semialgebraic set $X\subset \Z^m$ can always be described
by a finite  disjunction of conjunctions of conditions of the form
\begin{eqnarray*}
 f(x)\in \lambda P_n &&
\end{eqnarray*}
where $P_n$ consists of the $n$-th powers in
$K^\times=K\setminus\{0\}$ for $n>0$, $\lambda$ is in $K$ and $f$
is a polynomial over $\Q_p$, as proven by A.~Macintyre in
\cite{Mac}.
\par
III) One can enrich the structure given in II) by putting on the
$p$-adic numbers also an analytic structure, given by power series
which converge on $\Z_p^m$, $m>0$. This  gives rise to an algebra
of subanalytic sets. More precisely, a subanalytic  set is a
subset of $\Q_p^m$ for some $m$ which can be obtained by taking
finite unions, finite intersections, complements and projections
of zerosets $X\subset \Q_p^{m+e}$ of  polynomials over $\Q_p$ and
of zerosets $Y\subset\Z_p^{m+e'}\subset\Q_p^{m+e'}$ of converging
power series (see chapter \ref{chap:cell} for precise
definitions). In \cite{DvdD}, J.~Denef and L.~van den Dries give a
clear and simple description of subanalytic sets, which is stated
in chapter \ref{chap:cell}.
\par
IV) We also study fields, carrying a valuation to an ordered group
and containing a Henselian valuation ring. Examples are
$\F_q((t))$, $\Q_p$, $\Q_q((t))$, $\R((t))$, $\C((t))$, and so on.
Using first order languages containing the field operations and
possibly other symbols, we study the algebra of definable sets and
definable functions.
\end{quote}
\par
\section{Overview and context of the main results}
For the collection of semialgebraic sets  on the real numbers, the
following classification is known \cite{vdD}:
 \begin{quote} \textit{There exists a real semi-algebraic bijection
between two real semi-algebraic sets if and only if they have the
same dimension and Euler characteristic.}
\end{quote}
 This is one motivation to look for a non-trivial Euler
characteristic on the $p$-adic semialgebraic sets and is related
to the question posed by L.~B\'elair whether there exists a
semialgebraic bijection from $\Q_p$ to $\Q_p\setminus\{0\}$. In
joint work with D.~Haskell, such a bijection is constructed, see
\cite{CH} or chapter \ref{chap:CH}. It follows that there exists
no nontrivial Euler characteristic on the $p$-adic semialgebraic
sets! Equivalently, the Grothendieck ring (in the sense of logic)
of the $p$-adic numbers is trivial.
 Subsequently, a $p$-adic analogue of the classification mentioned
 above was proven by the author; it is as follows.
 \begin{quote} \textit{There exists a  semi-algebraic bijection
between two infinite $p$-adic semi-algebraic sets if and only if
they have the same dimension.}
\end{quote}
This is the content of chapter \ref{chap:classification}.
\par
In Igusa theory, one often encounters parametrized integrals, as
well $p$-adic, real as over $\F_q((t))$ and $\C$.  We study
$p$-adic parametrized integrals like
\[
I(\lambda)=\int_{(\lambda,y)\in A}|f(\lambda,y)||dy|,
\]
where $A\subset \Z_p^{m+e}$ is a semialgebraic set, $f$ a
polynomial over $\Q_p$ in the variables $(\lambda,y)$, $|dy|$ a
Haar measure and $|\cdot|$ stands for the $p$-adic norm. More
generally, one can take a semialgebraic function for $f$, or even
a subanalytic function. For example, Weil's local singular series
are functions which can be written as a parametrized integral, see
chapter \ref{chap:decay}. In the case that $\lambda$ stands for
one $p$-adic variable, the function $I$ can be studied by
classical methods like resolution of singularities. However, these
methods fail up to present to describe $I$ when $\lambda$ is a
tuple of variables. It was an important insight of J.~Denef to
formulate and prove a cell decomposition theorem for semialgebraic
sets and semialgebraic functions to study $I$ as well as other
questions related to $p$-adic integrals. For example, he obtains
the rationality of several Poincar\'e series  by calculating
integrals by means of the cell decomposition theorem \cite{Denef}.
Using cell decomposition, Denef \cite{Denef1} proves
 \begin{quote}\textit{ The
function $I$ is in the $\Q$-algebra generated by functions of the
form $\lambda\to |g(\lambda)|$ and $\lambda\to v(h(\lambda))$,
where $g$ and $h$ are semialgebraic functions, $h$ nowhere zero,
and $v(\cdot)$ stands for the $p$-adic order.}
 \end{quote}
If $f$ is a subanalytic function and $A$ is a subanalytic set, one
can define a parametrized integral $I_{\mathrm an}$ like $I$
above. The analogue result for $I_{\mathrm an}$ to belong to a
similar algebra was an open problem, conjectured to be true by
Denef in \cite{Denef1}. The underlying open problem was to find a
subanalytic cell decomposition theorem. This is obtained by the
author in \cite{Ccell}, see chapter \ref{chap:cell}, and
represents a crucial part of the thesis. The proof of this
subanalytic cell decomposition relies on recent work in
\cite{vdDHM} and uses rigid analysis. Accordingly, the conjecture
about $I_{\mathrm an}$ is solved:
 \begin{quote}\textit{ The
function $I_{\mathrm an}$ is in the $\Q$-algebra generated by
functions of the form $\lambda\to |g(\lambda)|$ and $\lambda\to
v(h(\lambda))$, where $g$ and $h$ are subanalytic functions.}
 \end{quote}
This subanalytic analogue on parametrized integrals is obtained
using the subanalytic cell decomposition, see chapter
\ref{chap:cell}.

\par
It was noticed by the author that some of the proofs and results
about $p$-adic numbers, like for example the classification
result, hold also for Presburger sets. In order to be able to
formulate a classification result, a good notion of dimension of
Presburger sets is needed. In \cite{CPres} (see chapter
\ref{chap:Pres}), the dimension of a Presburger set is defined
using the notion of algebraic closure in the sense of logic and
the independence property of algebraic closure for Presburger
sets, proven in \cite{Wagner}. The main result of chapter
\ref{chap:Pres} is a classification of Presburger sets up to
Presburger isomorphism.
 \begin{quote}\textit{There exists a Presburger bijection
between two infinite Presburger sets if and only if they have the
same dimension.}
 \end{quote}
Also elimination of imaginaries for Presburger sets is proven,
which can be stated as follows.
 \begin{quote}\textit{Given a Presburger-definable equivalence relation $\sim$ on $\Z^m$ for some $m$,
there exists a Presburger function $F:\Z^m\to \Z^k$ for some $k$
such that $F(x)=F(y)$ if and only  if $x\sim y$.}
 \end{quote}
\par
In chapters two and three, Grothendieck rings of fields, carrying
a valuation and containing a Henselian ring, are calculated. A
Grothendieck ring corresponds to a universal Euler characteristic,
in the sense that every Euler characteristic on definable sets
factors through the natural map to the Grothendieck ring. For
several fields and several first order languages, the triviality
of the Grothendieck ring is proven.
\par
In the last chapter, we study multivariate exponential sums. Let
$f=(f_1,\ldots,f_r)$ be a dominant polynomial mapping.  Let
$a=(a_1,\ldots,a_r)$ be a tuple of integers, relatively prime to
$p$ and let $m=(m_1,\ldots,m_r)$ be a tuple of positive integers.
For each $a$ and $m$ we define the normalized Kloosterman sum
\[
E(a,m)=\frac{1}{p^{\prod_i m_i}}\sum_{x_i\in \Z_{p}/p^{m_i}}\exp
(2\pi i(
\frac{a_1f_1(x)}{p^{m_1}}+\ldots+\frac{a_rf_r(x)}{p^{m_r}})).
\]
If $r=1$ this sum has been studied intensively by Igusa
\cite{Igusa3}, Loxton \cite{Lox}, Lichtin \cite{Lichtin},
\cite{Lichtin2}, and many others. If $r=1$ and if $f$ is
homogeneous, Igusa proved that $E$ goes to zero when $m_1$ goes to
infinity and that $E$ can be bounded by $cp^{\alpha m_1}$ where
$c$ and $\alpha<0$ are real numbers, ensuring a quick decay of the
exponential sum when $m_1$ goes to infinity. Igusa also relates
$\alpha$ to numerical data of a resolution of singularities of
$f=0$. Later on, several generalizations have been obtained, for
example, Lichtin obtains a non trivial decay for $E$ when $r=2$
and the decay $\alpha$ is linked to fine geometric invariants of
the polynomials. We study $E$ for arbitrary $r$ and obtain the
following general result
\begin{quote}\textit{
There exist real numbers $\alpha<0$ $c>0$ such that
$$
|E(a,m)|<c\min \{p^{ m\alpha},1\},\qquad\mbox{ for all $a$ and
$m$},
$$
where $m=\max\{m_1,\ldots,m_r\}$.}
\end{quote}
This is a generalization of Igusa's result on the decay of
exponential sums to the multidimensional case, but without
relating $\alpha$ to any geometric invariants. Instead of
resolution of singularities, we use integration techniques
involving  cell decomposition to prove this qualitative decay of
Kloosterman sums.
\par
 By the previous work on subanalytic sets and cell decomposition, all the
arguments hold also in the analytic setting. Thus, whenever
$f=(f_1,\ldots,f_r)$ is a dominant mapping given by convergent
power series $f_i$, an analogue qualitative decay of the analytic
Kloosterman sums holds.
\par
Future projects to pursue are, among others, to extend the
complete Igusa theory to the multi-variate setting, and to develop
the model theory of valued fields to a higher extend. I hope this
will come into reach in the next years. It is my strong believe
that the language of definable functions, definable sets, cell
decomposition, and other concepts of model theory are very
powerful tools to pursue these projects.
\section{References for general theory}
\par
 For a general introduction to model
theory, I refer to \cite{CK}, \cite{Hodges}, and \cite{Marker}.
\par
For a general introduction to the basic theory of $p$-adic
numbers, see \cite{Koblitz} and \cite{Prestel}.
\par
For a general introduction on Igusa theory, the reader may consult
\cite{Igusa1}, \cite{Igusa2}, \cite{Igusa3}, \cite{Igusa:intro}
and \cite{DenefBour}. For the first steps towards a multi-variate
Igusa-theory, I refer to \cite{Lichtin}, \cite{Lichtin2},
\cite{Denef1} and chapter \ref{chap:cell}.
\par
Of course, there are many other good books which are not mentioned
here, wherefore our apologies. Further, most of the  definitions
will be given and general context will be explained in the
chapters itself.

\clearemptydoublepage \mainmatter

\chapter{Presburger sets and p-minimal fields.}\label{chap:Pres}
\begin{abstract}
\footnote{This chapter corresponds to \cite{CPres}.} We prove a
cell decomposition theorem for Presburger sets and introduce a
dimension theory for $Z$-groups with the Presburger structure.
Using the cell decomposition theorem we obtain a full
classification of Presburger sets up to definable bijection. We
also exhibit a tight connection between the definable sets in an
arbitrary p-minimal field and Presburger sets in its value group.
We give a negative result about expansions of Presburger
structures and prove  uniform elimination of imaginaries for
Presburger structures within the Presburger language.
\end{abstract}
\section{Introduction}
At the ``Alg\`ebre, Logique et Cave Particuli\`ere'' meeting in
Lyon (1995), A.~Pillay posed the question of whether there exists
some dimension theory for $Z$-groups with the Presburger structure
which would give rise to a classification of all Presburger sets
up to definable bijection, possibly using other invariants as
well. In this paper we answer this question of Pillay: we classify
the Presburger sets up to definable bijection
(Thm.~\ref{classification2}), using as only classifying invariant
the (logical) algebraic dimension. In order to prove this
classification, we first formulate a cell decomposition theorem
for Presburger groups (Thm.~\ref{cell decomp}) and a
rectilinearisation theorem for the definable sets
(Thm.~\ref{recti}). Also a rectilinearisation theorem depending on
parameters is proven (Thm.~\ref{param recti}).

Expansions of Presburger groups have recently been studied
intensively. One could say that on the one hand one looks for
(concrete) expansions which remain decidable and have bounded
complexity, and on the other hand different kinds of minimality
conditions (like coset-minimality, etc.) are used to characterize
general classes of expansions (see e.g. \cite{Wagner},
\cite{Point}). In section \ref{section Michaux} we examine
expansions of Presburger groups satisfying natural kinds of
minimality conditions.

In \cite{Haskell}, D.~Haskell and D.~Macpherson defined the notion
of p-minimal fields, as a $p$-adic counterpart of o-minimal
fields. A p-minimal field always is a $p$-adically closed field,
and its value group is a $Z$-group. Interactions between definable
sets in a given $p$-adically closed field and Presburger sets in
its value group have been studied in the context of $p$-adic
integration for several p-minimal structures (see \cite{Denef1}).
In Theorem \ref{application}, we exhibit a close connection
between definable sets in arbitrary p-minimal fields and
Presburger sets in the corresponding value groups.

In the last section, we use the cell decomposition theorem in an
elementary way to obtain uniform elimination of imaginaries for
Z-groups
 without introducing extra sorts.
\\
\subsection*{Terminology and notation.}
In this paper $G$ always denotes a $Z$-group, i.e.~a group which
is elementary equivalent\footnote{Given a langauge $\cL$, two
structures $M$ and $M'$ for $\cL$ are called elementary equivalent
if every $\cL$-sentence holds in $M$ if and only if it holds in
$M'$.}
 to the integers $\Z$ in the
Presburger language $\Lp = \langle
+,\leq,\{\equiv(\mathrm{mod}{n})\}_{n>0},0,1 \rangle$ where
$\equiv(\mathrm{mod}{n})$ is the equivalence relation in two
variables modulo the integer $n>0$. We call $(G,\Lp)$ a Presburger
structure and we write $H$ for the nonnegative elements in $G$. By
a Presburger set, function, etc., we mean a $\Lp$-definable set,
function, etc., and by \emph{definable} we always mean definable
with parameters (otherwise we say $\emptyset$-definable,
$S$-definable, etc.). We call a set $X\subset G^m$ \emph{bounded}
if there is a tuple $z\in H^m$ such that $-z_i\leq x_i\leq z_i$
for each $x\in X$ and $i=1,\ldots,m$. For $k\leq m$ we write
$\pi_k:G^m\to G^k$ for the projection on the first $k$ coordinates
and for $X\subset G^{k+n}$ and $x\in\pi_k(X)$ we write $X_x$ for
the fiber $\{y\in G^n\mid (x,y)\in X\}$. We recall that the theory
$\mathrm{Th}(\Z,\Lp)$ has definable Skolem
functions\footnote{Given a structure $M$ for a language $\cL$, we
say that $M$ has definable Skolem functions if for any definable
set $A\subset M^{n+e}$ there exists a definable function $f$ from
the projection $\pi_n(A)$ to $M^e$, such that $(a,f(a))\in A$ for
each $a\in \pi_n(A)$. We say that a theory has definable Skolem
functions if and only if all of its models have.}, quantifier
elimination\footnote{We say that a structure $M$ for some language
$\cL$ has quantifier elimination if every $\cL$-formula is
equivalent in $M$ to a quantifier free $\cL$-formula.} in $\Lp$
and is decidable \cite{Presburger}.
\section{Cell Decomposition Theorem}
We prove a cell decomposition theorem for Presburger structures,
by first proving it in dimension 1 and subsequently using a
compactness argument. An elementary arithmetical proof can also be
given, using techniques like in the proof of Lemma 3.2 in
\cite{Denef3}, but our proof has the advantage that it goes
through in other contexts as well (see section \ref{section
Michaux} and \ref{section p}). As always, $G$ denotes a $Z$-group.
 \begin{definition}\label{linearity}
We call a function $f:X\subset G^m\to G$ \emph{linear} if there is
a constant $\gamma\in G$ and integers $a_i$,  $0\leq c_i < n_i$
for $i=1,\ldots,m$ such that $x_i-c_i\equiv 0 \pmod{n_i}$ and
\[
f(x)=\sum_{i=1}^m a_i(\frac{x_i-c_i}{n_i})+\gamma.
\]
for all $x=(x_1,\ldots,x_m)\in X$. We call $f$ \emph{piecewise
linear} if there is a finite partition $\Pm$ of $X$ such that all
restrictions $f|_A$, $A\in \Pm$ are linear. We speak analogously
of linear and piecewise linear maps $g:X\to G^n$.
 \end{definition}
The following definition fixes the notion of (Presburger) cells.
 \begin{definition}\label{def:cell2}
A cell of type $(0)$ (also called a $(0)$-cell) is a point
$\{a\}\subset G$. A $(1)$-cell is a set with infinite cardinality
of the form
 \begin{equation}\label{cell dim 1}
  \{x\in G\mid \alpha\sq_1 x\sq_2 \beta,\ x\equiv c \pmod{n}\},
 \end{equation}
with $\alpha,\beta$ in $G$, integers $0\leq c<n$ and $\sq_i$
either $\leq$ or no condition. Let $i_j\in\{0,1\}$ for
$j=1,\ldots,m$ and $x=(x_1,\ldots,x_m)$. A
$(i_1,\ldots,i_{m},1)$-cell is a set $A$ of the form
 \begin{equation}\label{cell}
  A = \{(x,t)\in G^{m+1} \mid x\in D,\ \alpha(x)\sq_1 t \sq_2
  \beta(x),\ t\equiv c \pmod{n}\},
 \end{equation}
with $D=\pi_{m}(A)$ a $(i_1,\ldots,i_m)$-cell,
 $\alpha,\ \beta:D\to G$ linear functions, $\sq_i$ either
$\leq$ or no condition and integers $0\leq c<n$ such that the
cardinality of the fibers $A_x=\{t\in G\mid (x,t)\in A\}$ can not
be bounded uniformly in $x\in D$ by an integer.
 \\
 A $(i_1,\ldots,i_{m},0)$-cell is a set of the form
 \[
\{(x,t)\in G^{m+1} \mid x\in D,\ \alpha(x)=t\},
 \]
with $\alpha:D\to G$ a linear function and $D\subset G^{m}$ a
$(i_1,\ldots,i_m)$-cell.
 \end{definition}
 \begin{remarks}\label{remark cell}
\begin{itemize}
\item[(i)]
Although we consider in Definition \ref{def:cell2} a condition on
the cardinality of fibers, the type of a cell does not alter if
one takes elementary extensions.\footnote{Given a model $M$ for
some language $\cL$, an elementary extension of $M$ is a model
$M'$ containing $M$ such that each $\cL$-sentence with
parameters from $M$ holds in $M$ if and only if it holds in $M'$.}\\
\item[(ii)]\label{remark cell item}
 To an infinite $(i_1,\ldots,i_m)$-cell $A\subset G^m$ we can associate
(as in \cite{vdD}) a projection $\pi_A:G^m\to G^k$ such that the
restriction of  $\pi_A$ to $A$ gives a bijection from $A$ onto a
$(1,\ldots,1)$-cell $A'\subset G^k$. Also,  a $(i_1, \ldots,
i_m)$-cell
is finite if and only if $i_1 = \cdots = i_m = 0$, and then it is a singleton.\\
\item[(iii)]
Let $A$ be a $(i_1,\ldots,i_{m},1)$-cell as in Eq.~(\ref{cell}),
then it is clear that a linear function $f:A\to G$ can be written
as
 \begin{equation}\label{lin dim 1}
 f(x,t)=a(\frac{t-c}{n})+\gamma(x),\quad  (x,t)\in A,
 \end{equation}
with $a$ an integer, $\gamma:D\to G$ a linear function and $c,n,D$
as in Eq.~(\ref{cell}).
 \end{itemize}\end{remarks}
 \begin{theorem}[Cell Decomposition]\label{cell decomp}
Let $X\subset G^m$ and $f:X\to G$ be $\Lp$-definable. Then there
exists a finite partition $\Pm$ of $X$ into cells, such that the
restriction $f|_A:A\to G$ is linear for each cell $A\in\Pm$.
Moreover, if $X$ and $f$ are $S$-definable, then also the parts
$A$ can be taken $S$-definable.
 \end{theorem}
 \begin{proof}[Proof by induction on $m$.]
If $X\subset G$, $f:X\to G$ are $\Lp$-definable, then Theorem
\ref{cell decomp} follows easily by using quantifier elimination
and elementary properties of linear congruences. Alternatively,
the more general Thm.~4.8 of \cite{Point} can be used to prove
this one dimensional version (see also Proposition~\ref{piecewise
linear} below).
 Let $X\subset G^{m+1}$ and  $f:X\to G$ be
$\Lp$-definable, $m>0$. We write $(\sq_1, \sq_2)\in
\{\leq,\emptyset\}^2$ to say that $\sq_1$, resp.~$\sq_2$,
represents either the symbol $\leq$ or no condition. Let  $\SSS$
be the set $\Z\times\{(n,c)\in\Z^2 \mid 0\leq c <
n\}\times\{\leq,\emptyset\}^2$.
 For any $d=(a,n,c,\sq_1, \sq_2) \in \SSS$
and $\xi=(\xi_1,\xi_2,\xi_3) \in G^3$ we define a Presburger
function $F_{(d,\xi)}$ as follows
\[
 \{t\in G\mid \xi_1 \sq_1 t \sq_2 \xi_2,\ t\equiv c
\pmod{n}\}\to G:t\mapsto a(\frac{t-c}{n})+\xi_3.
\]
The domain  $Dom(F_{(d,\xi)})$ of such a function $F_{(d,\xi)}$ is
either empty, a $(1)$-cell or a finite union of $(0)$-cells. For
fixed $k>0$ and  $d\in \SSS^k$, let $\varphi_{(d,k)}(x,\xi)$ be a
Presburger formula in the free variables $x=(x_1,\ldots,x_m)$ and
$\xi=(\xi_1,\ldots,\xi_k)$, with
$\xi_i=(\xi_{i1},\xi_{i2},\xi_{i3})$, such that $G\models
\varphi_{(d,k)}(x,\xi)$ if and only if the following are true:
 \begin{itemize}
\item[(i)] $x\in\pi_m(X)$,\\
\item[(ii)] the collection  of the domains $Dom(F_{(d_i,\xi_i)})$ for $i=1,\ldots,k$ forms
 a partition of the fiber $X_x\subset
G$,\\
\item[(iii)] $F_{(d_i,\xi_i)}(t)=f(x,t)$ for
each  $t\in Dom(F_{(d_i,\xi_i)})$ and $i=1,\ldots,k$.
 \end{itemize}
Now we define for each $k$ and $d \in \SSS^k$ the set
\[
B_{(d,k)}= \{x\in G^m\mid \exists \xi\quad
\varphi_{(d,k)}(x,\xi)\}.
\]
Each set $B_{(d,k)}$ is $\Lp$-definable. Also, the (countable)
collection $\{B_{(d,k)}\}_{k,d}$ covers $\pi_m(X)$ since each
$x\in\pi_m(X)$ is in some $B_{(d,k)}$ by the induction basis. We
can do this construction in any elementary extension of $G$, so by
logical compactness we must have that finitely many sets of the
form $B_{(d,k)}$ already cover $\pi_m(X)$. Consequently, we can
take Presburger sets $D_1,\ldots,D_s$ such that $\{D_i\}$ forms a
partition of $\pi_m(X)$ and each $D_i$ is contained in a set
$B_{(d,k)}$ for some $k$ and $k$-tuple $d$. For each
$i=1,\ldots,s$, fix a $k$ and $k$-tuple $d$ with $D_i\subset
B_{(d,k)}$, then we can define the Presburger set
\[
\Gamma_i=\{(x,\xi)\in D_i\times G^{3k}\mid
\varphi_{(d,k)}(x,\xi)\},
\]
with $\pi_m(\Gamma_i)=D_i$ by construction. Since the theory
$\mathrm{Th}(G,\Lp)$ has definable Skolem functions, we can choose
definably for each $x\in D_i$ tuples $\xi\in G^{3k}$ such that
$(x,\xi)\in\Gamma_i$. Combining it all, it follows that there
exists a finite partition $\Pm$ of $X$ consisting of Presburger
sets of the form
\[
A=\{(x,t)\in G^{m+1}\mid x\in C,\ \alpha(x)\sq_1 t \sq_2
\beta(x),\ t\equiv c \pmod{n}\},
\]
such that $f|_A$ maps $(x,t)\in A$ to
$a(\frac{t-c}{n})+\gamma(x)$; with $\alpha,\beta,\gamma:C\to G$
and $C\subset G^m$ $\Lp$-definable, $\sq_i$ either $\leq$ or no
condition, integers $a$, $0\leq c<n$ and $\pi_m(A)=C$. The theorem
now follows after applying the induction hypothesis to $C$ and
$\alpha,\beta,\gamma:C\to G$ and partitioning further.
\end{proof}
\section[Dimension theory]{Dimension theory for Presburger\\ structures}
Any Presburger structure satisfies the exchange property for
algebraic closure. This is a corollary of a more general result in
\cite{Wagner} but can also be proven  using the cell decomposition
theorem in an elementary way. In particular this yields an
algebraic dimension
function on the Presburger sets in the following (standard) way.\\
 \begin{definition}
Let $X\subset G^m$ be $A$-definable for some finite set $A$ by a
formula $\varphi(x,a)$ where $a=(a_1,\ldots,a_s)$ and
$A=\{a_1,\ldots, a_s\}$, then the (algebraic) dimension of $X$,
written $\mathrm{dim}(X)$, is the greatest integer $k$ such that
in some elementary extension $\bar G$ of $G$ there exists
$x=(x_1,\ldots,x_m)\in \bar G^{m}$ with $\bar G\models
\varphi(x,a)$ and
\[
\mathrm{rk}(x_1,\ldots,x_m,a_1,\ldots,a_s)-\mathrm{rk}(a_1,\ldots,a_s)=k,
\]
where $\mathrm{rk}(B)$ of a set $B\subset \bar G$ is the
cardinality of a maximal algebraically independent subset of $B$
(in the sense of model theory, see \cite{Hodges}).
\end{definition}
This dimension function is independent of the choice of a set of
defining parameters $A$ and  the following properties of algebraic
dimension are standard.
 \begin{proposition}\label{behaves well}
\begin{itemize}
\item[(i)]For Presburger sets $X,Y\subset G^m$ we have
$\mathrm{dim}(X\cup Y)=\max(\mathrm{dim}X,\mathrm{dim}Y)$.\\
\item[(ii)]Let $f:X\to G^m$ be $\Lp$-definable, then
\[
\mathrm{dim}(X)\geq \mathrm{dim}(f(X)).
\]
\end{itemize}
 \end{proposition}
The dimension of a cell $C$ is directly related to the type of $C$
(see Lemma \ref{algdim cell}). Also, if we have a Presburger set
$X$ and a finite partition $\Pm$ of $X$ into cells, the dimension
of $X$ is directly related to the types of the cells in $\Pm$ (see
Cor.~\ref{dim, partitie}).
 \begin{lemma}\label{algdim cell}
Let $C\subset G^m$ be a $(i_1,\ldots,i_m)$-cell, then
$\mathrm{dim}(C)=i_1+\ldots+i_m$.
 \end{lemma}
 \begin{proof} For a $(0)$- and a $(1)$-cell this is clear.
Possibly after projecting, we may suppose that $C\subset G^m$ is a
$(1,\ldots,1)$-cell. The Lemma follows now from the definition of
the type of a cell using induction on $m$ and a compactness
argument.
 \end{proof}
 \begin{cor}\label{dim, partitie}
For any Presburger set $X\subset G^m$ and any finite partition
$\Pm$ of $X$ into cells we have
 \begin{eqnarray*}
 \mathrm{dim}(X)
  & = &
  \max\{i_1+\ldots+i_m\mid C\in\Pm,\mbox{ $C$ is a
  $(i_1,\ldots,i_m)$-cell}\}\nonumber
  \\
  & = &
  \max\{i_1+\ldots+i_m\mid X \mbox{ contains a
  $(i_1,\ldots,i_m)$-cell}\}.\label{Eqdim1}
 \end{eqnarray*}
 \end{cor}
 \begin{proof} The first equality is a consequence of Lemma~\ref{algdim cell}
and Proposition~\ref{behaves well}. To prove the second equality
we take a $(i_1,\ldots,i_m)$-cell $C\subset X$ such that
$i_1+\ldots+i_m$ is maximal. By the cell decomposition we can
obtain a partition $\Pm$ of $X$ into cells such that $C\in\Pm$.
Now use the previous equality to finish the proof.
 \end{proof}
 \begin{remark}
It is also possible to take the second equality of corollary
\ref{dim, partitie} as the definition for the dimension of a
Presburger set and to proceed similarly as in \cite{vdD} by van
den Dries to develop a dimension theory for Presburger structures.
 \end{remark}
 \section{Classification of Presburger sets}\label{sectionclassification}
The cell decomposition theorem provides us with the technical
tools to classify the $\emptyset$-definable Presburger sets up to
$\Lp$-definable bijection. The key step to this classification is
a rectilinearisation theorem, which also has a parametric
formulation. We recall that $G$ denotes a $Z$-group and $H=\{x\in
G\mid x\geq0\}$, we also write $H^0=\{0\}$. Also notice that a set
$A$ is $\emptyset$-definable if and only if $A$ is $\Z$-definable,
to be precise, definable over $\Z\cdot 1\subset G$.
 \begin{theorem}[Rectilinearisation]\label{recti}
Let $X$ be a $\emptyset$-definable Presburger set, then there
exists a finite partition $\Pm$ of $X$ into $\emptyset$-definable
Presburger sets, such that for each $A\in \Pm$ there is an integer
$l\geq 0$ and a $\emptyset$-definable linear bijection $f:A\to
H^l$.
 \end{theorem}
 \begin{proof} We give a proof by induction on $\mathrm{dim}X$. If $\mathrm{dim} X=0$ then
$X$ is finite and the theorem follows, so we choose a Presburger
set $X$ with $\mathrm{dim} X=m+1$, $m\geq0$. Any $\Lp$-definable
object occurring in this proof will be $\emptyset$-definable; we
will alternately apply $\emptyset$-definable linear bijections and
partition further. By the cell decomposition theorem and possibly
after projecting (see the remark after Definition
\ref{def:cell2}), we may suppose that $X$ is a $(1,\ldots,1)$-cell
contained in $G^{m+1}$, so we can write
\[
X=\{(x,t)\in G^{m+1}\mid x\in D,\ \alpha(x)\sq_1 t\sq_2 \beta(x),\
t\equiv c\pmod{n}\},
\]
with $x=(x_1,\ldots,x_m)$, $\pi_m(X)=D\subset G^m$ a
$(1,\ldots,1)$-cell, integers $0\leq c<n$, $\alpha,\beta:D\to G$
$\emptyset$-definable linear functions and $\sq_i$ either $\leq$
or no condition. By induction we may suppose that $D=H^m$. If both
$\sq_1$ and $\sq_2$ are no condition, the theorem follows easily,
so we may suppose that one of the $\sq_i$, say $\sq_1$, is $\leq$.
Moreover, after a linear transformation
$(x,t)\mapsto(x,\frac{t-c}{n})$ we may assume that $c=0$ and
$n=1$, then we can apply the following linear bijection
\[
   f:X \to A:(x,t) \mapsto (x_1,\ldots,x_m,t-\alpha(x)),
\]
onto
 \begin{equation*}
A  =  \{(x,t)\in H^{m+1}\mid t \sq_2 \beta(x)-\alpha(x)\}.
 \end{equation*}
Because $\beta(x)-\alpha(x)$ is a linear function from $H^m$ to
$G$ there are integers $k_i$ such that
\begin{equation}\label{some cell}
A  =  \{(x,t)\in H^{m+1}\mid t \sq_2 k_0+\sum_{i=1}^m k_i x_i\}.
\end{equation}
Moreover, since $\pi_m(A)=H^m$, all integers $k_i$ must be
nonnegative. We proceed by induction on $k_1\geq0$. If $k_1=0$
then $A=H\times\{(x_2,\ldots,x_m,t)\in H^m\mid t\leq
k_0+\sum_{i=2}^m k_i x_i\}$ and the theorem follows by induction
on the dimension. Now suppose $k_1>0$, then we partition $A$ into
two parts
 \begin{eqnarray*}
 A_1 & = & \{(x,t)\in H^{m+1}\mid t\leq x_1-1\},\\
 A_2 & = & \{(x,t)\in H^{m+1}\mid x_1\leq t\leq k_0+\sum_{i=1}^m k_i x_i\},
 \end{eqnarray*}
where $\pi_m(A_2)= H^m$ and $\pi_m(A_1)=\{x\in H^m \mid 1\leq
x_1\}$. We apply the linear bijection
\[
A_2\to B:(x,t)\mapsto (x_1,\ldots,x_m,t-x_1)
\]
with
\[
B=\{(x,t)\in H^{m+1}\mid  t\leq k_0+(k_1-1)x_1+\sum_{i=2}^m k_i
x_i\}
\]
and the theorem  for $B$ follows by induction on $k_1$. We
conclude the proof by the following linear bijection:
\[
A_1\to H^{m+1}:(x,t)\mapsto(x_1-1-t,x_2,\ldots,x_m,t).
\]
 \end{proof}
 \begin{theorem}[Parametric Rectilinearisation]\label{param recti}
Let $X\subset G^{m+n}$ be a $\emptyset$-definable Presburger set,
then there exists a finite partition $\Pm$ of $X$ into
$\emptyset$-definable Presburger sets, such that for each $A\in
\Pm$ there is a set $B\subset G^{m+n}$ with $\pi_m(A)=\pi_m(B)$
and a $\emptyset$-definable family $\{f_\lambda\}_{\lambda\in
\pi_m(A)}$ of linear bijections $f_\lambda:A_\lambda\subset G^n\to
B_\lambda\subset G^n$ with $B_\lambda$ a set of the form
$H^l\times\Lambda_\lambda$ where $\Lambda_\lambda$ is a bounded
$\lambda$-definable set and the integer $l$ only depends on
$A\in\Pm$.
 \end{theorem}
 \begin{proof} We give a proof by induction on $n$, following the lines
of the proof of Theorem \ref{recti}. So we assume that $X$ is a
cell
\[
\begin{array}{l}
X =  \{(\lambda,x,t)\in G^{m+(n+1)}\mid
 \\
 (\lambda,x)\in D,\
\alpha(\lambda,x)\sq_1 t\sq_2 \beta(\lambda,x),\ t\equiv
c\pmod{n}\},
\end{array}
\]
with $\lambda=(\lambda_1,\ldots,\lambda_m)$, $x=(x_1,\ldots,x_n)$,
$D\subset G^{m+n}$ a cell, integers $0\leq c<n$,
$\alpha,\beta:D\to G$ $\emptyset$-definable linear functions and
$\sq_i$ either $\leq$ or no condition. By subsequently applying
the induction hypothesis to $D$, partitioning further and applying
linear bijections (similar as to obtain Eq.~(\ref{some cell}) in
the proof of Theorem \ref{recti}, keeping the parameters $\lambda$
fixed now), we may assume that $X$ has the form
 \[
X=\{(\lambda,x,t)\in G^{m+n+1}\mid (\lambda,x)\in D',\ 0\leq t
\leq \gamma(\lambda,x)\},
 \]
with $\pi_{m+n}(X)=D'\subset G^{m+n}$ a Presburger set such that
for each $\lambda\in\pi_m(D')$ $D'_\lambda=H^l\times
\Gamma_\lambda$ where $\Gamma_\lambda$ is a $\lambda$-definable
bounded set, $l$ a fixed positive integer and $\gamma:D'\to G$ a
$\emptyset$-definable linear function. If $l=0$, $X_\lambda$ is a
bounded set for each $\lambda$ and the theorem follows
immediately. Let thus $l\geq 1$, i.e.~the projection of $X$ on the
$x_1$-coordinate is $H$, then the function $\gamma$ can be written
as $(\lambda,x)\mapsto k_1x_1+\gamma'(\lambda,x_2,\ldots, x_m)$
with $k_1$ an integer, necessarily nonnegative because the
projection of $X$ on the $x_1$-coordinate is $H$ and $\gamma'$ is
a linear function. The reader can finish the proof by induction on
$k_1\geq 0$, similar as in the proof of Theorem \ref{recti}.
 \end{proof}
 \begin{theorem}[Classification]\label{classification2}
Let $X$ be a $\emptyset$-definable Presburger set with
$\mathrm{dim}X=m>0$, then there exists a $\emptyset$-definable
Presburger bijection $f:X\to G^m$. In other words, there exists a
$\emptyset$-definable Presburger bijection between two infinite
$\emptyset$-definable Presburger sets $X,Y$ if and only if
$\mathrm{dim} X=\dim Y$.
 \end{theorem}
 \begin{proof} Let $X$ be $\emptyset$-definable and infinite.  We use induction on
$\mathrm{dim}X=m$. We say for short that two Presburger sets $X,Y$
are isomorphic if there exists a $\emptyset$-definable Presburger
bijection between them and write $X\cong Y$. If $m=1$, then
Theorem \ref{recti} yields a partition $\Pm$ of $X$ such that each
part is either a point or isomorphic to $H$. Consider the
bijections
\[
 \begin{array}{l}
f_1:H\to G :\left\{\begin{array}{rcl} 2x & \mapsto & x,
 \\ 2x+1 & \mapsto& -x,\end{array}\right.\\

f_2:H\cup\{-1\}\to H:x\mapsto x+1,\\

f_3:(\{0\}\times H) \cup (\{1\}\times H)\to
H:\left\{\begin{array}{rcl} (0,x) & \mapsto & 2x,\\ (1,x) &
\mapsto& 2x+1;\end{array}\right.
 \end{array}\]
the bijections $f_1,f_2$, applied repeatedly to (isomorphic copies
of) parts in $\Pm$ yield a definable bijection from $X$ onto $H$
and thus $G\cong X$ by applying $f_1$ (in the obvious way). Now
let $\mathrm{dim}X=m>1$. Using Theorem \ref{recti} we find a
partition $\Pm$ of $X$ such that each part is isomorphic to $H^l$
and thus to $G^l$ since $H\cong G$ by $f_1$. Since
$\mathrm{dim}X=m$, at least one part is isomorphic to $G^m$. Take
$A,B\in\Pm$ with $A\cong G^m$ and $B\cong G^l$, then it suffices
to prove that $A\cup B\cong G^m$. If $l=0$ this is clear and if
$l>0$ then $A\cup B\cong G\times(A'\cup B')$ for some disjoint and
$\emptyset$-definable sets $A',B'$ with $A'\cong G^{m-1}$ and
$B'\cong G^{l-1}$. The induction hypothesis applied to $A'\cup B'$
finishes the proof.
 \end{proof}
 \section{Expansions of $Z$-groups}\label{section Michaux}
We define the notion of Presburger minimality ($\Lp$-minimality)
for expansions of Presburger structures $(G,\Lp)$. This notion of
$\Lp$-minimality is a concrete instance of the general notion of
$\Lm$-minimality as in \cite{Macpherson} and has already been
studied in \cite{Point}.
 \begin{definition}
Let $G$ be a $Z$-group and $\Lm$ an expansion of the language
$\Lp$, then we say that $(G,\Lm)$ is $\Lp$-minimal if every
$\Lm$-definable subset of $G$ is already  $\Lp$-definable
(allowing parameters as always). We say that $\mathrm{Th}(G,\Lm)$
is $\Lp$-minimal if every model of this theory is $\Lp$-minimal.
 \end{definition}
Comparing this notion with the terminology of \cite{Point}, a
structure $(G,\Lm)$ is $\Lp$-minimal if and only if it is a
\emph{discrete coset-minimal group without definable proper convex
subgroups} (see \cite{Point}). Theorem 4.8 of \cite{Point} says
that a definable function in one variable between such groups is
piecewise linear. We reformulate this result with our terminology.
 \begin{proposition}[\cite{Point}, Thm.~4.8]\label{piecewise linear}
Let $(G,\Lm)$ be $\Lp$-minimal, then any definable function
$f:G\to G$ is piecewise linear.
 \end{proposition}
Proposition \ref{piecewise linear} allows us to repeat without any
change the compactness argument of the proof of the cell
decomposition theorem for any model of a $\Lp$-minimal theory.
This leads to the following description of $\Lp$-minimal theories.
 \begin{theorem}\label{strong}
Let $(G,\Lm)$ be an expansion of a Presburger structure $(G,\Lp)$,
then the following are equivalent:
\begin{itemize}
 \item[(i)] $\mathrm{Th}(G,\Lm)$ is $\Lp$-minimal;
 \item[(ii)] $(G,\Lm)$ is a definitional expansion of $(G,\Lp)$; precisely,
 any $\Lm$-definable set $X\subset G^m$ is already $\Lp$-definable.
 \end{itemize}
Thus, the theory $\mathrm{Th}(G,\Lp)$ does not admit any proper
$\Lp$-minimal expansion.
 \end{theorem}
 \begin{proof} Any Presburger minimal theory has definable Skolem
functions. For if $X\subset G^{m+1}$ is a definable set in some
model $G$, we can choose definably for any $x\in\pi_m(X)$  the
smallest nonnegative element in $X_x$ if there is any, and the
largest  negative element otherwise (this is well-defined by
Presburger minimality). This implies the definability of Skolem
functions by induction. Now replace in the statement of the cell
decomposition Theorem (theorem \ref{cell decomp}) the word
\emph{$\Lp$-definable} by \emph{$\Lm$-definable}. Then repeat the
case $m=1$ of the proof of Theorem \ref{cell decomp}, using now
the $\Lp$-minimality and Proposition \ref{piecewise linear}. Using
the same compactness argument as in the proof of Theorem \ref{cell
decomp} we find that any $\Lm$-definable set $X\subset G^m$ is a
finite union of Presburger cells, thus a fortiori, $X$ is
$\Lp$-definable.
 \end{proof}
 \begin{remark}
For an arbitrary expansion $(G,\Lm)$ of $(G,\Lp)$ it is, as far as
I know, an open problem whether the statements (i) and (ii) of
Thm.~\ref{strong} are equivalent with the following:
\begin{itemize}
\item[(iii)] $(G,\Lm)$ is $\Lp$-minimal.
\end{itemize}
In the special case  $G=\Z$, statements (i), (ii) and (iii) are
indeed equivalent, proven by C.~Michaux and R.~Villemaire
in~\cite{Michaux}.
 \end{remark}
 \section{Application to p-minimal fields}\label{section p}
In this section, we let $K$ be a $p$-adically closed field with
value group $G$. Recall that a $p$-adically closed field is a
field $K$ which is elementary equivalent to a finite field
extension of the field $\Q_p$ of $p$-adic numbers; in particular,
the value group $G$ is a $Z$-group and $K$ has quantifier
elimination in the Macintyre language
$\Lmac=\langle+,-,.,0,1,\{P_n\}_{n\geq 1}\rangle$ where $P_n$
denotes the set of $n$-th powers in $K^\times$. We write $v:K\to
G\cup\{\infty\}$ for the valuation map and for any $m>0$ we write
$\bar v$ for the map $\bar v:(K^\times)^ m\to G^m:x\mapsto
(v(x_1),\ldots,v(x_m))$. We give a definition of p-minimality,
extending the original definition of \cite{Haskell} slightly.
 \begin{definition}
Let $K$ be a $p$-adically closed field and let $(K,\Lm)$ be an
expansion of $(K,\Lmac)$. We say that the structure $(K,\Lm)$ is
p-minimal if any $\Lm$-definable subset of $K$ is already
$\Lmac$-definable (allowing parameters). The theory
$\mathrm{Th}(K,\Lm)$ is called p-minimal if every model of this
theory is p-minimal.
 \end{definition}
Examples of p-minimal fields known at this moment are $p$-adically
closed fields with the semialgebraic structure and with
subanalytic structure with restricted power series (see
\cite{vdDHM} and chapter \ref{chap:cell}). Theorem
\ref{application} exhibits a close connection between definable
sets in a p-minimal field $K$ and Presburger sets in the value
group $G$ of $K$; to prove it, we use Lemma \ref{Mac}, which is a
reformulation of the interpretability of $(G,\Lp)$ in $(K,\Lmac)$.
 \begin{lemma}\label{Mac}
Let $K$ be a $p$-adically closed field with value group $G$, then
for any $\Lp$-definable set $S\subset G^m$ the set
\[
\bar v^{-1}(S)=\{(x_1,\ldots,x_m)\in (K^\times)^m\mid\bar v(x)\in
S\}
\]
is $\Lmac$-definable.
 \end{lemma}
 \begin{proof} Let $S\subset G^m$ be $\Lp$-definable. By Theorem \ref{cell
decomp} we may suppose that $S$ is a Presburger cell. The Lemma
follows now inductively  from the fact that conditions imposed on
$(x_1,\ldots,x_{m-1},t)\in (K^\times)^m$ of the form $\pm
v(t)\leq\frac{1}{e}(\sum_{i=1}^{m-1}a_iv(x_i))+d$ or $v(t)\equiv
c\pmod{n}$ are $\Lmac$-definable for any integers $a_i$,
$e\not=0$, $0\leq c< n$ and $d\in G$ (see e.g.~\cite[Lemma
2.1]{Denef2}).
 \end{proof}
 \begin{theorem}\label{application}
Let $(K,\Lm)$ be a p-minimal field with p-minimal theory and let
$G$ be the value group of $K$. Then for any $\Lm$-definable set
$X\subset (K^\times)^m$ the set
\[
\bar v(X)=\{(v(x_1),\ldots,v(x_m))\in G^m\mid (x_1,\ldots,x_m)\in
X\} \subset G^m
\]
is $\Lp$-definable.
 \end{theorem}
 \begin{proof} Put $S_m=\{\bar v(X)\subset G^m\mid X\subset (K^\times)^m, \mbox{ $X$
is $\Lm$-definable}\}$, then it is easy to see that the collection
$(S_m)_{m\geq 0}$  determines  a  structure on $G$ (i.e.~the
collection $\cup_m S_m$ is precisely the collection of
$\Lm'$-definable sets for some language $\Lm'$). We first show
that this structure is in fact $\Lp$-minimal. Choose a
$\Lm$-definable set $X\subset K^\times$, then, by p-minimality,
$X$ is $\Lmac$-definable. We can thus apply the $p$-adic
semialgebraic cell decomposition (\cite{Denef2}, in the
formulation of \cite[Lemma 4]{C}, see also chapter
\ref{chap:cell})
 to the set $X$ to obtain that $X$ is a
finite union of $p$-adic cells, i.e.~sets of the form
\[
\{x\in K\mid  v(a_1)\sq_1 v(x-c)\sq_2 v(a_2),\ x-c\in\lambda
P_n\}\subset K^\times,
\]
with $a_1,a_2,c,\lambda\in K$ and $\sq_i$ either $\leq,<$ or no
condition. The image under $v$ of such a cell is either a finite
union of $(0$)-cells or a $(1)$-cell and thus a $\Lp$-definable
subset of $G$. By consequence,  the structure $(S_m)_{m\geq_0}$ is
$\Lp$-minimal. By the Presburger minimality of $(S_m)_{m\geq 0}$,
the p-minimality of $\mathrm{Th}(K,\Lm)$, and Lemma \ref{Mac} to
interpret $G$ into $K$, we can repeat the compactness argument of
the proof of the cell decomposition theorem \ref{cell decomp} for
the structure $(S_m)_m$ on $G$ to find that each $A\in \cup_m S_m$
is a finite union of Presburger cells. This proves the theorem.
 \end{proof}
 \begin{remark}
For a $p$-adically closed field $K$ it is proven by Cluckers in
\cite{C} (chapter \ref{chap:classification}) that there exists a
$\Lmac$-definable bijection $X\to Y$ between two infinite
parameter free $\Lmac$-definable sets $X,Y$ if and only if $\dim
X=\dim Y$. This is proven by reducing to a Presburger problem
similar to the classification theorem in section
\ref{sectionclassification}. An analogous classification for
$\emptyset$-definable sets in arbitrary p-minimal structures is
not known up to now.
 \end{remark}
 \section{Elimination of imaginaries}
As a last application of the cell decomposition theorem we prove
uniform elimination of imaginaries for Presburger structures. We
say that a structure $(M,\Lm)$ has uniform elimination of
imaginaries if for any $\emptyset$-definable equivalence relation
on $M^k$ there exists a $\emptyset$-definable function $F:M^k\to
M^r$ for some $r$ such that two tuples $x,y\in M^k$ are equivalent
if and only if $F(x)=F(y)$.
 \begin{theorem}
 The theory $Th(\Z,\Lp)$ has uniform elimination of imaginaries,
meaning that any Presburger structure $(G,\Lp)$ eliminates
imaginaries uniformly.
 \end{theorem}
 \begin{proof} Since $\mathrm{Th}(\Z,\Lp)$ has definable Skolem functions,
we only have to prove the following statement for an arbitrary
$Z$-group $G$ (see e.g.~\cite[Lemma 4.4.3]{Hodges}). For any
$\emptyset$-definable Presburger set $X\subset G^{m+1}$ there
exists a $\emptyset$-definable Presburger function $F:G^m\to G^n$
for some $n$, such that $F(x)=F(x')$ if and only if $X_x=X_{x'}$
(if $x\not\in \pi_m(X)$ then we put $X_x=\emptyset$). So let
$X\subset G^{m+1}$ be a $\emptyset$-definable Presburger set.
Apply the cell decomposition theorem to obtain a partition $\Pm$
of $X$ into cells. For each cell $A\in \Pm$ of the form
$A=\{(x,t)\in G^{m+1}\mid x\in D,\ \alpha(x)\sq_{1A} t\sq_{2A}
\beta(x),\ t\equiv c\pmod{n}\}$ (as in Eq.~\ref{cell}) and each
$\xi=(\xi_1,\xi_2)\in G^2$ we define a set
\[
C_A(\xi)=\{t\in G\mid \xi_1\sq_{1A} t\sq_{2A} \xi_2,\ t\equiv
c\pmod{n}\}.
\]
Notice that for each $x\in \pi_m(X)$ we have at least one
partition of $X_x$ into sets of the form $C_A(\xi)$ with $A\in\Pm$
and $\xi\in G^2$. For $x,y\in G$ we write $x\lhd y$ if and only if
one of the following conditions is satisfied
 \begin{itemize}
\item[(i)] $0\leq x< y$,\\
\item[(ii)]$0< x\leq -y$,\\
\item[(iii)]$0<-x<y$,\\
\item[(iv)]$0<-x<-y$.
 \end{itemize}
This gives a new ordering $0 \lhd 1\lhd -1\lhd 2\lhd -2\lhd\ldots$
on $G$  with zero as its smallest element. For each $k>0$ we also
write $\lhd$ for the lexicographical order on $G^k$ built up with
$\lhd$. The order $\lhd$ is $\Lp$-definable and each Presburger
set has a unique $\lhd$-smallest element. For each $x\in G^m$ and
each $I\subset \Pm$ with cardinality $|I|=s\geq0$ we let
$y_I(x)=(\xi_A)_{A\in I}$, $\xi_A\in G^2$, be the $\lhd$-smallest
tuple in $G^{2s}$ such that $\cup_{A\in I}C_A(\xi_A)=X_x$ if there
exists at least one such tuple and we put $y_I(x)=(0,\ldots,0)\in
G^{2s}$ otherwise. One can reconstruct the set $X_x$ given all
tuples $y_I(x)$, $I\subset \Pm$. Let $F$ be the function mapping
$x\in \pi_m(X)$ to $y=(y_I(x))_{I\subset \Pm}$. Since the
lexicographical order $\lhd$ is $\Lp$-definable it is clear that
$F$ is $\Lp$-definable and that $F(x)=F(x')$ if and only if
$X_x=X_{x'}$ for each $x,x'\in G^m$.
 \end{proof}

\subsection*{Acknowledgment}
I would like to thank J.~Denef, D.~Haskell, C.~Michaux, F.~Point,
F.~ Wagner and K.~Zahidi for fruitful discussions during the
preparation of this chapter.

\clearemptydoublepage

\chapter{Grothendieck rings of $\Z$-valued fields}\label{chap:CH}
\begin{abstract}\footnote{This chapter
corresponds to \cite{CH}.}
 We prove the triviality of the
Grothendieck ring of a $\Z$-valued field $K$ under slight
conditions on the logical language and on $K$. We construct a
definable bijection from the plane $K^2$ to itself minus a point.
When we specialize to local fields with finite residue field, we
construct a definable bijection from the valuation ring to itself
minus a point.
\end{abstract}
\section{Introduction}
At the Edinburgh meeting on the model theory of valued fields in
May 1999, Luc B\'elair posed the question of whether there is a
definable bijection between the set of $p$-adic integers and the
set of $p$-adic integers with one point removed. At the same
meeting, Jan Denef asked what is the Grothendieck ring of the
$p$-adic numbers, as did Jan Kraj\'{\i}\v{c}ek independently in
\cite{K}. A general introduction to Grothendieck rings of logical
structures was recently given in \cite{KS} and in \cite{DL}, par.
3.7. Calculations of non-trivial Grothendieck rings and related
topics such as motivic integration can be found in \cite{DLinvent}
and \cite{DL}. The logical notion of the Grothendieck ring of a
structure is analogous to that of the Grothendieck ring in the
context of algebraic $K$-theory and has analogous elementary
properties (see \cite{S}). Here we recall the definition.

\begin{definition} Let $\M$ be a structure and $\Def(\M)$ the set
of definable subsets of $M^n$ for every positive integer $n$. For any
$X,Y\in\Def(\M)$, write $X\cong Y$ iff there is a definable bijection (an
isomorphism) from $X$ to $Y$. Let $F$ be the free abelian group whose
generators are isomorphism classes $\lfloor X\rfloor$ with $X\in\Def(\M)$
(so $\lfloor X\rfloor =\lfloor Y\rfloor$ if and only if $X\cong Y$) and
let $E$ be the subgroup generated by all expressions $\lfloor X\rfloor
+\lfloor Y\rfloor -\lfloor X\cup Y\rfloor -\lfloor X\cap Y\rfloor$ with
$X,Y\in\Def(\M)$. Then the Grothendieck group of $\M$ is the quotient
group $F/E$. Write $[X]$ for the image of $X\in\Def(\M)$ in $F/E$. The
Grothendieck group has a natural structure as a ring with multiplication
induced by $[X]\cdot[Y]=[X\times Y]$ for $X,Y\in\Def(\M)$. We call this
ring the Grothendieck ring $K_0(\M)$ of $\M$.
\end{definition}
\par
It is easy to see that the above questions are related: the
Gro\-then\-dieck ring is trivial if and only if there is a
definable bijection between $M^k$ and itself minus a point for
some $k$, which happens if and only if the Grothendieck group is
trivial. Moreover, if we find such a $k$ then we have for any
$X\in\Def(\M)$ a definable bijection from the disjoint union of
$M^k\times X$ and $X$ to $M^k\times X$; if there is a definable
injection from $M^k$ into $X$ we find a definable bijection from
$X$ to itself minus a point.
\par
In this paper we answer the questions posed by B\'elair and Denef.
Furthermore, we prove the triviality of the Grothendieck ring of any
$\Z$-valued field which satisfies some slight conditions and give in this
general setting an explicit bijection from the plane to itself minus a
point. For the fields $\Q_p$ and $\F_q((t))$ we explicitly construct a
definable bijection from the valuation ring to itself minus a point.
\par
Dave Marker independently produced a definable bijection from
$\Z_p$ to $\Z_p\setminus\{0\}$, after it was noticed by Lou~van
den Dries that its existence followed from unpublished notes of
D.~Haskell. R.~Cluckers has proved further that there is a
definable bijection between any two definable sets in the
$p$-adics if and only if they have the same dimension, see
\cite{C} and chapter \ref{chap:classification}.
\subsection*{Notation and terminology}
Fix a $\Z$-valued field $K$, that is,  a field with a valuation
$v:K^\times\to Z$ to an ordered group $Z$ which is elementary
equivalent to the integers in the Presburger language. Let
$R=\{x\in K|v(x)\geq0\}$ be the valuation ring,
$R^*=R\setminus\{0\}$ and $\bar K=R/m$ the residue field, with $m$
the maximal ideal of $R$ and natural projection $R\to \bar
K:x\to\bar x$. An angular component map is a homomorphism
$ac:K^\times\to\bar K^\times$ such that $ac(x)=\bar x$ if
$v(x)=0$. We extend $ac$ to a map $ac:K\to\bar K$ by putting
$ac(0)=0$ (for the existence of angular component maps, see
\cite{P} and \cite{Belair}).
\begin{definition} Let $\Lm$ be an extension of the
language of rings with $K$ as a model. We say that the structure
$(K,\Lm)$ satisfies condition ($*$) if we can choose an angular component
map $ac$ and an $\Lm$-definable element $\pi\in R$ with $v(\pi)=1$ and
$ac(\pi)=1$ such that the sets $R$ and $R^{(1)}=\{x\in R|ac(x)=1\}$ are
$\Lm$-definable.
\end{definition}
Notice that if condition ($*$) is satisfied, the set $\{(x,y)\in
K^2|v(x)\leq v(y)\}$ is $\Lm$-definable by the formula $\exists
z\in R\,(zx=y)$. A bijection $X\to Y$ with $X,Y\in\Def{(K,\Lm)}$
with $\Lm$-definable graph will be called an isomorphism.
\par
Let $X\subset K^m$ and $Y\subset K^n$ be definable sets, $m\geq n$. Let
$X'=\{0\}\times X$ and $Y'=\{1\}^{m-n+1}\times Y$. Then we define the
disjoint union $X\sqcup Y$ of $X$ and $Y$ up to isomorphism to be $X'\cup
Y'$.  We say that a set $W$ is isomorphic to $X\sqcup Y$ if $W$ is
isomorphic to $X'\cup Y'$ and then obviously $[W]=[X]+[Y]$.  If $(K,\Lm)$
satisfies condition ($*$) then we can find $W\subset R^m$ with $W\cong
X\sqcup Y$ as follows. The map $i:K\to R$ which sends $x$ to $\pi x$ if
$v(x)\geq0$ and to $1+1/x$ if $v(x)<0$ is a definable injection. For
$m=n=1$, put $X''=\pi.i(X)$ and $Y''=1+\pi.i(Y)$. Then $X''\cong X$,
$Y''\cong Y$ and $X''\cap Y''=\phi$, so  $W=X''\cup Y''$ is isomorphic to
$X\sqcup Y$. For $m>1$, use the same method in each coordinate.
\section{Calculations of Grothendieck rings}
\begin{proposition}\label{bijections} Let $K$ be a $\Z$-valued field, which
is a model for the language $\Lm$. If the structure $(K,\Lm)$ satisfies
condition ($*$), then the following holds:
\item{(i)} The disjoint union of $R$ and $R^{(1)}$ is isomorphic to
$R^{(1)}$ and thus $[R]=0.$
\item{(ii)} The disjoint union of two copies of $R^{* 2}$ is isomorphic to
$R^{* 2}$ itself, and hence $[R^{* 2}]=0.$
\end{proposition}
\begin{proof}
(i) The map
\[
 \{0\}\times R\cup \{1\}\times R^{(1)}\to R^{(1)}:\left\{\begin{array}{rcl}
              (0,x) & \mapsto & 1+\pi x,\\
              (1,x) & \mapsto & \pi x,\end{array}\right.\
\]
is easily seen to be an isomorphism as required. This yields in the
Grothendieck ring $[R] + [R^{(1)}]= [R^{(1)}]$, so $[R]=0$.

(ii)
 Define the sets
\begin{eqnarray*}
X_1=\{(x,y)\in R^{* 2}|v(x)\leq v(y)\},\\
X_2=\{(x,y)\in R^{* 2}|v(x)>v(y)\},
\end{eqnarray*}
then $X_1,X_2$ form a partition of $R^{* 2}$. The isomorphisms
\begin{eqnarray*}
\{0\}\times R^{* 2}\to X_1: (0,x,y)\mapsto (x,xy),\\
\{1\}\times R^{* 2}\to X_2: (1,x,y)\mapsto (\pi xy,y),
\end{eqnarray*}
imply that $R^{* 2}\sqcup R^{* 2}$ is isomorphic to $X_1\cup X_2=R^{* 2}$.
It follows that $2[R^{* 2}]=[R^{* 2}]$, so $[R^{* 2}]=0$. Notice that the
proof of (ii) does not use the full power of ($*$), only that $R$ is
definable.
\end{proof}
\begin{theorem}\label{Gring=0}
Let $K$ be a $\Z$-valued field, which is a model for the language $\Lm$.
If the structure $(K,\Lm)$ satisfies condition ($*$), then the
Grothen\-dieck ring $K_0(K)$ is trivial and there exists an isomorphism
from $R^2\setminus\{(0,0)\}$ to $R^2$.
\end{theorem}
\begin{proof} Since $0=[R]=[R^*]+[\{0\}]$ we have $[R^*]=-1$. Together with
$0=[R^{* 2}]=[R^*]^2$ this yields $1=0$, so $K_0(K)$ is trivial.
\par
Define the isomorphisms $\psi:R^2\to \pi^3 R^2:(x,y)\mapsto (\pi^3
x,\pi^3 y)$ and $\varphi_i:R^2\to (\pi^i+\pi^3 R)\times(\pi^i+\pi^3
R):(x,y)\mapsto (\pi^i+\pi^3 x,\pi^i+\pi^3 y)$ for $i=1,2$.
\par
Since clearly $\psi(R^*\times R^*)\cup \varphi_1(R^*\times R^*)$ is
isomorphic to $R^{*2}\sqcup R^{*2}$, we can find by
Proposition~\ref{bijections}(ii) an isomorphism
\[
f_1:\varphi_1(R^*\times R^*)\to \psi(R^*\times R^*)\cup
\varphi_1(R^*\times R^*).\] Define $f_2$ by
\[
f_2:\psi(R\times R^*)\cup\varphi_2(R^{(1)}\times R^*) \to
\varphi_2(R^{(1)}\times R^*):
\]
\[
 \left\{\begin{array}{rcl}
              \psi(x,y) & \mapsto & \varphi_2(1+\pi x,y),\\
              \varphi_2(x,y) & \mapsto & \varphi_2(\pi x,y).\end{array}\right.\
\]
Analogously, we can modify the function given in the proof of
Proposition~\ref{bijections}(i) to get an isomorphism
\[
f_3:\varphi_2(\{0\}\times R^{(1)})\to\varphi_2 (\{0\}\times
R^{(1)})\cup\psi(\{0\}\times R).
\]
Finally,
\[
  g:R^2\setminus\{(0,0)\}\to R^2:x\mapsto\left\{
\begin{array}{ll}
 f_1(x) & \mbox{if }x\in \varphi_1(R^*\times R^*),\\
 f_2(x) & \mbox{if }x\in \psi(R\times R^*)\cup \\
        & \qquad  \varphi_2(R^{(1)}\times R^*),\\
 f_3(x) & \mbox{if }x\in \varphi_2(\{0\}\times R^{(1)}),\\
 x      & \mbox{else},\\
\end{array}\right.
\]
is the required isomorphism.
\end{proof}
We give some examples for the conditions of Theorem~\ref{Gring=0} to be
satisfied. Let $\Lm_{ac}$ be the language of rings with an extra constant
symbol to denote $\pi$ and a relation symbol to denote the set $R^{(1)}$.
Let $\Lm_{ac,R}$ be the language $\Lm_{ac}$ with an extra relation symbol
to denote $R$.
\begin{itemize}
\item Let $K$ be a valued field with valuation to the integers $\Z$. Then we can define an
angular component as follows. Choose $\pi\in K$ with $v(\pi)=1$ and put
$ac(x)=\overline{\pi^{-v(x)}x}$ for $x\not=0$. Then clearly $ac(\pi)=1$
and $(K,\Lm_{ac,R})$ satisfies condition ($*$).
\item Let $K$ be a Henselian field with valuation to the integers $\Z$.
Then $R$ is already definable in the language of rings: \\
if char$(\bar{K})\not=2$ we have $R=\{x\in K|\exists y\in
K,y^2=1+\pi x^2\}$ and if char$(\bar{K})=2$ then we use the
formula $\exists y\in K,y^3=1+\pi x^3$ to define $R$. This implies
that $(K,\Lm_{ac})$ satisfies condition ($*$).
\item For definability of
the valuation ring in fields of rational functions within the
language of rings, see \cite{D} and \cite{KR}.
\end{itemize}
\par
Now we specialize our attention to local fields with finite residue field.
\begin{theorem}\label{FF}
Let $K=\F_q((t))$ be the formal Laurent series over the finite field
$\F_q$ and $\Lm_t$ the language of rings with a constant symbol to denote
$t$. Then $K_0(K)$ is trivial and we have an isomorphism $R\to R^*$.
\end{theorem}
\begin{proof}
We first show that $K$ satisfies condition ($*$). Since $K$ is a Henselian
field, $R$ is definable as shown above. For each $x\in\F_q$ we have
$x^{q-1}=1$, so we can define $R^{(1)}$ as
\[R^{(1)}=\{x\in R|\exists y\in R^*,\ \bigvee_{n=0}^{q-2} t^ny^{q-1}=x\},\]
 again by Hensel's lemma.
\par
By Theorem~\ref{bijections} we have an isomorphism $f:R^2\to
R^2\setminus\{(0,0)\}$. For a Laurent series $H(t)\in K$ we have
$H(t)^p=H(t^p)$. Consequently, the map
\[g:K^2\to K:(x,y)\mapsto x^p+ty^p\]
is an injection from the plane into the line. We obtain the isomorphism
\[R\to R^*:x\mapsto \left\{
\begin{array}{ll}
 g\circ f\circ g^{-1}(x) & \mbox{if } x\in g(R^2),\\
 x & \mbox{else.}
\end{array}\right.
\]
\end{proof}
\par
Now let $\Q_p$ be the field of $p$-adic numbers and $K$ a fixed finite
field extension of $\Q_p$. Choose an element $\pi$ with $v(\pi)=1$, then
$ac(x)=\pi^{-v(x)}x\ \mathrm{ mod }(\pi)$  defines an angular component for
$x\not=0$. We work with $\Lm_\pi$, the language of rings with an extra
constant symbol to denote $\pi$. For a definable set $X\subset K$ and
$k\in\N_0$ we write
\begin{eqnarray*}
X^{(k)} & = & \{x\in X|x\not=0 \mbox{ and } v(\pi^{-v(x)}x-1)\ge k\},
\end{eqnarray*}
which corresponds with our previous definition of $R^{(1)}$. The set $R$
and each $X^{(k)}$ is definable by the same argument as in the proof of
Theorem \ref{FF}, so $(K,\Lm_\pi)$ satisfies condition $(*)$. We put
$P_n=\{x\in K^\times|\exists y\in K,\ y^n=x\}$ and $\bar P_n=P_n\cap R$.
Recall that $P_n$ is a subgroup of finite index in $K^\times$ for each
$n$.
\par
For convenience, we recall the following easy corollary of Hensel's Lemma.
\begin{cor}\label{corhensel2}
Let $n>1$ be a natural number. For each $k>v(n)$, and $k'=k+ v(n)$ the
function
$$
K^{(k)}  \to P_n^{(k')}:x\mapsto x^n
$$
is an isomorphism.
\end{cor}
\par
In the next proposition we exhibit some isomorphisms between definable
sets.
\begin{proposition}\label{p-adic} Let $K$ be a finite field extension of
the $p$-adic numbers and $\Lm_\pi$ the language of rings with an extra
constant symbol to denote $\pi$. Then we have
\item{(i)} for each $k>0$, the union of two disjoint
copies of $R^{(k)}$ is isomorphic to $R^{(k)}$;
\item{(ii)} the union of two disjoint
copies of $R^*$ is isomorphic to $R^*$.
\end{proposition}
\begin{proof}
(i) \textbf{Case 1: $p\not=2$.} The map $R^{(k)}\to\bar
P_2^{(k)}:x\mapsto x^2$ is an isomorphism for each $k>0$ by
Corollary \ref{corhensel2}. By Hensel's Lemma, $R^{(k)}=\bar
P_2^{(k)}\cup\pi\bar P_2^{(k)}$ is a partition. Hence the function
\[
 \{0\}\times R^{(k)}\cup \{1\}\times R^{(k)}\to R^{(k)}:\left\{\begin{array}{rcl}
                    (0,x) & \mapsto & x^2,\\
                    (1,x) & \mapsto & \pi x^2,\end{array}\right.\
\]
is an isomorphism.
\par
\textbf{Case 2: $p=2$.} The map $R^{(k)}\to \bar
P_3^{(k)}:x\mapsto x^3$ is an isomorphism by Corollary
\ref{corhensel2}, and by Hensel's Lemma $R^{(k)}=\bar
P_3^{(k)}\cup\pi\bar P_3^{(k)}\cup\pi^2\bar P_3^{(k)}$ is a
partition. Explicitly, we see that cubing and multiplying by $1$,
$\pi$ or $\pi^2$ is an isomorphism from three disjoint copies of
$R^{(k)}$ to $R^{(k)}$.  First suppose that $k>v(2)$ and put
$k'=k+v(2)$, then $R^{(k)}\to\bar P_2^{(k')}:x\mapsto x^2$ is an
isomorphism  by Corollary \ref{corhensel2}. By Hensel's lemma, we
have a partition $R^{(k)}=\bigcup_{i=1}^{2^l}\alpha_i\bar
P_2^{(k')}$ for some $l\in\N_0$. Thus we can say there are
isomorphisms from $R^{(k)}$ to $2^l$ disjoint copies of $R^{(k)}$
and to three disjoint copies of $R^{(k)}$. Some arithmetic on the
number of disjoint copies yields the required isomorphism for
$k>v(2)$.
\par
 If $k\leq v(2)$ then $R^{(k)}$ admits a finite partition into parts of the form
$\alpha R^{(v(2)+1)}$, with $v(\alpha)=0$, and hence that the required
isomorphism exists follows from property (i) for $R^{(v(2)+1)}$.
\par
(ii) Since  $R^*$ admits a finite partition with parts of the form
$\alpha R^{(1)}$ with $v(\alpha)=0$, this follows from (i).
\end{proof}
Now we give the solution of the problems raised by J. Denef and L.
B\'elair.
\begin{theorem}
Let $K$ be a finite field extension of $\Q_p$ and $\Lm_\pi$ the language
of rings with an extra constant symbol to denote $\pi$. Then $K_0(K)=0$
and we have an isomorphism from $R$ to itself minus a point.
\end{theorem}
\begin{proof}
The triviality of the Grothendieck ring follows from
Theorem~\ref{Gring=0}.
\par
 We write the isomorphism explicitly in the case
$p\ne 2$. First let
\[ W = 1+ \pi^2R^* \cup \pi^2R \cup \pi + \pi^2R^{(1)}.
\]
As in the proof of Proposition~\ref{p-adic}, we can write
\[ R^* =\bigcup_{i=1}^l \alpha_iR^{(1)} =
       \bigcup_{i=1}^l\big( \alpha_i{\bar P_2^{(1)}}
                             \cup \pi\alpha_i{\bar P_2^{(1)}}\big)
\]
as a partition for some $l\in \N_0$. Thus the function
\[
   f_1:\pi^2 R^* \cup 1+\pi^2 R^*\to 1+\pi^2 R^*:
                          \left\{\begin{array}{rcl}
                 \pi^2\alpha_i x & \mapsto & 1+\pi^2(\alpha_i x^2),\\
                 1+\pi^2\alpha_i x & \mapsto & 1+\pi^2(\pi\alpha_i x^2),
\end{array}\right.
\]
where $x\in R^{(1)}$, is a well-defined isomorphism. Modify the function
given in the proof of Proposition~\ref{bijections}(i), to get
\[
  f_2:\pi^2 R\cup \pi+\pi^2 R^{(1)}\to \pi+\pi^2 R^{(1)}:
                       \left\{\begin{array}{rcl}
               \pi^2 x & \mapsto & \pi+\pi^2(1+\pi x),\\
               \pi+\pi^2 x & \mapsto & \pi+\pi^2(\pi x).
\end{array}\right.
\]
Then the function
\[
  f:W\to W\setminus\{0\}:x\mapsto\left\{
                          \begin{array}{ll}
           f_1^{-1}(x) & \mbox{if }x\in 1+\pi^2 R^*,\\
           f_2(x) & \mbox{if }x\in\pi^2 R\cup \pi+\pi^2 R^{(1)},\\
                          \end{array}\right.
\]
is an isomorphism. Finally,
\[
  g:R\to R^*:x\mapsto\left\{
\begin{array}{ll}
 f(x) & \mbox{if }x\in W\\
 x      & \mbox{if } x\notin W\\
\end{array}\right.
\]
is an isomorphism.
\par
In the case $p=2$, we know from Proposition \ref{p-adic}(i) that there is
a function which plays the role of $f_1$. The rest is as above.
\end{proof}
\begin{remark}\label{remark:CPres:Euler}
\begin{itemize}
\item
The construction of the bijection $\Z_p\to\Z_p\setminus\{0\}$ also works
for the field $\F_q((t))$ if $2\nmid q$. If $2|q$ the proof of
Proposition~\ref{p-adic}(i) collapses since the index of the squares in
$\F_q((t))^\times$ is infinite.
\item The triviality of the Grothendieck ring of a structure $\M$
implies that every  Euler characteristic on the definable sets is trivial.
An Euler characteristic is a map $\chi:\Def(\M)\to R_\chi$ with  $R_\chi$
a ring, such that $\chi(X)=\chi(Y)$ if $X\cong Y$, $\chi(X\cup
Y)=\chi(X)+\chi(Y)$ if $X\cap Y=\phi$ and $\chi(X\times
Y)=\chi(X)\chi(Y)$. In general an  Euler characteristic on $\Def(\M)$
factorizes through $\Def(\M)\to K_0(\M):X\mapsto [X]$.
\end{itemize}
\end{remark}
%
%
%
%
%
%
%
%

\clearemptydoublepage

\chapter{Grothendieck rings of Laurent series fields}
\begin{abstract}
\footnote{This chapter corresponds to \cite{Cgroth}.} We study
Grothendieck rings (in the sense of logic) of fields. We prove the
triviality of the Grothendieck rings of certain fields by
constructing definable bijections which imply the triviality. More
precisely, we consider valued fields, for example, fields of
Laurent series over the real numbers, over $p$-adic numbers and
over finite fields, and construct definable bijections from the
line to the line minus one point.
\end{abstract}
 \section{Introduction}
Recently, the Grothendieck ring of a structure, in the sense of
logic, has been introduced in \cite{DL} and independently in
\cite{KS}. The Gro\-then\-dieck ring of a model-theoretical
structure is built up as a quotient of the definable sets by
definable bijections (see below), and thus, depends both on the
model and the language. For $(M,\cL)$ a structure with the
signature of a language $\cL$ we write $K_0(M,\cL)$ for the
Grothendieck ring of
 $(M,\cL)$. In \cite{CH} (cfr.~chapter \ref{chap:CH}) and \cite{vdD}, the following explicit
calculations of Grothendieck rings of fields are made:
\begin{quote}
$K_0(\R,\Lr)$ is isomorphic to $\Z$,\\
$K_0(\Q_p,\Lr)$ is trivial,\\
$K_0(\F_p((t)),\Lr)$ is trivial.
\end{quote}
Here, $\Lr$ is the language $(+,-,\cdot,0,1)$. In \cite{DL} and
\cite{KS} it is shown that the Grothendieck ring $K_0(\C,\Lr)$ is
extremely big and complicated; $K_0(\C,\Lr)$, and many other
Grothendieck rings, are not explicitly known.
\par Any Euler characteristic (in the sense of remark \ref{remark:CPres:Euler} of chapter
\ref{chap:CH}, see also \cite{CH} or \cite{KS}), defined on the
definable sets, factors through the natural projection of
definable sets into the Grothendieck ring, and, in this sense, to
know a Grothendieck ring is to know a universal Euler
characterictic. Nevertheless, it happens that a Grothendieck ring
is trivial.
\par
 The triviality of a
Grothendieck ring can be proven by  constructing a definable
bijection from $X$ to $X\setminus \{a\}$, where $X$ is a definable
set and $\{a\}$ a point on $X$. We develop general techniques to
obtain definable bijections $K\to K^\times$, where $K$ is a valued
field and $K^\times=K\setminus\{0\}$. In section
\ref{sect:laurent} we explain iterated Laurent series fields. In
the present paper we prove:
\begin{theorem}\label{thm:bijection:Qp}
Let $L$ be either $\Q_p$ or a finite field extension of $\Q_p$,
and let $K$ be one of the fields
\[
L,\ L((t_1)),\ L((t_1))((t_2)),\ L((t_1))((t_2))((t_3)),
\]
and so on. Then $K_0(K,\Lr)=0$ and there exists a $\Lr$-definable
bijection $K\to K^\times$.\footnote{Here, as always, definable
means definable with parameters. Given a model $M$ for some
language $\cL$, any formula obtained by replacing variables in a
$\cL$-formula by elements of $M$ is called a formula with
parameters from $M$. }\footnote{For $K=\Q_p$, this result is
proven in \cite{CH}, see chapter \ref{chap:CH}.}
\end{theorem}
\begin{theorem}\label{thm:bijection:Fq}
Let $L$ be $\F_q$ where $\F_q$ is the finite field with $q=p^l$
elements, $p$ a prime, and let $K$ be one of the fields
\[
L((t_1)),\ L((t_1))((t_2)),\ L((t_1))((t_2))((t_3)),
\]
and so on. Then $K_0(K,\Lr)=0$ and there exists a $\Lr$-definable
bijection $K\to K^\times$.\footnote{For $K=\F_q((t))$ this is
proven in \cite{CH}, see chapter \ref{chap:CH}.}
\end{theorem}
Central in the proofs of this paper is a subgroup $H(K,\cL)$ of
$\Z$, associated to a field $K$ and a language $\cL$, which is
sensitive to some elementary arithmetical properties of  the
indices of $n$-th powers in $K^\times$ and of the number of $n$-th
roots  of $1$ in $K^\times$ (see section \ref{sect:calcul}). Using
the definition, it follows immediately that, for example,
\begin{quote}
$H(\R,\Lr)$ is $\Z$,\\
$H(\Q_p,\Lr)$ is $\Z$, and\\
$H(\C,\Lr)$ is $\{0\}$.\\
\end{quote}
We give two criteria for valued fields, for which the value group
has a well-determined minimal strictly positive element, to have a
trivial Grothendieck ring (proposition \ref{prop:HK=Z} and
\ref{prop:crit2}).
\par
We also consider Laurent series fields over $\R$ and over fields
of characteristic $p>0$, using the language of Denef - Pas.
The language of Denef - Pas \cite{Pas} was introduced to study
uniform $p$-adic integrals for all primes $p$, and is now still
used in, for example, the theory of motivic integration (see
\cite{DL} and  \cite{DLinvent}).
\par
For any $\Z$-valued field $K$ with angular component map $\ac$,
the Grothendieck ring $K_0(K,\Lpas)$ is trivial, and there exists
a $\Lpas$-definable bijection from $K^2$ onto $K^2\setminus
\{(0,0)\}$, see \cite{CH}, Thm.~1 (cfr.~Thm.~\ref{Gring=0} of
chapter \ref{chap:CH}) and proposition \ref{prop:crit2} below.
(Proposition \ref{prop:crit2} is more general than \cite{CH},
Thm.~1, see also Thm.~\ref{Gring=0} of chapter \ref{chap:CH}). The
following theorems give stronger results for iterated Laurent
series fields.
\begin{theorem}\label{thm:bijection:R((t))}
Let $K$ be one of the fields $\R((t_1))$, $\R((t_1))((t_2))$,
$\R((t_1))((t_2))((t_3))$, and so on. We have
\[H(K,\Lr)=\Z.\]
Endow $K$ with a valuation onto a group of the form $\Z^k$ with
lexicographical order, $k>0$, and with the natural angular
component map (as in section \ref{sect:Pas}). Then
\[K_0(K,\Lpas)=\{0\},\]
and there exist a bijection $K\to K^\times$, definable in the
language $\Lpas$ of Denef - Pas.
\end{theorem}
\begin{theorem}\label{thm:bijection3}
Let $L$ be an arbitrary field of characteristic $p>0$. Let $K$ be
one of the fields $L((t_1))$, $L((t_1))((t_2))$,
$L((t_1))((t_2))((t_3))$ and so on. Endow $K$ with a valuation
onto a group of the form $\Z^k$ with lexicographical order, $k>0$,
and with the natural angular component map (as in section
\ref{sect:Pas}). Then $K_0(K,\Lpas)=\{0\}$ and there exists a
bijection $K\to K^\times$ definable in the language of Denef -
Pas.
\end{theorem}
\subsection{Valued fields}
Fix a field $K$. We call $K$ a valued field if there is an ordered
group\footnote{Here, an ordered group is a totally ordered
non-trivial abelian group $G$  such that $x<y$ implies $x+z<y+z$
for all $x,y,z$ in $G$.} $(G,+,\leq)$  and a surjective valuation
map $v:K\to G\cup\{\infty\}$ such that
\begin{itemize}
\item[(i)] $v(x)=\infty$ if and only if $x=0$;
\item[(ii)] $v(xy)=v(x)+v(y)$ for all $x,y\in K$;
\item[(iii)] $v(x+y)\geq \min\{v(x),v(y)\}$ for all $x,y\in K$.
\end{itemize}
We write $R$ for the  valuation ring $\{x\in K\mid v(x)\geq0\}$ of
$K$, $M$ for its unique maximal ideal and we write $k$ for the
residue field  $R/M$ and $p:R\to k$ for the natural projection. If
$G=\Z$ we call $K$ a $\Z$-valued field. Whenever $K$ is a
$\Z$-valued field, the valuation ring $R$ is a discrete valuation
ring, and a generator $\pi$ of the maximal ideal of $R$ is called
a uniformizer.
\par
A valued field often carries an angular component map modulo $M$,
or angular component map for short; it is a group homomorphism
$\ac:K^\times\to k^\times$, extended by putting $\ac(0)=0$, and
satisfying $\ac(x)=p(x)$ for all $x$ with $v(x)=0$ (see \cite{P}).
\subsection{Iterated Laurent series fields}\label{sect:laurent}
 We define iterated
Laurent series fields by induction. Let $L((t_1))$ be the field of
(formal) Laurent series in the variable $t_1$ over $L$ and let
$L((t_1))\ldots((t_{n-1}))((t_n))$ be the field of (formal)
Laurent series in the variable $t_n$ over
$L((t_1))\ldots((t_{n-1}))$. On a field $L((t_1))\ldots((t_n))$ we
can put many valuations, for example the valuation  $v_n$ taking
values in the lexicographically ordered $n$-fold product of $\Z$,
defined as follows. If $n=1$, then we put $v_1(x)=s\in\Z$ whenever
$x=\sum_{i\geq s}a_st_1^i$ with $a_s\not=0$ and $a_i\in L$. For
general $n$, and $x=\sum_{i\geq s}a_st_n^i$, where $a_s\not=0$ and
$a_i\in L((t_1))\ldots((t_{n-1}))$, we put
$v_n(x)=(s,v_{n-1}(a_s))\in\Z^n$. Remark that the valuation ring
with respect to the valuation $v_n$ is Henselian.
 \subsection{Grothendieck rings}
Let $\cL$ be a language and let $M$ be a model for $\cL$. For
$\cL$-definable sets $X\subset M^m$, $Y\subset M^n$, $m,n>0$, a
$\cL$-definable bijection $X\to Y$ is called an $\cL$-isomorphism
and we write $X\cong_{\cL} Y$, or $X\cong Y$ if the context is
clear, if $X$ and $Y$ are $\cL$-isomorphic. (Definable always
means definable with parameters.) For definable $X$ and $Y$, we
can choose disjoint definable sets $X',Y'\subset K^{m'}$ for some
$m'>0$, such that $X\cong X'$ and $Y\cong Y'$, and then we define
the \emph{disjoint union} $X\sqcup Y$ of $X$ and $Y$ up to
isomorphism as $X'\cup Y'$. By the Grothendieck group $K_0(M,\cL)$
of the structure $(M,\cL)$ we mean the group generated by symbols
$[A]$, for $A$ a $\cL$-definable set, with the relations
$[A]=[A']$ if $A\cong_{\cL} A'$ and $[A]=[B]+[C]$ if $A$ is the
disjoint union of $B$ and $C$. The group $K_0(M,\cL)$ carries a
multiplicative structure induced by $[A\times B]=[A][B]$, where
$A\times B$ is the Cartesian product of definable sets. The
so-obtained ring is called the Grothendieck ring and for a
$\cL$-definable set $X$ we write $[X]$ for the image of $X$ in
$K_0(M,\cL)$.
\par
Let $T$ be a theory in some language $\cL$. A formula $\varphi$
with free variables $x_1,\ldots x_n$ determines a set in $M^n$ for
any model $M$ of $T$. On these sets we can define a disjoint union
operation and Cartesian products in the natural way. The
Grothendieck group $K_0(T,\cL)$ is the group generated by symbols
$[\varphi]$, for $\varphi$ a $\cL$-formula, with the relations
$[\varphi]=[\varphi']$ if there is a formula which yields the
graph of a $\cL$-definable bijection between the sets defined by
$\varphi$ and $\varphi'$ in any model of $T$, and the relation
$[\varphi]=[\psi]+[\psi']$ if $\varphi$ is the disjoint union of
$\psi$ and $\psi'$.  This group carries a multiplicative structure
induced by the Cartesian product of definable sets and the so
obtained ring is called the Grothendieck ring of $(T,\cL)$.
\section{Languages of Denef - Pas}\label{sect:Pas}
Let $K$ be a valued field, with a valuation map $v:K\to
G\cup\{\infty\}$ for some ordered group $G$, and an angular
component map $\ac:K\to k$, where $k$ is the residue field. Let
$\cL_k$ be an arbitrary expansion of $\Lr$ and let $\cL_G$ be an
arbitrary expansion of the language of ordered groups with
infinity, namely $(+,-,0,\infty,\leq)$. A language of Denef - Pas
can in fact be either language in a wide variety of languages; it
is always a three-sorted language of the form
$(\cL_k,\Lr,\cL_G,v,\ac)$, with as sorts:
\begin{itemize}
\item[(i)] a $k$-sort for the residue field-sort,
\item[(ii)] a $K$-sort for the valued field-sort, and
\item[(iii)] a $G$-sort for the value group-sort.
\end{itemize}
The language $\Lr$ is used for the $K$-sort, $\cL_k$ for the
$k$-sort and $\cL_G$ for the $G$-sort. The function symbol $v$
stands for the valuation map $K\to G\cup\{\infty\}$ and $\ac$
stands for an angular component map $K\to k$ (in fact, this is an
angular component map modulo the maximal ideal $M$). A structure
for a language of Denef - Pas is denoted $(k,K,G\cup\{\infty\})$,
where $k, K$ and $G$ are as above.
\par
Remark that if $G=\Z$, namely if $K$ is a $\Z$-valued field, there
exists a natural angular component map $\ac:K\to k$ sending
$x\not=0$ to $t^{-v(x)}x\bmod M$, where $t$ is a uniformizer of
the valuation ring. More generally, if the value group of $K$ is
$\Z^n$, and $t_1,\ldots,t_n$ are field elements such that
$v(t_1)=(1,0,\ldots,0),\ldots,v(t_n)=(0,\ldots,0,1)$ forms a set
of generators of $\Z^n$, there is a natural angular component map
$\ac:K\to k$ given by $\ac(x)=x\prod_i t_i^{-r_i}\bmod M$, where
$v(x)=(r_1,\ldots,r_n)$. These angular component maps are
canonical up to the choice of $t$ and $t_i$. Languages of Denef -
Pas are denoted $\Lpas$.
\section{Calculations of Grothendieck rings}\label{sect:calcul}
Let $K$ be a field and $\cL$ an expansion of $\Lr$. We  write
$P_n(K)$ or $P_n$ for the $n$-th powers in $K^\times$. For $n>1$
we put
\[r_n(K)=\sharp\{x\in K\mid x^n=1\}\]
and
\[s_n(K)=[K^\times:P_n(K)]
\]
which is either a positive integer or $\infty$.
\begin{definition}\label{def:power}
For $n>1$ we put
\[
\lambda_n(K,\cL)=\frac{s_n(K)}{r_n(K)}
\]
 if the following conditions are satisfied
\begin{itemize}
\item  $s_n(K)<\infty$ and $\frac{s_n(K)}{r_n(K)}\in\Z$ ;\\
\item there exists a $\cL$-definable $n$-th root function. This means that
there exists a definable set $\sqrt[n]{P_n}$ and a definable
bijection $\sqrt[n]{}:P_n(K)\to\sqrt[n]{P_n}$ such that
$(\sqrt[n]{x})^n=x$ for each $x\in P_n(K)$.
\end{itemize}
If one of the above conditions is not satisfied, we put
$\lambda_n(K,\cL)=1$. We define $H(K,\cL)$ as the subgroup of $\Z$
generated by the numbers
\[
\lambda_n(K,\cL)-1
\]
for all $n>1$.
\end{definition}
Remark that if $\cL'$ is an expansion of $\cL$, then there is a
group inclusion $H(K,\cL)\to H(K,\cL')$. Let $\Lv=(\Lr,R)$ be the
language of rings with an extra $1$-ary relation symbol $R$ which
corresponds to a valuation ring inside the model. If the model is
a valued field, we take the natural interpretations.
\begin{lemma}\label{lem:HK}
Let $K$ be a field and $\cL$ an expansion of $\Lr$. For each
positive number $m\in H(K,\cL)$, there exists a $\cL$-definable
bijection
\[
\bigsqcup_{i=1}^{m+1} K^\times\to K^\times,
\]
and thus, in $K_0(K,\cL)$,
\[m [K^\times] =0.
\]
 Moreover, if $K$ is a valued field and $\cL$ is an expansion of $\Lv$, there
exists a $\cL$-definable bijection
\[
\bigsqcup_{i=1}^{m+1} (R\setminus\{0\})\to R\setminus\{0\},
\]
and thus, in $K_0(K,\cL)$,
\[m [R\setminus\{0\}] =0.
\]
\end{lemma}
\begin{proof} We first prove that $K^\times\cong \bigsqcup_{i=1}^{\lambda_n} K^\times$
for all $\lambda_n=\lambda_n(K,\cL)$, $n>1$. If $\lambda_n=1$
there is nothing to prove, so let $\lambda_n>1$. Remark that for
each $x\in K^\times$ and each definable set $A\subset K$ there is
a $\cL$-isomorphism $xA\cong A$. With the notation of definition
\ref{def:power}, the sets $x\sqrt[n]{P_n}$ form a partition of
$K^\times$ when $x$ runs over the $n$-th roots of unity. This
gives $\bigsqcup_{i=1}^{r_n}\sqrt[n]{P_n}\cong K^\times$. Since
$K^\times$ is the disjoint union of all cosets of $P_n$ inside
$K^\times$, we find $\bigsqcup_{i=1}^{s_n} P_n\cong K^\times$.
Combining with the isomorphism $P_n\cong\sqrt[n]{P_n}$ we
calculate:
 \[
 K^\times \cong \bigsqcup_{i=1}^{s_n} P_n \cong \bigsqcup_{i=1}^{s_n} \sqrt[n]{P_n}
\cong
\bigsqcup_{i=1}^{\lambda_n}(\bigsqcup_{i=1}^{r_n}\sqrt[n]{P_n})
  \cong
\bigsqcup_{i=1}^{\lambda_n}K^\times,
  \]
where $s_n=s_n(K)$ and $r_n=r_n(K)$. Now let $m>0$ be in
$H(K,\cL)$ and let $n>1$, $s>0$ be integers. By what we just have
shown, we can add $\lambda_n-1$ disjoint copies of $K^\times$ to
$\bigsqcup_{i=1}^s K^\times$, in the sense that $\bigsqcup_{i=1}^s
K^\times\cong\bigsqcup_{i=1}^{s+\lambda_n-1}K^\times$. Similarly,
if $s>\lambda_n-1$, we can subtract $\lambda_n-1$ disjoint copies
of $K^\times$ from $\bigsqcup_{i=1}^s K^\times$, to be precise,
$\bigsqcup_{i=1}^s
K^\times\cong\bigsqcup_{i=1}^{s-\lambda_n+1}K^\times$. The lemma
follows since the numbers $\lambda_n -1$ generate $H(K,\cL)$.
\par
If $\cL$ is an expansion of $\Lv$, we have the same isomorphisms
and the same arguments for $R\setminus\{0\}$ instead of
$K^\times$, working with $R\cap P_n$ and $R\cap \sqrt[n]{P_n}$
instead of $P_n$ and $\sqrt[n]{P_n}$.
\end{proof}
\begin{definition} Let $R$ be a valuation ring such that the value group
has a minimal strictly positive element. Let $\pi\in R$ have
minimal strictly positive valuation. Write $M$ for the maximal
ideal of $R$. Let $\ac$ be an angular component map $K\to k$,
where $k$ is the residue field. We define the set $R^{(1)}$ as
\[
R^{(1)}=\{x\in R\mid \ac(x)=1\}.
\]
\end{definition}
The set $R^{(1)}$ is not necessarily definable in the language
$\Lr$. If $R^{(1)}$ is definable in some language $\cL$ we have
the following criterion. Remark also that a minimal strictly
positive element in the value group necessarily is unique.
\begin{proposition}\label{prop:HK=Z} Let $K$ be a valued field. Suppose
that the value group has a minimal strictly positive element and
let $\pi\in R$ have this minimal strictly positive valuation. Let
$\cL$ be an expansion of $\cL_{\rm v}$ and let $\ac$ be an angular
component map $K\to k$. If $R^{(1)}$ is $\cL$-definable and
$H(K,\cL)=\Z$, then
\[K_0(K,\cL)=0,\qquad K\cong_\cL K^\times,\ \mbox{ and }\ R\cong_\cL
R\setminus\{0\}.
\]
\end{proposition}
\begin{proof} We first prove that $K_0(K,\cL)=0$. We may suppose
that $\ac(\pi)=1$, otherwise we could replace $\pi$ by $\pi/a$
where $a$ is an arbitrary element with $v(a)=0$ and
$\ac(a)=\ac(\pi)$. The following is a $\cL$-isomorphism
\[
R\sqcup R^{(1)}\to R^{(1)}:\left\{\begin{array}{lcl}x\in R & \mapsto & 1+\pi x,\\
x\in R^{(1)} & \mapsto & \pi x.\end{array}\right.
 \]
This implies, in $K_0(K,\cL)$, that $[R]+[R^{(1)}]=[R^{(1)}]$, and
thus after cancellation, $[R]=0$. By lemma \ref{lem:HK} and
because $H(K,\cL)=\Z$, also $[R\setminus\{0\}]=0$. The following
calculation implies $K_0(K,\cL)=0$:
\[
0=[R]=[R\setminus\{0\}]+[\{0\}]=[\{0\}]=1.
\]
We have
\[
[\{0\}]=1
\]
because $[\{0\}]$ is the multiplicative unit in $K_0(K,\cL)$.
 \par
Next we prove $R\cong R\setminus\{0\}$, by taking translates and
applying homotheties to the occurring sets. We make all occurring
disjoint unions explicit. Write $f_1$ for the isomorphism
\[f_1:1+\pi^2(R\setminus\{0\})\to \pi^2(R\setminus\{0\})\cup1+\pi^2(R\setminus\{0\}),\]
given by lemma  \ref{lem:HK}, it is an isomorphism from one copy
of $R\setminus\{0\}$ onto two disjoint copies of
$R\setminus\{0\}$.
Define the function $f_2$ on $\pi^2 R\cup \pi+\pi^2 R^{(1)}$ by
\[f_2:\pi^2 R\cup \pi+\pi^2 R^{(1)}\to \pi+\pi^2 R^{(1)}:
\left\{\begin{array}{rcl}
 \pi^2 x & \mapsto & \pi+\pi^2(1+\pi x),\\
 \pi+\pi^2 x & \mapsto & \pi+\pi^2(\pi x),
\end{array}\right.\]
then $f_2$ is an isomorphism from the disjoint union of $R$ and
$R^{(1)}$ to a copy of $R^{(1)}$. Finally, we find
$\cL$-isomorphisms:
\[f:R\to R\setminus\{0\}:x\mapsto\left\{
\begin{array}{ll}
 f_1(x) & \mbox{if }x\in 1+\pi^2(R\setminus\{0\}),\\
 f_2(x) & \mbox{if }x\in\pi^2 R\cup \pi+\pi^2 R^{(1)},\\
 x      & \mbox{else}
\end{array}\right.\]
and
\[K\to K^\times :x\mapsto\left\{
\begin{array}{ll}
 f(x) & \mbox{if }x\in R,\\
  x   & \mbox{else.}
\end{array}\right.\]
\end{proof}
Proposition \ref{prop:HK=Z} immediately yields the triviality of
the Gro\-then\-dieck rings of $\Q_p$ and of $\F_q((t))$ with
characteristic different from 2, which was originally proven in
\cite{CH}, cfr.~chapter \ref{chap:CH}. Theorems
\ref{thm:bijection:Qp} and \ref{thm:bijection:Fq} of the present
paper are generalizations. In case that $H(K,\cL)$ is different
from $\Z$, we formulate the following criterion. (The argument of
this criterion is similar to the proof of \cite{CH}, Thm.~1, see
also Thm.~\ref{Gring=0} of chapter \ref{chap:CH}).
\begin{proposition}\label{prop:crit2}
Let $K$ be a valued field. Suppose that the value group has a
minimal, strictly positive element and let $\pi\in R$ have this
minimal strictly positive valuation. Let $\cL$ be an expansion of
$\cL_{\rm v}$  and let $\ac$ be an angular component map $K\to k$.
If $R^{(1)}$ is $\cL$-definable, then
\[K_0(K,\cL)=0,\qquad K^2\cong_\cL K^2\setminus\{(0,0)\},\ \mbox{ and }\ R^2\cong_\cL
R^2\setminus\{(0,0)\}.
\]
\end{proposition}
\begin{proof}
We first prove that $K_0(K,\cL)=0$. As above we may suppose that
$\ac(\pi)=1$. The following is a $\cL$-isomorphism
\[
g_1:R\sqcup R^{(1)}\to R^{(1)}:\left\{\begin{array}{lcl}x\in R & \mapsto & 1+\pi x,\\
x\in R^{(1)} & \mapsto & \pi x.\end{array}\right.
 \]
As above, this implies that $[R]=0$ in $K_0(K,\cL)$.
\par
We show that the disjoint union of two copies of
$(R\setminus\{0\})^2$ is $\cL$-isomorphic to $(R\setminus\{0\})^2$
itself. Define the sets
\begin{eqnarray*}
X_1=\{(x,y)\in (R\setminus\{0\})^2|v(x)\leq v(y)\},\\
X_2=\{(x,y)\in (R\setminus\{0\})^2|v(x)>v(y)\},
\end{eqnarray*}
then $X_1,X_2$ form a partition of $(R\setminus\{0\})^2$. The
isomorphisms
\begin{eqnarray*}
(R\setminus\{0\})^2\to X_1: (x,y)\mapsto (x,xy),\\
(R\setminus\{0\})^2\to X_2: (x,y)\mapsto (\pi xy,y),
\end{eqnarray*}
imply that $(R\setminus\{0\})^2\sqcup (R\setminus\{0\})^2$ is
isomorphic to $X_1\cup X_2$ which is exactly
$(R\setminus\{0\})^2$. After cancellation, it follows that
$[(R\setminus\{0\})^2]=0$.
\par
Since $0=[R]=[R\setminus\{0\}]+[\{0\}]=[R\setminus\{0\}]+1$ we
have $[R\setminus\{0\}]=-1$. Together with
$0=[(R\setminus\{0\})^2]=[R\setminus\{0\}]^2$ this yields $1=0$,
so $K_0(K,\cL)$ is trivial. Write $g_2$ for the isomorphism
\[g_2:(R\setminus\{0\})^2\to (R\setminus\{0\})^2\sqcup (R\setminus\{0\})^2\]
\par
Now take the disjoint union of $R^{(1)}\times (R\setminus\{0\})$,
$R\times (R\setminus\{0\})$ and $(R\setminus\{0\})^2$ inside $R^2$
in some way, meaning that we take disjoint isomorphic copies
inside $R^2$ of the mentioned sets. Using the above isomorphisms
$g_1$ and $g_2$ in a clever way on these disjoint copies, it is
clear that we can remove one copy of $R\times (R\setminus\{0\})$
from $R^2$ and put one copy of $(R\setminus\{0\})^2$ back instead,
hence we find an isomorphism from $R^2$ to itself minus a point.
For details of this construction, we refer to the proof of
\cite{CH}, Thm.~1, cfr.~Thm.~\ref{Gring=0} of chapter
\ref{chap:CH}.
\end{proof}
\section{The proofs of theorems 1, 2, and 3}
\prooftitle{Proof of theorem \ref{thm:bijection:Qp}} Fix a field
$K$ as in the statement. Using Hensel's lemma, it is elementary to
calculate for each $n$ the numbers $r_n(K)$ and $s_n(K)$, and to
find that, for $n$ a prime number, $s_n(K)/r_n(K)$ is a power of
$n$. Further, it is not difficult to check that taking $n$-th
roots is definable. Therefore, the generator $\lambda_2-1$ of
$H(K,\Lr)$ is odd and $\lambda_3-1$ is even. This implies that
$H(K,\Lr)=\Z$.
\par
We calculate explicitly for $K=\Q_p((t))$, for the other fields of
the statement the arguments are completely similar, although,
notation can get more complicated.
\par
Let $v$ be the valuation on $K$ into $\Z\times\Z\cup\{\infty\}$
with lexicographical order, determined by: $v(x)=(s,r)$ for
$x=\sum_{i\geq s}a_it^i$ with $a_i\in\Q_p$ and $a_s\not=0$, the
$p$-adic valuation of $a_s$ being $r$. The valuation ring $R$ is
definable and can be described  by
\[
R=\{x\in K\mid 1+tx^2\in P_2(K)\ \&\ 1+px^2\in P_2(K)\}
\]
if $p\not=2$ and by
\[
R=\{x\in K\mid 1+tx^3\in P_3(K)\ \&\ 1+px^3\in P_3(K)\}
\]
if $p=2$. Write $M$ for the maximal ideal of $R$. Let $\ac:K\to
\F_p$ be the angular component $\ac(x)=p^{-r}t^{-s}x\bmod M$ for
nonzero $x$ with $v(x)=(s,r)$.
%
The set $R^{(1)}=\{x\in R\mid \ac(x)=1\}$ is definable since it is
the union of the sets
\[
p^it^jP_{p-1}(K),
\]
for $i,j=0,\ldots, p-2$.
\par
Now we can use proposition \ref{prop:HK=Z},
to find a $\Lr$-definable bijection $K\to K^\times$ and to find
that $K_0(K,\Lr)$ is trivial. This proves the proposition for
fields of the form $\Q_p((t))$. When $L$ is a finite field
extension of $\Q_p$ and $K$ an iterated Laurent series field over
$L$, there are $\Lr$-formula's playing the role of the formula's
above in the obvious way and the reader can make the adaptations.
 \qed
\prooftitle{Proof of Theorem \ref{thm:bijection:Fq}} Suppose for
simplicity that $K$ is the field $\big(\F_p((t_1))\big)((t_2))$,
where $p$ is a prime. The other cases are completely similar.
%
Let $v_2$ be the valuation on $K$ into $\Z\times\Z$ as in section
\ref{sect:laurent}; this is a valuation determined by:
$v(x)=(s,r)$ for a Laurent series $\sum_{i\geq s}a_it_2^i$ with
$a_i\in\F_q((t_1))$ and $a_s\not=0$ and $a_s=\sum_{i\geq r}
b_it_1^i$ with $b_r\not=0$ and $b_i\in\F_q$. Write $M$ for the
maximal ideal with respect to $v_2$. The valuation ring $R=\{x\mid
v_2(x)\geq 0\}$ is $\Lr$-definable because of the following
observation:
\[
R=\{x\in K\mid 1+t_2x^2\in P_2(K)\ \&\ 1+t_1x^2\in P_2(K)\},
\]
if $\char(K)\not=2$ and
\[
R=\{x\in K\mid 1+t_2x^3\in P_3(K)\ \&\ 1+t_1x^3\in P_3(K)\},
\]
if $\char (K)=2$.
\par
Let $\ac:K\to \F_p$ be the angular component map $x\mapsto
t_2^{-s}t_1^{-r}x\bmod M$ for nonzero $x$ with $v_2(x)=(s,r)$.
The set $R^{(1)}=\{x\in R\mid \ac(x)=1\}$ is definable since it is
the  union of the sets
\[
t_1^jt_2^iP_{p-1}(K)
\]
for $i,j=0,\ldots,p-2$. Now use proposition \ref{prop:crit2} to
find that $K_0(K,\Lr)$ is trivial and to find a $\Lr$-definable
bijection
\[g_3:R^2\to R^2\setminus\{(0,0)\}.\]
The following is a definable injection:
\[g_4:R^2\to R: (x,y)\mapsto x^p+t_1y^p,\]
and thereby, we can define the $\Lr$-isomorphism
\[
g_5:R\to R\setminus\{0\}: x\mapsto\left\{\begin{array}{ll}
g_4g_3(g_4^{-1}(x)) & \mbox{if } x\in g_4(R^2)\\
x & \mbox{else}.\end{array}\right.
\]
This finishes the proof.
 \qed
\prooftitle{Proof of theorem \ref{thm:bijection:R((t))}} Let $K$
be the field $\R((t_1))...((t_{n-1}))((t_n))$. Taking $m$-th roots
is clearly $\Lr$-definable, and, using the notation of definition
\ref{def:power} we have that $\lambda_2(K,\Lr)$ is a power of two
and $\lambda_3(K,\Lr)$ is a power of three. Therefore,
$H(\R((t)),\Lr)=\Z$. The existence of a $\Lpas$-definable
bijection $K\to K^\times$ and the triviality of $K_0(K,\Lpas)$ are
formal consequences of proposition \ref{prop:HK=Z}, because
$H(K,\Lr)=\Z$, the value group clearly has a unique minimal
strictly positive element and $R^{(1)}$ is $\Lpas$-definable.
 \qed
\prooftitle{Proof of theorem \ref{thm:bijection3}} Let $K$ be
$L((t_1))\ldots((t_n))$ where $L$ is a field of characteristic
$p>0$. The statement follows immediately from proposition
\ref{prop:crit2} using the definable injection
\[R^2\to R: (x,y)\mapsto x^p+t_1y^p\]
as in the proof of theorem \ref{thm:bijection:Fq}.
\qed

\clearemptydoublepage

\chapter[Classification of semialgebraic sets]{Classification of semialgebraic $p$-adic sets up to
semialgebraic bijection.}\label{chap:classification}
\begin{abstract}
\footnote{This chapter corresponds to \cite{C}.} We prove that two
infinite $p$-adic semialgebraic sets are isomorphic (i.e. there
exists a semialgebraic bijection between them) if and only if they
have the same dimension.
\end{abstract}
\section{Introduction}
In real semialgebraic geometry (as opposed to $p$-adic
semialgebraic geometry) the following classification is well-known
\cite{vdD}:
 \begin{quote}\textit{There exists a real semialgebraic bijection
between two real semialgebraic sets if and only if they have the
same dimension and Euler characteristic.}
\end{quote}
More generally L. van den Dries \cite{vdD} gave such a
classification for o-minimal expansions of the real field, using
the dimension and Euler characteristic as defined for $o$-minimal
structures. Since the semialgebraic Euler characteristic
 $\chi$ is in fact the canonical map from
the real semialgebraic sets onto the Grothendieck ring (see
\cite{CH} and chapter \ref{chap:CH}) of $\R$ (which is $\Z$), we
see that the isomorphism class of a real semialgebraic set only
depends on its image in the Grothendieck ring and its dimension.
\par
In this paper we treat the $p$-adic analogue of this
classification. The
 Grothendieck ring of $\Q_p$ is recently proved to be trivial by D.
Haskell and the author in \cite{CH} (cfr.~chapter \ref{chap:CH}),
so the analogue of the real case is a classification of the
$p$-adic semialgebraic sets up to semialgebraic bijection using
only the dimension. We give such a classification for the $p$-adic
semialgebraic sets and for finite field extensions of $\Q_p$,
using explicit isomorphisms of \cite{CH} (cfr.~chapter
\ref{chap:CH}), and the $p$-adic Cell Decomposition Theorem of J.
Denef \cite{Denef2}. The most difficult part in giving this
classification is to prove that for any semialgebraic set $X$
there is a finite partition into semialgebraic sets, such that
each part is isomorphic to a Cartesian product of one dimensional
sets, in other words semialgebraic sets have a rectilinearization.
Since all arguments hold also for finite field extensions of
$\Q_p$, we work in this more general setting.
\subsection*{Notation and terminology}
Let $p$ denote a fixed prime number, $\Q_p$ the field of $p$-adic
numbers and $K$ a fixed finite field extension of $\Q_p$. For
$x\in K$ let $v(x)\in\Z\cup\{+\infty\}$ denote the valuation of
$x$. Let $R=\{x\in K\mid v(x)\geq0\}$ be the valuation ring,
$K^\times=K\setminus\{0\}$ and for $n\in\N_0$ let $P_n$ be the set
$\{x\in K^\times\mid\exists y\in K\ y^n=x\}$. We call a subset of
$K^n$ {\em semialgebraic } if it is a Boolean combination (i.e.
obtained by taking finite unions, complements and intersections)
of sets of the form $\{x\in K^m\mid f(x)\in P_n\}$, with $f(x)\in
K[X_1,\ldots,X_m]$. The collection of semialgebraic sets is closed
under taking projections $K^m\to K^{m-1}$, even more: it consists
precisely of Boolean combinations of projections of the zero locus
of polynomials over $K$. Further we have that sets of the form
$\{x\in K^m \mid v(f(x))\leq v(g(x))\}$ with $f(x),g(x)\in
K[X_1,\ldots,X_m]$ are semialgebraic (see \cite{Denef2} and
\cite{Mac}). A function $f:A\to B$ is semialgebraic if its graph
is a semialgebraic set; if further $f$ is a bijection, we call $f$
 an {\em isomorphism } and we write $A\cong B$.
\par
 Let $\pi$ be a fixed element of $R$ with
$v(\pi)=1$, thus $\pi$ is a uniformizing parameter for $R$. For a
semialgebraic set $X\subset K$ and $k>0$ we write
\[
X^{(k)}=\{x\in X\mid x\not=0 \mbox{ and } v(\pi^{-v(x)}x-1)\geq
k,\ x\not=0\},
\]
which is semialgebraic (see \cite{Denef2}, Lemma 2.1); $X^{(k)}$
consists of those points $x\in X$ which have a $p$-adic expansion
$x=\sum_{i=s}^\infty a_i\pi^i$ with $a_s=1$ and $a_i=0$ for
$i=s+1,\ldots,s+k-1$. By a \emph{finite partition} of a
semialgebraic set we mean a partition into finitely many
semialgebraic sets. Let $X\subset K^n$, $Y\subset K^m$ be
semialgebraic. Choose disjoint semialgebraic sets $X',Y'\subset
K^k$ for some $k$, such that $X\cong X'$ and $Y\cong Y'$, then we
define the disjoint union of $X$ and $Y$ up to isomorphism as
$X'\cup Y'$. In the introduction of \cite{CH} (cfr.~chapter
\ref{chap:CH}), it is shown that we can take $k=\max(m,n)$, i.e.
we can realize the disjoint union without going into higher
dimensional affine spaces.
\par
\section[Preliminary results]{Preliminary and well-known results}
We recall some well-known facts.
\begin{lemma}[Hensel]\label{Hensel}
Let $f(t)$ be a polynomial over $R$ in one variable $t$, and let
$\alpha\in R,\ e\in\N$. Suppose that
$f(\alpha)\equiv0\bmod\pi^{2e+1}$ and $
  v(f'(\alpha))\leq e$, where $f'$ denotes the derivative of $f$. Then there
exists a unique $\bar\alpha\in R$ such that $f(\bar\alpha)=0$ and
$\bar\alpha\equiv\alpha\bmod\pi^{e+1}$.
\end{lemma}
\begin{cor}\label{corhensel}
Let $n>1$ be a natural number. For each $k> v(n)$, and $k'=k+ v(n)$ the
function $$ K^{(k)}  \to P_n^{(k')}:x\mapsto x^n $$ is an isomorphism.
\end{cor}
\par
The next theorem gives some concrete isomorphisms between one
dimensional sets.
\begin{proposition}[\cite{CH}, Prop. 2, cfr.~Prop.~\ref{p-adic} of chapter \ref{chap:CH}]\label{gring}
\item{(i)} The union of two disjoint
copies of $R\setminus\{0\}$ is isomorphic to $R\setminus\{0\}$.
\item{(ii)}  For each $k>0$ the union of two disjoint
copies of $R^{(k)}$ is isomorphic to $R^{(k)}$.
\item{(iii)} $R\cong R\setminus\{0\}$.
\end{proposition}
We deduce an easy corollary, also consisting of concrete isomorphisms.
\begin{cor}\label{corgring} For each $k$ we have isomorphisms
\item{(i)} $R^{(k)}\cong R\setminus\{0\}$,
\item{(ii)} $R\setminus \{0\}\cong K$.
\end{cor}
\begin{proof} (i) There is a finite partition
$R\setminus\{0\}=\bigcup_\alpha \alpha R^{(k)}$ with $v(\alpha)=0$, say
with $s$ parts. Then $R\setminus\{0\}$ is a fortiori isomorphic to the
union of $s$ disjoint copies of $R^{(k)}$,
which is by Proposition \ref{gring}(ii) isomorphic to $R^{(k)}$.\\
 (ii) The map
\[
 (\{0\}\times R)\cup (\{1\}\times R\setminus\{0\})\to K:\left\{\begin{array}{rcl}
                    (0,x) & \mapsto & x,\\
                    (1,x) & \mapsto & 1/(\pi x),\end{array}\right.\
\]
is a well-defined isomorphism. It follows that $K$ is isomorphic to the
disjoint union of $R$ and $R\setminus\{0\}$. Now use (i) and (iii) of
Proposition~\ref{gring}.
\end{proof}
\par Give $K^m$ the topology induced by the norm
$|x|=\max(|x_i|_p)$ with $|x_i|_p=p^{- v(x_i)}$ for $x=(x_1,\ldots
x_m)\in K^m$. P. Scowcroft and L. van den Dries \cite{SvdD} proved
there exists no isomorphism from an open set $A\subset K^m$ onto
an open set $B\subset K^n$ with $n\not= m$, so we can define the
dimension of semialgebraic sets as follows.
\begin{definition}[\cite{SvdD}]\textup{
 The \emph{dimension} of a semialgebraic set $X\not=\phi$ is the greatest
natural number $n$ such that we have a nonempty semialgebraic
subset $A\subset X$ and an isomorphism from $A$ to a nonempty
semialgebraic open subset of $K^n$. We put $\dim(\phi)=-1$.
 }\end{definition}
\par
P. Scowcroft and L. van den Dries \cite{SvdD} proved many good
properties of this dimension, for example that it is invariant
under isomorphisms.
\begin{proposition}[\cite{SvdD}]\label{classif:prop:dim}
Let $A$ and $B$ be semialgebraic sets, then the following is true:
\item{(i)} If $A\cong B$ then $\dim(A)=\dim(B)$,
\item{(ii)} $\dim(A\cup
B)=\max(\dim(A),\dim(B))$.
\item{(iii)} $\dim (A)=0$ if and only if $A$ is finite and
nonempty.
\end{proposition}
We will prove the converse of (i) for infinite semialgebraic sets.
\begin{lemma}\label{embed}
For any semialgebraic set $X$ of dimension $m\in\N_0$ there exists
a semialgebraic injection $X\to K^m$.
\end{lemma}
\begin{proof} By \cite{SvdD},~Cor. 3.1 there is a finite partition of
$X$ such that each part $A$ is isomorphic to a semialgebraic open
$A'\subset K^k$ for some $k\leq m$. Now realize the disjoint union
of the sets $A'$ without going into higher embedding dimension
(see the introduction).
\end{proof}
\par
We formulate the $p$-adic Cell Decomposition Theorem by J. Denef
\cite{Denef, Denef2}, which is the analogue of the real
semialgebraic Cell Decomposition Theorem.
\begin{theorem}[Cell Decomposition \cite{Denef,Denef2}] Let $x=(x_1,\ldots,x_m)$ and
$\hat x=(x_1,\ldots,x_{m-1}),\ m>0$. Let $f_i(\hat x,x_m)$,
$i=1,\ldots,r$, be polynomials in $x_m$ with coefficients which
are semialgebraic functions from $K^{m-1}$ to $K$. Let $n\in\N_0$
be fixed. Then there exists a finite partition of $K^m$ into sets
$A$ of the form
\[A=\{x\in K^m\mid \hat x\in D\mbox{ and }  v(a_1(\hat x))\sq_1
 v(x_m-c(\hat x))\sq_2 v(a_2(\hat x))\},\]
 such that
\[f_i(x)=u_i(x)^nh_i(\hat x)(x_m-c(\hat x))^{\nu_i},\mbox{ for each }\ x\in A,
\ i=1,\ldots,r,\] with $u_i(x)$ a unit in $R$ for each $x$,
$D\subset K^{m-1}$ semialgebraic, $\nu_i\in\N$, $h_i,a_1,a_2,c$
semialgebraic functions from $K^{m-1}$ to $K$ and $\sq_1$, $\sq_2$
either $\leq,<,$ or no condition.
\end{theorem}
The next lemma is also due to J. Denef \cite{Denef}.
\begin{lemma}[\cite{Denef}, Cor. 6.5]\label{order}
Let $b:K^m\to K$ be a semialgebraic function. Then there exists a
finite partition of $K^m$ into semialgebraic sets such that for
each part $A$ there are $e>0$ and polynomials $f_1, f_2\in
R[X_1,\ldots X_m]$ such that
\[
v(b(x))=\frac{1}{e} v(\frac{f_1(x)}{f_2(x)}), \mbox{ for each }\ x\in A,
\]
with $f_2(x)\not=0$ for each $x\in A$.
\end{lemma}
\section[Rec\-ti\-lin\-eariza\-tion]{Definable bijections and \\ rec\-ti\-lin\-eariza\-tion}
We give an application of the Cell Decomposition Theorem and Lemma
\ref{order}, inspired by similar applications in \cite{Denef}. For
details of the proof we refer to the proof of \cite{Denef}, Thm.
7.4. By $\lambda P_n$ with $\lambda=0$ we mean $\{0\}$.
\begin{lemma}\label{partition}
Let $X\subset K^m$ be semialgebraic and $b_j:K^m\to K$
semialgebraic functions for $j=1,\ldots,r$. Then there exists a
finite partition of $X$ s.t. each part $A$ has the form
\begin{eqnarray*}
A & = & \{x\in K^m\mid \hat x\in D, \  v(a_1(\hat
    x))\sq_1 v(x_m-c(\hat x))\sq_2 v(a_2(\hat x)),\\
 & & \qquad\qquad\qquad x_m-c(\hat x)\in \lambda P_n\},\end{eqnarray*}
and such that for each $x\in A$ we have
\[ v(b_j(x))=\frac{1}{e_j} v((x_m-c(\hat x))^{\mu_j}d_j(\hat x)),\]
with $\hat x=(x_1,\ldots, x_{m-1})$, $D\subset K^{m-1}$
semialgebraic, $e_j>0$, $\mu_j\in\Z$, $\lambda\in K$, $c,a_i,d_j$
semialgebraic functions from $K^{m-1}$ to $K$ and $\sq_i$ either
$<,\leq$ or no condition.
\end{lemma}
\begin{proof}
By Lemma \ref{order} we have a finite partition of $X$ such that
for each part $A_0$ there are $e_j>0$ and polynomials $g_j,
g'_j\in R[X_1,\ldots,X_m]$ with
\[
v(b_j(x))=\frac{1}{e_j} v(\frac{g_j(x)}{g'_j(x)}), \mbox{ for each }x\in
A_0,\ j=1,\ldots,r.
\]
Let $f_i$ be the polynomials which appear in a description of
$A_0$ as a Boolean combination of sets of the form $\{x\in K^m\mid
f(x)\in P_n\}$.. Apply now the Cell Decomposition Theorem as in
the proof of \cite{Denef}, Thm. 7.4 to the polynomials $f_i,g_j$
and $g_j'$ to obtain the lemma.
\end{proof}
\par
The proof of the next proposition is an application of both the
Cell Decomposition Theorem and some hidden Presburger arithmetic
in the value group of $K$; it is the technical heart of this
paper. If $l=0$ then $\prod_{i=1}^l R^{(k)}$ denotes the set
$\{0\}$.
\begin{definition}\textup{
We say that a semialgebraic function $f:B\to K$ satisfies
condition (\ref{condition}) (with constants $e,\mu_i,\beta$) if we
have constants $e\in\N_0,\mu_i\in\Z,\beta\in K$ such that each
$x=(x_i)\in B$ satisfies
\begin{equation}\label{condition}
 v(f(x))=\frac{1}{e} v(\beta\prod_i x_i^{\mu_i}).
 \end{equation}
 }\end{definition}
\begin{proposition}[Rectilinearization]\label{resolution}
 Let $X$ be a semialgebraic set and $b_j:X\to K$ semialgebraic functions for
$j=1,\ldots,r$. Then there exists a finite partition of $X$ into
semialgebraic sets such that for each part $A$ we have constants
$l\in\N$, $k\in\N_0$, $\mu_{ij}\in\Z$, $\beta_j\in K$, and an
isomorphism
\[f:\prod_{i=1}^l R^{(k)}\to A,\]
such that for each $x=(x_1,\ldots,x_l)\in\prod_{i=1}^l R^{(k)}$ we have
\[ v(b_j\circ f(x))= v(\beta_j\prod_{i=1}^{l}x_i^{\mu_{ij}}).\]
\end{proposition}
\begin{proof}
We work by induction on $m=\dim(X)$. Let $\dim(X)=1$ and $b_j:X\to
K$ semialgebraic functions, $j=1,\ldots,r$. By Lemma \ref{embed}
we may suppose that $X\subset K$. We reduce first to the case that
$X$ and $b_j$ have the special form (\ref{eq dim1}) (see below).
By Lemma \ref{partition} there is a partition such that each part
$A$ is either a point or of the form
\[A=\{x\in
K\mid v(a_1)\sq_1 v(x-c)\sq_2 v(a_2),\
 x-c\in \lambda P_n\},\]
and such that for each $x\in A$ we have $ v(b_j(x))=\frac{1}{e_j}
v(\beta_j(x-c)^{\mu_j}),$ with $a_i,c,\lambda,\beta_j\in K$,
$e_j>0$ and $\mu_j\in\Z$. We may assume that $\lambda\not=0$,
$a_1\not=0\not=a_2$, $\sq_i$ is either $\leq$ or no condition and
since the translation
\[
\{x\in K\mid v(a_1)\sq_1 v(x)\sq_2 v(a_2),\ x\in \lambda P_n\}\to
A:x\mapsto x+c
\]
is an isomorphism, we may also assume that $c=0$. If both $\sq_1$
and $\sq_2$ are no condition we can partition $A$ into parts
$\{x\in A \mid 0\leq v(x)\}$ and $\{x\in A\mid v(x)\leq-1\}$. It
follows that if $\sq_1$ is no condition we may suppose that
$\sq_2$ is $\leq$, then we can apply the isomorphism
 \[
 \{x\in K\mid v(\frac{1}{a_2})\leq v(x),\ x\in
\frac{1}{\lambda}P_n\}\to A:
 x\mapsto\frac{1}{x},
 \]
 and replace $\mu_j$ by $-\mu_j$.
This shows we can reduce to the case that $X$ has the form
\begin{equation}\label{eq dim1}
X=\{x\in K\mid v(a_1)\leq v(x)\sq_2 v(a_2),\ x\in\lambda P_n\},
\end{equation}
with $a_1\not=0\not=a_2$, $\lambda\not=0$, $\sq_2$ either $\leq$
or no condition and $v(b_j(x))=\frac{1}{e_j} v(\beta_j
x^{\mu_j})$ for each $x\in X$.\\

 \textbf{Case 1:} $\sq_2$ is $\leq$ (in equation (\ref{eq dim1})).\\
By Hensel's Lemma we can partition $X$ into finitely many parts of
the form $y+\pi^sR$ for some fixed $s>v(a_2)$ and with $v(a_1)\leq
v(y)\leq v(a_2)$ for each $y$. For each such part there is a
finite partition
$y+\pi^sR=\bigcup_{\gamma\in\Gamma}A_\gamma\cup\{y\}$, with
$A_\gamma=y+\pi^s\gamma R^{(1)}$ and $v(\gamma)=0$ for each
$\gamma$. The functions $f_\gamma:R^{(1)}\to A_\gamma:x\mapsto
y+\pi^s\gamma x$ are isomorphisms which satisfy $v(b_j\circ
f_\gamma(x))=\frac{1}{e_j}v(\beta_jy^{\mu_j})$ for all $x\in
R^{(1)}$. This last expression is independent of $x$, so there
exists $\beta_j'\in K$ such that $v(b_j\circ f_\gamma(x)) =
v(\beta'_j)$ for all $x\in R^{(1)}$. This proves Case 1.\\

 \textbf{Case 2:} $\sq_2$ is no condition (in equation (\ref{eq dim1})).\\
The map
\[f_1:R\cap\lambda' P_n\to X:x\mapsto a_1x,\]
with $\lambda'=\lambda/a_1$ is an isomorphism. Let $n'$ be a
common multiple of $e_1,\ldots,e_r$ and $n$. Choose $k> v(n')$ and
put $k'=k+ v(n')$. Let  $R\cap\lambda' P_n=\bigcup_\gamma
B_\gamma$ be a finite partition, with $B_\gamma=\gamma(R\cap
P_{n'}^{(k')})$ and $0\leq v(\gamma)<n'$. Now we have that the map
$f_\gamma:R^{(k)}\to B_\gamma:x\mapsto\gamma x^{n'}$ is an
isomorphism by corollary \ref{corhensel}. Let $g_\gamma$ be the
semialgebraic function $f_1\circ f_\gamma$, which is an
isomorphism from $R^{(k)}$ onto a semialgebraic set
$A_\gamma\subset X$. The sets $A_\gamma$ form a finite partition
of $X$. Put $\mu'_j=\mu_jn'/e_j$, then we have for each $x\in
R^{(k)}$ that
\[
v(b_j\circ g_\gamma(x))=\frac{1}{e_j} v(\beta_j(a_1\gamma
x^{n'})^{\mu_j})=\frac{1}{e_j}v(\beta_j(a_1\gamma)^{\mu_j})+v(x^{\mu_j'})
\]
\[
= v(\beta'_j)+v(x^{\mu'_j})=v(\beta'_jx^{\mu'_j}),
\]
with $\beta'_j$ in $K$ such that $ v(\beta'_j)=\frac{1}{e_j}
v(\beta_j(a_1\gamma)^{\mu_j})$. This proves case 2.
\par
Now let $\dim(X)=m>1$ and let $b_j:X\to K$ be semialgebraic
functions, $j=1,\ldots,r$. By Lemma \ref{embed} we may suppose
that
$X\subset K^m$.\\

 \textbf{Claim.}\textit{ We can
partition $X$ into finitely many semialgebraic sets such that for
each part $A$ we have an isomorphism of the form $f:D_1\times
D_{m-1}\to A$, with $D_1\subset K$ and $D_{m-1}\subset K^{m-1}$
semialgebraic, such that the functions $b_j\circ f$ satisfy
condition (\ref{condition}), i.e. there are constants
$e_j\in\N_0,\mu_{ij}\in\Z,\beta_j\in K$ such that each $x=(x_i)\in
D_1\times D_{m-1}$ satisfies
\[
 v(b_j\circ f(x))=\frac{1}{e_j} v(\beta_j\prod_i x_i^{\mu_{ij}}).
 \]}

If the claim is true, we can apply the induction hypotheses once to $D_1$
and the functions $x_1\mapsto x_1^{\mu_{1j}}$ and once to $D_{m-1}$ and
the functions $(x_2,\ldots,x_m)\mapsto\beta_j\prod_{i=2}^mx_i^{\mu_{ij}}$
for $j=1,\ldots, r$. It follows easily that we can partition $X$ such that
for each part $A$ there is an isomorphism $f:\prod_iR^{(k)}\to A$ such
that all $f\circ b_j$ satisfy condition (\ref{condition}) with constants
$e'_j,\mu'_{ij}$ and $\beta'_j$. Now we can proceed as in Case 2 for $m=1$
to make all $e'_j\in\N_0$ occurring
in condition (\ref{condition}) equal to 1. The proposition follows now immediately.\\

\textbf{Proof of the claim.} First we show we can reduce to the
case described in equation (\ref{eq dim>1}) below. Using Lemma
\ref{partition} and its notation, we find a finite partition of
$X$ such that each part $A$ has the form
\begin{eqnarray*}
A & = & \{x\in K^m\mid\hat x\in D,\  v(a_1(\hat
    x))\sq_1 v(x_m-c(\hat x))\sq_2 v(a_2(\hat x)),\\
 & &\qquad\qquad\qquad x_m-c(\hat x)\in \lambda P_n\},
\end{eqnarray*}
 and such that for each $x\in A$ we have
 \[
  v(b_j(x))=\frac{1}{e_j}
 v((x_m-c(\hat x))^{\mu_{mj}}d_j(\hat x)), \]
 with $\mu_{mj}\in\Z$.
Similar as for $m=1$, we may suppose that $c(\hat x)=0$ for all
$\hat x$. Apply now the induction hypotheses to the set $D\subset
K^{m-1}$ and the functions $a_1,a_2,d_j$. We find a finite
partition of $A$ such that for each part $A'$ we have an
isomorphism $f:B\to A'$, where $B$ is a set of the form
\[
\begin{array}{ll}
\{x\in K^m\mid
 &
  \hat x\in D',\
v(\alpha_1\prod_{i=1}^lx_i^{\eta_i})
\sq_1 v(x_m)\sq_2 v(\alpha_2\prod_{i=1}^lx_i^{\varepsilon_i}),
 \\
 & x_m\in\lambda P_n\},
 \end{array}
\]
with $D'=\prod_{i=1}^lR^{(k)}$, $l\leq m-1$, such that each
$b_j\circ f$ satisfies condition (\ref{condition}). We will
alternately partition further and apply isomorphisms to the parts
which compositions with $b_j$ will always satisfy condition
(\ref{condition}). By the induction hypotheses we may suppose that
$\lambda\not=0$ and $\dim(D')=m-1$, i.e.
$D'=\prod_{i=1}^{m-1}R^{(k)}$. Analogously as for $m=1$ we may
suppose that $\alpha_1\not=0\not=\alpha_2$, $\sq_2$ is either
$\leq$ or no condition and $\sq_1$ is the symbol $\leq$ (possibly
after partitioning or applying
$x\mapsto(x_1,\ldots,x_{m-1},1/x_m)$).
\par
Choose $\bar k>v(n)$ and put $k'=\bar k+v(n)$. We may suppose that
$k'>k$, so we have a finite partition $B=\bigcup_\gamma B_\gamma$
with $\gamma=(\gamma_1,\ldots,\gamma_m)\in K^m$, $0\leq
v(\gamma_i)<n$ and $B_\gamma=\{x\in B \mid x_i\in\gamma_i
P_n^{(k')}\}$. Now we have isomorphisms
\[
f_\gamma:C_\gamma\to
B_\gamma:x\mapsto(\gamma_1x_1^n,\ldots,\gamma_mx_m^n),
\]
with
\[
\begin{array}{ll}
C_\gamma=\{x\in (\prod_{i=1}^{m-1}R^{(\bar k)})\times K^{(\bar
k)}\mid  & v(\alpha'_1\prod_{i=1}^{m-1}x_i^{\eta_i}) \leq
v(x_m)\sq_2
 \\
 &
  v(\alpha'_2\prod_{i=1}^{m-1}x_i^{\varepsilon_i})\},
  \end{array}
\]
for appropriate choice of $\alpha_i'\in K$. Put
$\nu_i=\varepsilon_i-\eta_i$, $\beta=\alpha'_2/\alpha'_1$, then we have
the isomorphism
\[
\begin{array}{c}
\{x\in \prod_{i=1}^mR^{(\bar k)}\mid v(x_m)\sq_2
 v(\beta\prod_{i=1}^{m-1} x_i^{\nu_i})\} \to C_\gamma:\\
 x\mapsto(x_1,\ldots,x_{m-1},\alpha'_1x_m\prod_{i=1}^{m-1}x_i^{\eta_i}).
\end{array}
\]
\par
If $\sq_2$ is no condition, the claim is trivial. It follows that
we can reduce to the case that we have an isomorphism
\begin{equation}\label{eq dim>1}
 f:E=\{x\in \prod_{i=1}^mR^{(\bar k)}\mid v(x_m)\leq
 v(\beta\prod_{i=1}^{m-1} x_i^{\nu_i})\}\to X
\end{equation}
with $\beta\not=0$, $\bar k>0$, and $\nu_i\in\Z$, such that each
$b_j\circ f$  satisfies condition (\ref{condition}).
\par
Suppose we are in the case described in (\ref{eq dim>1}). If
$\nu_i\leq0$ for $i=1,\ldots,m-1$ then we have a finite partition
$E=\cup_sE^{=s}$, with $s\in\{0,1,\ldots,v(\beta)\}$ and
$E^{=s}=\{x\in E\mid v(x_m)=s\}$. Also,
$E^{=s}=\{(x_1,\ldots,x_{m-1})\mid\exists x_m\ (x_1,\ldots,x_m)\in
E^{=s}\}\times \{x_m\in R^{(\bar k)}\mid v(x_m)=s\}$ and the claim
follows.
\par
Suppose now that $\nu_1>0$ in (\ref{eq dim>1}). First we prove the
proposition when $\nu_1=1$, using some implicit Presburger arithmetic on
the value group. We can partition $E$ into parts $E_1$ and $E_2$, with
\begin{eqnarray*}
 E_1 & = & \{x\in
 E\mid v(x_m)< v(\beta\prod_{i=2}^{m-1}x_i^{\nu_i})\},\\
 E_2 & = & \{x\in
 E\mid v(\beta\prod_{i=2}^{m-1}x_i^{\nu_i})\leq v(x_m)\}\\
 & = & \{x\in\prod_{i=1}^m R^{(\bar k)}\mid v(\beta\prod_{i=2}^{m-1}x_i^{\nu_i})\leq
  v(x_m)\leq v(\beta x_1\prod_{i=2}^{m-1}x_i^{\nu_i})\}.
\end{eqnarray*}
Since $ v(\beta\prod_{i=2}^{m-1}x_i^{\nu_i})\leq
 v(x_1\beta\prod_{i=2}^{m-1}x_i^{\nu_i})$ for $x\in E_1$, it follows that
 \[
 E_1=R^{(\bar k)}\times\{(x_2,\ldots,x_m)\in\prod_{i=2}^mR^{(\bar k)}\mid v(x_m)<
 v(\beta\prod_{i=2}^{m-1}x_i^{\nu_i})\},
 \]
and the restrictions $b_j\circ f|E_1$ satisfy condition
(\ref{condition}).\\
 As for $E_2$, let $D_{m-1}$ be the set
\[
D_{m-1}=\{(x_2,\ldots,x_m)\in\prod_{i=2}^mR^{(\bar k)}\mid v
 (\beta\prod_{i=2}^{m-1}x_i^{\nu_i})\leq v(x_m)\}.
\]
We may suppose that $\beta\in K^{(\bar k)}$, then the map
\[
R^{(\bar k)}\times D_{m-1}\to E_2:
x\mapsto(\frac{x_1x_m}{\beta\prod_{i=2}^{m-1}x_i^{\nu_i}},x_2,\ldots,x_m),
\]
can be checked by elementary Presburger arithmetic to be an isomorphism.
This proves the claim when $\nu_1=1$.
\par
Suppose now that $X$ is of the form described in (\ref{eq dim>1}) and
$\nu_1>1$. We prove we can reduce to the case $\nu_1=1$ by partitioning
and applying appropriate power maps. Choose $\tilde k>v(\nu_1)$ and put
$\tilde k'=\tilde k+v(\nu_1)$. We may suppose that $\tilde k\geq\bar k$,
so we have a finite partition $E=\bigcup_\alpha E_\alpha$, with
$\alpha=(\alpha_1,\ldots,\alpha_m)\in K^m$, $v(\alpha_1)=0$, $0\leq
v(\alpha_i)<\nu_1$ for $i=2,\ldots,m$ and
\[E_\alpha=\{x\in E\mid x_1\in \alpha_1R^{(\tilde k)},\ x_i\in \alpha_i
P_{\nu_1}^{(\tilde k')} \mbox{ for } i=2,\ldots,m\}.\] By corollary
\ref{corhensel} we have isomorphisms
\[f_\alpha:C_\alpha\to E_\alpha:x\mapsto(\alpha_1x_1,\alpha_2x_2^{\nu_1},\ldots,
\alpha_mx_m^{\nu_1}),\]
 with $C_\alpha=\{x\in\prod_{i=1}^m R^{(\tilde k)}\mid v(x_m)\leq v
(\beta' x_1\prod_{i=2}^{m-1}x_i^{\nu_i})\},$ where $\beta'\in K^\times$
depends on $\alpha$. This reduces the problem to the case described in
(\ref{eq dim>1}) with $\nu_1=1$ and thus the proposition is proved.
\end{proof}
\begin{remark} Proposition~\ref{resolution} can be strengthened, adding conditions on
the Jacobians of the isomorphisms, to become useful for $p$-adic
integration.\end{remark}
\section{Classification of semialgebraic sets}
\begin{theorem}\label{dim}
Let $X$ be a semialgebraic set, then either $X$ is finite or there
exists a semialgebraic bijection $X\to K^k$ with $k\in\N_0$ the
dimension of $X$.
\end{theorem}
\begin{proof} We give a proof by induction on dim$(X)=m$. Let dim$(X)=1$. Use
Proposition \ref{resolution} to partition $X$ such that each part is
isomorphic to either $R^{(k)}$ or a point. By combining the isomorphisms
of Proposition~\ref{gring} and Corollary~\ref{corgring}, it follows that
$X\cong K$.
\par
 Now suppose dim$(X)=m>1$. Proposition \ref{resolution} together with the case $m=1$
implies that we can finitely partition $X$ such that each part is
isomorphic to $K^l$, for some $l\in\{0,\ldots,m\}$, with
$K^0=\{0\}$. By proposition \ref{classif:prop:dim} at least one
part must be isomorphic to $K^m$. Suppose that $A$ and $B$ are
disjoint parts, such that $A\cong K^l$ and $B\cong K^m$, with
$l\in\{0,\ldots,m\}$. It is enough to prove that $A\cup B\cong
K^m$. First suppose that $l=0$, so $A$ is a singleton $\{a\}$.
Since $m>1$ there exists an injective semialgebraic function
$i:R\to A\cup B$ such that $i(R\setminus\{0\})\subset B$ and
$i(0)=a$. It follows that $A\cup B\cong B\cong K^m$ since $R\cong
R\setminus\{0\}$ (Proposition \ref{gring}). If $1\leq l$ we have
$A\cup B \cong K\times(A'\cup B')$, for some disjoint sets
$A'\cong K^{l-1}$ and $B'\cong K^{m-1}$.  By induction we find
$A'\cup B'\cong K^{m-1}$ and thus $A\cup B\cong K^m$. This proves
the theorem.
\end{proof}
We obtain as a corollary of Theorem \ref{dim} the following
classification of the $p$-adic semialgebraic sets.
\begin{cor}
Two infinite semialgebraic sets are isomorphic if and only if they
have the same dimension.
\end{cor}
\begin{remark}
The field $\F_q((t))$ of Laurent series over the finite field is
often considered to be the characteristic $p$ counterpart of the
$p$-adic numbers. D. Haskell and the author \cite{CH}
(cfr.~chapter \ref{chap:CH}), proved that the Grothendieck ring of
$\F_q((t))$ in the language of rings is trivial, analogously as
for $\Q_p$. In the proof of Thm.~2 in \cite{CH}
(cfr.~Thm.~\ref{FF} of chapter \ref{chap:CH}), it is also shown
that there is a definable injection from the plane $\F_q((t))^2$
into the line $\F_q((t))$, which makes it plausible we cannot
define an invariant for this field which behaves like a (good)
dimension. It is an open question to classify the definable sets
of $\F_q((t))$ up to definable bijection.
\end{remark}

\clearemptydoublepage

\chapter[Analytic cell decomposition and integrals]{Analytic p-adic Cell
Decomposition and Integrals}\label{chap:cell}
%
 \begin{abstract}
\footnote{This chapter corresponds to \cite{Ccell}.}
  Roughly speaking, the semialgebraic cell
decomposition theorem for $p$-adic numbers describes piecewise the
$p$-adic valuation of polynomials  (and more generally of
semialgebraic $p$-adic functions), the pieces being geometrically
simple sets, called cells. In this paper we prove a  cell
decomposition theorem which describes piecewise the valuation of
analytic functions (and more generally of subanalytic functions),
the pieces being subanalytic cells. We use this cell decomposition
theorem to solve a conjecture of Denef on $p$-adic subanalytic
integration; the conjecture describes the parameter dependence of
analytic $p$-adic integrals. We also classify subanalytic sets up
to subanalytic bijection.
\end{abstract}
\section{Introduction}
Let $p$ denote a fixed prime number, $\Z_p$ the ring of $p$-adic
integers, and $\Q_p$ the field of $p$-adic numbers.
\par
Let $f=(f_1,\ldots,f_r)$ be a map of restricted power series over
$\Z_p$ in the variables
$(\lambda,x)=(\lambda_1,\ldots,\lambda_s,x_1,\ldots, x_m)$; these
are power series converging on $\Z_p^{s+m}$. To $f$ one can
associate a parametrized $p$-adic integral
\begin{equation}\label{integral}
I(\lambda)=\int\limits_{\Z_p^m}\abs{f(\lambda,x)}\abs{dx},
\end{equation}
where $|dx|$ is the Haar measure on $\Z_p^m$ normalized so that
$\Z_p^m$ has measure $1$ and $|\cdot|$ denotes the $p$-adic norm.
In the case that the functions $f_i$ are polynomials the map $I$
has been studied by Igusa for $r=1$, by Lichtin for $r=2$,
and by Denef for arbitrary $r$ (see \cite{Igusa1, Igusa2, Igusa3},
\cite{Lichtin}
 and \cite{Denef1}). We prove the following
conjecture of Denef \cite{Denef1}:
\begin{theorem} The function $I$
is a simple $p$-exponential function. (A simple $p$-exponential
function is a $\Q$-linear combination of products of functions of
the form $v(h)$ and $p^{v(h')}$, where $h, h'$ are subanalytic
functions, $h(x)\not=0$, and $v(\cdot)$ is the $p$-adic
valuation.)
\end{theorem}
More generally, we prove that the algebra of simple
$p$-exponential functions is closed under $p$-adic integration
(theorem~\ref{thm:basic} and corollary~\ref{Thm:conjecture}),
where the integration operator assigns the value zero when the
integrated function is not integrable.
\par
The rationality of the Serre $p$-adic Poincar\'e series in the
analytic case, conjectured by Oesterl\'e \cite{Oest} and Serre
\cite{Serre}, and proven in \cite{DvdD}, can immediately be
obtained as a corollary of these integration theorems. This is
because it is well-known how to express this Poincar\'e series as
a $p$-adic integral (see \cite{Denef1}). The theory of $p$-adic
integration has also served as an inspiring example for the theory
of motivic integration and there are many connections to it (see
e.g.~\cite{DL} and \cite{DLinvent}).
 We extensively use the theory
of $p$-adic subanalytic sets, developed by Denef and van den Dries
in \cite{DvdD} in analogy to the theory of real subanalytic sets
(see e.g.~Hironaka~\cite{Hir}). A function $h:X\subset
\Z_p^n\to\Z_p^m$ is called subanalytic if its graph is a
subanalytic set and a set $X\subset \Z_p^m$ is called subanalytic
if $X$ is the image under the natural projection $\Z_p^{m+s}\to
\Z_p^m$ (for some $s\geq 0$) of a finite union of sets of the form
\begin{equation}\label{eq:sub}
\{x\in \Z_p^{m+s} \mid f(x)=0,\ g_1(x)\in P_{n_1},\ \ldots,\
g_m(x)\in P_{n_m}\},
\end{equation}
where $f$ and $g_i$ are restricted analytic functions and $P_n$,
$n>0$, stands for the set of $n$-th powers in
$\Q_p^\times=\Q_p\setminus\{0\}$. Remark that a set like in Eq.
(\ref{eq:sub}) is itself a projection of analytic subsets of some
$\Z_p^{m+s+s'}$ (in the sense of \cite{Bour}), hence the name
subanalytic.
\par
The key result of the present paper is a cell decomposition
theorem for subanalytic sets, in perfect analogogy to cell
decomposition for semialgebraic sets (see \cite{Denef},
\cite{Denef2}). In \cite{Denef} cell decomposition is used to
prove the rationality of Igusa's local zeta function and of many
Poincar\'e series. Cell decompositions are perfectly suited to
study parameter dependence of $p$-adic integrals, often even
better than for example resolution of singularities is; many of
the applications (in for example ~\cite{Denef1} and \cite{Denef3})
can, up to now, not be proven with other techniques. In this
paper, we completely avoid the use of resolutions of
singularities. Roughly speaking, $p$-adic cell decompositions
describe the norm of given functions after partitioning the domain
of the functions in finitely many basic sets, called cells. The
cell decomposition is strong enough to reprove many results of
\cite{DvdD} in a very short and direct way and to develop, for
example, the dimension theory of subanalytic sets in a way (common
in model theory), based on cell decompositions.
 \par In section \ref{section:classification} we use cell
decomposition to prove the following remarkable classification:
\begin{theorem}
\label{classification} Let $X\subset \Q_p^m$ and $Y\subset \Q_p^n$
be infinite subanalytic sets, then there exists a subanalytic
bijection $X\to Y$ if and only if $\dim(X)=\dim(Y)$.
\end{theorem}
This classification of subanalytic sets is similar to the
classification of semialgebraic sets in \cite{C}, see chapter
\ref{chap:classification}. All results of the paper also hold in
finite field extensions of $\Q_p$.
\subsection*{Terminology and Notation} Let $p$ denote a fixed prime number, $\Q_p$ the
field of $p$-adic numbers and $K$ a fixed finite field extension
of $\Q_p$ with valuation ring $R$. For $x\in K^\times$ let
$v(x)\in\Z$ denote the $p$-adic valuation of $x$ and
$|x|=q^{-v(x)}$ the $p$-adic norm, with $q$ the cardinality of the
residue class field. For $n>0$ we let $P_n$ be the multiplicative
group of the $n$-th powers in $K^\times$ and $\lP$ denotes
$\{\lambda x\mid x\in P_n\}$ for $\lambda\in K$.
\par
For $X=(X_1,\ldots,X_m)$ let $K\langle X\rangle$ be the ring of
restricted power series over $K$ in the variables $X$; it is the
ring of power series $\sum a_iX^i$ in $K[[X]]$ such that $|a_i|$
tends to $0$ as $|i|\to\infty$. (Here, we use the multi-index
notation where $i=(i_1,\ldots,i_m)$, $|i|=i_1+\ldots+i_m$ and
$X^i=X_1^{i_1}\ldots X_m^{i_m}$.) For $x\in R^m$ and $f=\sum
a_iX^i$ in $K\langle X\rangle$ the series $\sum a_ix^i$ converges
to a limit in $K$, thus, one can associate to $f$ a
\emph{restricted analytic function} given by
\[
f:K^m\to K:x\mapsto
 \left\{\begin{array}{ll} \sum_i a_i x^i & \mbox{ if
}x\in R^m,\\
0 & \mbox{ else.}
 \end{array}\right.
\]
We define the algebra $\BK$ of \emph{basic subanalytic functions}
as the smallest $K$-algebra of functions in $\cup_m\{f:K^m\to K\}$
satisfying:
\begin{itemize}
\item[(i)]
for each $f\in\TK$ the associated restricted analytic function $f$
is in $\BK$;
\item[(ii)]
each rational function
 \[f/g:K^m\to K:x\mapsto \left\{\begin{array}{ll}
 f(x)/g(x) & \mbox{ if }g(x)\not=0,\\
 0 & \mbox{else},
 \end{array}\right.
 \]
  with $f,g$ polynomials over $K$ is in
$\BK$;
\item[(iii)] for each $f\in\BK$ in $n$ variables and
each $g_1,\ldots,g_n\in\BK$ in $m$ variables, the function
$f(g_1,\ldots,g_n)$ is in $\BK$.
\end{itemize}
A set $X\subset K^m$  is called (globally) subanalytic if it is a
finite union of sets of the form
\[\{x\in K^m\mid f(x)=0,\ g_1(x)\in P_{n_1},\ldots,g_k(x)\in P_{n_k}\},\]
where the functions $f$ and $g_i$ are in $\BK$ and $n_i>0$; this
corresponds to the definition given in the introduction by the
quantifier elimination result of Denef and van den Dries
\cite{DvdD}. We call a function $g:A\to B$ subanalytic if its
graph is a subanalytic set. We refer to \cite{DvdD},
\cite{Denef1}, and \cite{vdDHM} for the theory of subanalytic
$p$-adic geometry and to \cite{Lip} for the theory of rigid
subanalytic sets. In section \ref{proof cell decomp} we will use
the framework of model theory. We let $\Ldan$ be the first order
language consisting of the symbols
\[
+,\ -,\ \cdot,\ \{P_n\}_{n>0},
\]
the function symbol $^{-1}$ and an extra function symbol $f$ for
each restricted analytic function associated to restricted power
series in $\TK$. We consider $K$ as a $\Ldan$-structure using the
natural interpretations of the symbols of $\Ldan$.
\par
We mention the following fundamental result in the theory of
subanalytic sets.
\begin{theorem}[\cite{DvdD}, Corollary (1.6)]\label{Gabrielov}
The collection of subanalytic sets is closed under taking
complements, finite unions, finite intersections and projections
$K^{m+n}\to K^{m}$, thus, the image of a subanalytic set under a
subanalytic function is subanalytic.
\par
 Also, the following are equivalent for
$X\subset K^m$:
\item[(i)] $X$ is subanalytic;
\item[(ii)] $X$ is $\Ldan$-definable;
\item[(iii)] $X\subset K^m\subset \PPP^m$ is a subanalytic subset
of the $p$-adic manifold $\PPP^m$ in the sense of \cite[p 81]
{DvdD}. Here $\PPP^m$ is the set of $K$-rational points on the
$m$-dimensional projective space over $K$.
 \end{theorem}
 \begin{remark}
A $p$-adic subanalytic subset of a $p$-adic manifold, as referred
to in (iii) above, is similar to the notion of a real subanalytic
subset of a real analytic manifold as in  \cite{Hir}.
 \end{remark}

\par
A set $X\subset K^m$ is called \emph{semialgebraic} if it is a
finite union of sets of the form
\[\{x\in K^m\mid f(x)=0,\ g_1(x)\in P_{n_1},\ \ldots,\ g_k(x)\in P_{n_k}\},\]
where the functions $f$ and $g_i$ are polynomials over $K$,
$n_i>0$; we call a function semialgebraic if its graph is
semialgebraic. It is well-known that also the collection of
semialgebraic sets is closed under taking complements, finite
unions and intersections, and projections (see \cite{Mac},
\cite{Denef2}).
 \section{Subanalytic Cell Decomposition}\label{proof
 cell decomp}
The cell decomposition theorem makes use of basic sets called
(subanalytic) cells, which we define inductively. For $m,l>0$ we
write $\pi_m:K^{m+l}\to K^m$ for the linear projection on the
first $m$ variables and, for $A\subset K^{m+l}$ and $x\in
\pi_m(A)$, we write $A_x$ for the fiber $\{t\in K^{m+l}\mid
(x,t)\in A\}$.
 \\

 \begin{definition}\label{def:cell}
A cell $A\subset K$ is either a point or a set of the form
 \begin{equation}
\{t\in K\mid |\aaa|\sq_1 |t-\gamma|\sq_2 |\bbb|,\
  t-\gamma\in \lP\},
\end{equation}
with constants $n>0$, $\lambda,\gamma\in K$, $\aaa,\bbb\in
K^\times$, and $\square_i$ either $<$ or no condition.
 A cell $A\subset K^{m+1}$, $m\geq0$, is a set  of the
form
 \begin{equation}\label{Eq:cell}
 \begin{array}{ll}
\{(x,t)\in K^{m+1}\mid
 &
 x\in D,\,  |\aaa(x)|\sq_1 |t-\gamma(
 x)|\sq_2 |\bbb(x)|,\\
 &
 t-\gamma(x)\in \lP\},
  \end{array}
\end{equation}
 with $(x,t)=(x_1,\ldots,
x_m,t)$, $n>0$, $\lambda\in K$, $D=\pi_m(A)$ a cell, subanalytic
functions $\aaa,\bbb:K^m\to K^\times$, $\gamma:K^m\to K$, and
$\square_i$ either $<$ or no condition. We call $\gamma$ the
center and $\lambda P_n$ the coset of the cell $A$.
 \end{definition}
Remark that a cell $A\subset K_1^{m+1}$ is either the graph of a
subanalytic function defined on $\pi_m(A)$ (if $\lambda=0$), or,
for each $x\in \pi_m(A)$, the fiber $A_x\subset K$ contains a
nonempty open.

\par
Theorem \ref{Thm:Cell} below is a subanalytic analogue of the
$p$-adic semialgebraic cell decomposition theorem (see
\cite{Denef} and \cite{Denef2}); it is a perfect analogue of the
reformulation \cite[lemma 4]{C}. In \cite{Denef1}, an overview is
given of applications of $p$-adic cell decomposition theorems,
going from a description of local singular series to profinite
$p$-groups (see \cite{DuS}).
 \begin{theorem}[Subanalytic Cell Decomposition]\label{Thm:Cell}
Let $X\subset K^{m+1}$ be a subanalytic set, $m\geq0$, and
$f_j:X\to K$ subanalytic functions for $j=1,\ldots,r$. Then there
exists a finite partition of $X$ into subanalytic cells $A$ with
center
$\gamma:K^m\to K$ and coset $\lambda P_n$ such that for each $(x,t)\in A$ 
 \begin{equation}
 |f_j(x,t)|=
 |\ddd_j(x)|\, |(t-\gamma(x))^{a_j}\lambda^{-a_j}|^\frac{1}{n},\qquad
 \mbox{ for each } j=1,\ldots,r,
 \end{equation}
with $(x,t)=(x_1,\ldots, x_m,t)$, integers $a_j$, and
$\ddd_j:K^m\to K$  subanalytic functions. If $\lambda=0$, we
understand $a_j=0$, and we use the convention $0^0=1$.
 \end{theorem}
For the proof of Thm.~ \ref{Thm:Cell} we use techniques from model
theory. (For general notions of model theory we refer to
\cite{Hodges}.)
\par
Let $(K_1,\Ldan)$ be an elementary extension
of $(K,\Ldan)$ and let $R_1$ be its valuation ring. In view of
Thm.~\ref{Gabrielov}, we can call a set $X\subset K_1^m$
subanalytic if it is $\Ldan$-definable (with parameters from
$K_1$)
and analogously for subanalytic functions, cells, and so on.
Expressions of the form $|x|<|y|$ for $x,y\in K_1$ are
abbreviations for the corresponding $\Ldan$-formula's expressing
$|x|<|y|$ for $x,y\in K$, as in lemma 2.1 of
\cite{Denef2}\footnote{For example, if $K=\Q_p$ with $p\not=2$,
the property $|x|<|y|$ is equivalent to $y^2+\frac{x^2}{p}\in
P_2$.}. One can make similar definitions for semialgebraic subsets
of $K_1$.\footnote{This can be done using the language of
Macintyre, consisting of $+,-,\cdot,0,1,$ and the collection of
predicates $\{P_n\}$ for $n>0$. 'Semialgebraic' then means
'definable with parameters'.}
\par
We recall some definitions and results of \cite{vdDHM} in a
slightly more simple form than they were originally formulated.
\begin{definition}[\cite{vdDHM}]
Let $K_1$ and $R_1$ be as above. A set $F\subset R_1$ is called a
\emph{$R_1$-affinoid} if it is of the form
\[
\{t\in R_1\mid |t-a_1|\leq |\pi_1|\mbox{ and } |t-a_i|\geq
|\pi_i|\mbox{ for } i=2,\ldots, k \},
\]
where $a_i,\pi_i\in R_1$, $\pi_i\not=0$ and every ball $\{t\in
K\mid |t-a_i|< |\pi_i|\}$, $i=2,\ldots, k$, lies inside $\{t\in
K\mid |t-a_1|\leq |\pi_1|\}$. For such affinoid $F$, we define the
algebra $\OF$ of holomorphic functions on $F$, consisting of all
functions $f:F\to K_1$ of the form
\[
f(t) = c\cdot
g(x_1,\ldots,x_l,\frac{t-a_1}{\pi_1},\frac{\pi_2}{t-a_2},\ldots,\frac{\pi_k}{t-a_k}),
\]
where $x\in R_1^l$ and $c\in K_1$ are fixed constants and $g$ is
in $\TK$. 
\end{definition}
Proposition \ref{cor to vdDHM 5.6} below is a reformulation of
\cite[Thm. 5.6]{vdDHM} and proposition \ref{vdDHM cor 3.4} is an
immediate corollary of Proposition 3.4 and Lemma 2.10 of
\cite{vdDHM}.
\begin{proposition}[\cite{vdDHM}]\label{cor to vdDHM 5.6}
Let $K_1$ and $R_1$ be as above. Let $f:R_1\to K_1$ be a
subanalytic function, then there exists a finite partition $\PP$
of $R_1$ into semialgebraic sets such that for each $A\in P$ we
can find a $R_1$-affinoid $F$ containing $A$, holomorphic
functions $h_1,\ldots, h_r\in\OF$ and a semialgebraic function
$g:K_1^r\to K_1$ such that
$f(t)=g\left(h_1(t),\ldots,h_r(t)\right)$ for each $t\in A$.
\end{proposition}
\begin{proposition}[\cite{vdDHM}]
\label{vdDHM cor 3.4} Let $F$ be a $R_1$-affinoid and let
$f\in\OF$, then there is a rational function $r(x)\in K_1(x)$ with
no poles in $F$ such that $|f(x)|=|r(x)|$ for all $x\in F$.
\end{proposition}
\begin{theorem}[theorem B of \cite{vdDHM}]
Each subanalytic subset of $K_1$ is semialgebraic.
\end{theorem}
Next we state and prove a one-dimensional version of
Thm.~\ref{Thm:Cell}.
\begin{lemma}\label{cell decomp dim 1}
Let $K_1$ and $R_1$ be as above. Let $X\subset R_1$ be a
subanalytic set and $f:X\to K_1$ a subanalytic function. Then
there exists a finite partition $\PP$ of $X$ into cells, such that
for each cell $A\in\PP$ with center $\gamma\in K$ and coset
$\lambda P_n$
\[
|f(t)|=|\ddd|\, |(t-\gamma)^a\lambda^{-a}|^\frac{1}{n}\mbox{ for
each } t\in A,
\]
with $\ddd\in K_1$ and an integer $a$. We use again the convention
that if $\lambda=0$ then $a=0$.
\end{lemma}
\begin{proof} Replace $f$ by the function $R_1\to K_1$ by putting
$f(x)=0$ if $x\not\in X$. By proposition \ref{cor to vdDHM 5.6} we
find a finite partition $\PP_1$ of $R_1$ into semialgebraic sets
such that for each $A\in\PP_1$ there is a $R_1$-affinoid
$F\supseteq A$, holomorphic functions $h_1,\ldots,h_r\in\OF$ and a
semialgebraic function $g:K_1^r\to K_1$ such that $f(x)=g(h(x))$
for each $x\in A$, where $h=(h_1,\ldots,h_r)$. By \cite{Denef},
Corollary 6.5, we can take a finite semialgebraic partition $\PP'$
of $K_1^r$ and for each part $B$ polynomials $p_B,q_B\in
K_1[T_1,\ldots,T_r]$ and an integer $e>0$ such that $q_B(y)\not=0$
for each $y\in B$
 and
\[
|g(y)|=|\frac{p_B(y)}{q_B(y)}|^{1/e},\quad \mbox{ for each }y\in B
\]
implying that
\[
|f(x)|=|\frac{p_B\left(h(x)\right)} {q_B\left(h(x)\right)}|^{1/e}
\mbox{ for each }x\in A\cap h^{-1}(B).
\]
 Since $\OF$ is an algebra, the functions
$p_B(h)$ and $q_B(h))$ are again in $\OF$. By Theorem B of
\cite{vdDHM}, the sets $h^{-1}(B)$ are semialgebraic for each
$B\in\PP'$, moreover, the sets $A\cap h^{-1}(B)$ form a finite
semialgebraic partition of $A$. By proposition \ref{vdDHM cor
3.4}, applied to the functions $p_B(h)$ and $q_B(h)$ for each
$B\in\PP'$, and by repeating this process for each $A\in \PP$, we
find a partition $\PP_2$ of $R_1$ into finitely many semialgebraic
sets $C$ and for each $C$ an integer $e\not=0$ and a rational
function $r_C(x)\in K_1(x)$ without poles in $C$, such that
\[
|f(x)|=|r_C(x)|^{1/e}\quad \mbox{ for each }x\in C.
\]
Now apply the semialgebraic cell decomposition theorem (in the
formulation of  \cite[Lem.~4]{C}) to the parts $C\in\PP_2$ and the
(semialgebraic) functions $r_C:C\to K_1$. Refine the obtained
partition such that for each cell $A\subset C$ with coset $\lambda
P_n$, the number $n$ is a multiple of $e$, then the lemma follows.
 \end{proof}
We will use the previous lemma and a model-theoretical compactness
argument to prove the following generalization of
Thm.~\ref{Thm:Cell}.
\begin{theorem}
\label{thm:cell2} Let $(K_1,\Ldan)$ be an arbitrary elementary
extension of $(K,\Ldan)$ with valuation ring $R_1$.
 Let $X\subset K_1^{m+1}$ be subanalytic and $f_j:X\to K_1$
subanalytic functions for $j=1,\ldots,r$. Then there exists a
finite partition of $X$ into subanalytic cells $A$ with center
$\gamma:K_1^m\to K_1$ and coset $\lambda P_n$ such that for each
$(x,t)\in A$
\[
 |f_j(x,t)|= |\ddd_j(x)|\,|(t-\gamma(x))^{a_j}\lambda^{-a_j}|^\frac{1}{n},
 \]
with $(x,t)=(x_1,\ldots, x_m,t)$, integers $a_j$, and
$\ddd_j:K_1^m\to K_1$ subanalytic functions, $j=1,\ldots,r$. Here
we use the convention that $a_j=0$ if $\lambda=0$ and that
$0^0=1$.
\end{theorem}
\begin{proof} The proof goes by induction on $m\geq 0$. It is enough to prove
the theorem for $r=1$ (the theorem then follows after a
straightforward further partitioning, see for example
\cite{Denef2}). When $m=0$, we start with a subanalytic set
$X\subset K_1$ and a subanalytic function $f:X\to K_1$. Partition
$X$ into parts $X_1=X\cap R_1$ and $X_2=X\cap (K_1\setminus R_1)$.
We prove the statement for $X_i$ and $f|_{X_i}:X_i\to K_1$,
$i=1,2$. For $X_1$ and $f|_{X_1}$ we can apply Lemma \ref{cell
decomp dim 1}. For $X_2$ and $f|_{X_2}$ we define $X'_2=X_2^{-1}=
\{w\mid w^{-1}\in X_2\}$ and $f':X_2'\to K_1:w\mapsto f(w^{-1})$,
and apply Lemma \ref{cell decomp dim 1} to $X'_2$ and the function
$f'$. This way we find a partition $\PP_1$ of $X_2'$ into cells
 such that for each cell $A\in \PP_1$ with center
$\gamma\in K_1$ and coset $\lambda P_n$ and each $w\in A$ we have
$|f'(w)|=|\ddd|\,|(w-\gamma)^a\lambda^{-a})|^\frac{1}{n}$, with an
integer $a$ and $\ddd\in K_1$. Remark that the sets $A^{-1}$ form
a partition of $X_2$ and
$|f(t)|=|\ddd|\,|(t^{-1}-\gamma)^a\lambda^{-a}|^\frac{1}{n}$ for
each $t\in A^{-1}$. Apply now the semialgebraic cell decomposition
theorem for each $A\in\PP_1$ to the (semialgebraic) function
$A^{-1}\to K_1:t\mapsto \sqrt[n]{(t^{-1}-\gamma).\lambda^{-1}}$,
where $\sqrt[n]{}$ is a semialgebraic $n$-th root function defined
on $P_n$ (such a function exists by the definability of Skolem
functions).
\par
Let $X\subset K_1^{m+1}$ and $f:X\to K_1$ be subanalytic, $m>0$.
We write $(x,t)=(x_1,\ldots,x_m,t)$ and know by the previous that
for each fixed $x\in K_1^m$ we can decompose the fiber $X_x$ and
the function $t\mapsto f(x,t)$ on this fiber. We will measure the
complexity of given decompositions on which $|f(x,\cdot)|$ has a
nice description and see that this must be uniformly bounded when
$x$ varies.
 \\
To measure the complexity of such decompositions, we define a
countable set $\SSS=\{\lP\mid \lambda\in K,\ n>0\}\times\Z\times
\{<, \emptyset\}^2$ and $\SSS'=(K_1^\times)^2\times K_1^2$. To
each $d=(\lP,a,\sq_1,\sq_2)$ in $\SSS$ and
$\xi=(\xi_1,\xi_2,\xi_3,\xi_4)\in \SSS'$ we associate a set
$Dom_{(d,\xi)}$ as follows
 \[
Dom_{(d,\xi)} = \{t\in K_1\mid |\xi_1|\sq_1 |t-\xi_3|\sq_2
|\xi_2|,\ t-\xi_3\in\lP\}
 \]
The set $Dom_{(d,\xi)}$ is either empty or a cell and is
independent of $\xi_4$ and $a$. For fixed $k>0$ and tuple
$d=(d_1,\ldots,d_k)\in \SSS^k$, let $\varphi_{(d,k)}(x,\xi)$ be a
$\Ldan$-formula in the free variables $x=(x_1,\ldots,x_m)$ and
$\xi=(\xi_1,\ldots,\xi_k)$, with
$\xi_i=(\xi_{i1},\xi_{i2},\xi_{i3},\xi_{i4})$, such that the tuple
$(x,\xi)\in K_1^{m+4k}$ satisfies $\varphi_{(d,k)}$  if and only
if the following are true:
\begin{itemize}
\item[(i)] $x\in\pi_m(X)$ and $\xi\in (\SSS')^k$,\\
\item[(ii)] the collection  of the sets $Dom_{(d_i,\xi_i)}$ for $i=1,\ldots,k$ forms
 a partition of the fiber $X_x=\{t\in K_1\mid (x,t)\in X\}\subset
K_1$,\\
\item[(iii)] $|\xi_{i4}|\,|(t-\xi_{i3})^{a_i}\lambda_i^{-a_i}|^\frac{1}{n_i}=\abs{f(x,t)}$ for
each  $t\in Dom_{(d_i,\xi_i)}$ and each $i=1,\ldots,k$.
 \end{itemize}
Now we define for each $k>0$ and $d \in \SSS^k$ the set
\[
B_d= \{x\in K_1^m\mid 
\exists \xi\quad \varphi_d(x,\xi)\}.
\]
Each set $B_d$ is subanalytic and the (countable) collection
$\{B_d\}_{k,d}$ covers $\pi_m(X)$, because each $x\in\pi_m(X)$ is
in some $B_d$ by the induction. Since $K_1$ is an arbitrary
elementary extension of $K$, finitely many sets of the form $B_d$
must already cover $\pi_m(X)$ by model-theoretical compactness.
Consequently, we can take subanalytic sets $D_1,\ldots,D_s$ such
that $\{D_i\}$ forms a partition of $\pi_m(X)$ and each $D_i$ is
contained in a set $B_d$ for some $k>0$ and $k$-tuple $d$. For
each $i=1,\ldots,s$, fix such a $d$ with $D_i\subset B_d$, and let
$\Gamma_i$ be the subanalytic set
\[
\Gamma_i=\{(x,\xi)\in D_i\times (\SSS')^k\mid 
\varphi_d(x,\xi)\}.
\]
Then $\pi_m(\Gamma_i)=D_i$ by construction ($\pi_m$ is the
projection on the $x$-coordinates). By theorem 3.6 \cite{DvdD} on
definable Skolem functions, there is a subanalytic function
$D_i\to K_1^{4k}$ associating to $x$ a tuple $\xi(x)\in (\SSS')^k$
such that $(x,\xi(x))\in\Gamma_i$ for each $x\in D_i$. The theorem
follows now by partitioning further with respect to the
$x$-variables and using the induction hypothesis.
 \end{proof}
\section{Classification of Subanalytic Sets}\label{section:classification}
For $X\subset K^m$ subanalytic and nonempty, the dimension
$\dim(X)$ of $X$ is defined as the biggest integer $n$ such that
there is a $K$-linear map $\pi:K^m\to K^n$ and a nonempty
$U\subset \pi(X)$, open in $K^n$ with respect to the norm topology
(for alternative definitions, see \cite{DvdD}).
\begin{theorem}\label{prop:bijectie}
For any subanalytic set $X\subset K^m$ and subanalytic functions
$f_i:X\to K$, $i=1,\ldots,r$, there is a semialgebraic set $Y$, a
subanalytic bijection $F:X\to Y$ and there are semialgebraic maps
$g_i:Y\to K$ such that
\[|g_i( F(x))|=|f_i(x)|\qquad  \mbox{for each }x\in X.
\]
\end{theorem}
\begin{proof} 
We will give a proof by induction on $m$. Suppose that $X\subset
K^{m+1}$ is subanalytic and that $f_i:X\to K$ are subanalytic
functions, $m\geq0$. Apply the cell decomposition theorem to $X$
and the functions $f_i$ to obtain a finite partition $\PP$ of $X$.
For $A\in\PP$ and $(x,t)\in A$, suppose that
$|f_i(x,t)|=|\ddd_i(x)|\,|(t-\gamma(x))^{a_i}\lambda^{-a_i}|^\frac{1}{n}$,
$i=1,\ldots,r$, and suppose that $A$ is a cell of the form
 \[
 \begin{array}{ll}
 \{(x,t)\in K^{m+1} \mid
 &
 x\in D,\, |\aaa(x)|\sq_1
 |t-\gamma(x)|\sq_2 |\bbb(x)|,\\
 &
 t-\gamma(x)\in\lP \},
 \end{array}
 \]
 like in Eq~(\ref{Eq:cell}).
 After the translation $(x,t)\mapsto (x,t-\gamma(x))$ we
may suppose that $\gamma$ is zero on $D$. Apply the induction
hypotheses to the sets $D$ and the subanalytic functions
$\aaa,\bbb$, and $\ddd_i$. Repeating this process for every $A\in
\PP$, the proposition follows after taking appropriate disjoint
unions inside $K^m$ of the occurring semialgebraic sets.
 \end{proof}
We prove the following generalization of Theorem
\ref{classification}.
\begin{theorem}[Classification of subanalytic sets]
\label{thm:classification} Let $X\subset K^m$ and $Y\subset K^n$
be infinite subanalytic sets, then there exists a subanalytic
bijection $X\to Y$ if and only if $\dim(X)=\dim(Y)$.
\end{theorem}
\begin{proof} By theorem \ref{prop:bijectie}
there are subanalytic bijections  $X\to X'$ and $Y\to Y'$ with
$X'$ and $Y'$ semialgebraic, but then there exists a semialgebraic
bijection $X'\to Y'$ if and only if $\dim(X')=\dim(Y')$  by
theorem 2 of \cite{C}, see chapter \ref{chap:classification}.
Since the dimension of a subanalytic set is invariant under
subanalytic bijections (see \cite{DvdD}), the theorem follows.
 \end{proof}
\section{Parametrized Analytic Integrals}\label{section integrals}
We show that certain algebra's of functions from $\Q_p^m$ to the
rational numbers $\Q$  are closed under $p$-adic integration.
These functions are called simple $q$-exponential functions and
they come up naturally when one calculates parametrized $p$-adic
integrals like for example Eq.~(\ref{integral}).
\par
For $x=(x_1,\ldots,x_m)$ a $m$-tuple of variables, we will write
$|dx|$ to denote the Haar measure on $K^m$, so  normalized that
$R^m$ has measure $1$.
\begin{definition}\label{basic algebra's} Let $\OOL$ be the algebra
of (subanalytic) simple $q$-exponential functions, to be precise,
$\OOL$ is the $\Q$-algebra generated by the functions $K^m\to\Q$
for $m\geq 0$ of the form $x\mapsto v(h(x))$ and $x\mapsto
|h'(x)|$ where $h:K^m\to K^\times$ and $h':K^m\to K$ are
subanalytic functions.
\par
To any function $f:K^{m+n}\to \Q$ in $\OOL$, $m,n\geq 0$, we
associate a function $I_m(f):K^m\to \Q$ by putting
 \begin{equation}\label{I_l}
I_m(f)(\lambda)= \int\limits_{K^n}f(\lambda,x)|dx|
 \end{equation}
if the function $x\mapsto f(\lambda,x)$ is absolutely integrable
and $I_m(f)(\lambda)=0$ otherwise.
\end{definition}
 \begin{theorem}[Basic Theorem on $p$-adic Analytic Integrals]\label{thm:basic}
The algebra $\OOL$ is closed under $p$-adic integration, to be
precise, for any function $f\in\OOL$ in $m+n$ variables, the
function $I_m(f):K^m\to \Q$ is in $\OOL$.
 \end{theorem}
 \begin{proof}
By induction it is enough to prove that for a function
$f:K^{m+1}\to \Q$ in $\OOL$ in the variables
$(\lambda_1,\ldots,\lambda_m,t)$ the function $I_m(f)$ is in
$\OOL$. Suppose that $f$ is a $\Q$-linear combination of products
of functions $|f_i|$ and $v(g_j)$, $i=1,\ldots,r$, $j=1,\ldots,s$
where $f_i$ and $g_j$ are subanalytic functions $K^{m+1}\to K$ and
$g_j(\lambda,t)\not=0$. Applying the cell decomposition theorem to
$K^{m+1}$ and the functions $f_i$ and $g_j$, we obtain a partition
$\cP$ of $K^{m+1}$ into cells such that $I_m(f)(\lambda)$ is a sum
of integrals over $A_\lambda=\{t\mid(\lambda,t)\in A\}$ for each
cell $A\in\cP$, where the integrands on these pieces $A_\lambda$
have a very simple form. More precisely, on each piece $A_\lambda$
the integrand is a $\Q$-linear combination of functions of the
form
\begin{equation}\label{int}
\ddd(\lambda) |(t-\gamma)^a\mu^{-a}|^\frac{1}{n}
v(t-\gamma(\lambda))^l
\end{equation}
where $A$ is a cell with center $\gamma:K^m\to K$ and coset $\mu
P_n$, and with integers $a$ and $0\leq l$, and a function
$\ddd:K^m\to\Q$ in $\OOL$. Partitioning further, we may suppose
that $\ddd(\lambda)\not=0$ on $\pi_m(A)$. Regroup all such terms
where the same exponents $a$ and $l$ appear, possibly by replacing
the functions $\ddd(\lambda)$ by other functions in $\OOL$.
Whether or not such a $\Q$-linear combination is absolutely
integrable only depends on the integers $a$ and $l$ occurring in
each of the terms as in Eq.~(\ref{int}) and on the symbols $\sq_i$
occurring in the description of the cell $A$, being $<$ or no
condition, and hence does not depend on the particular choice of
$\lambda$ in the projection $\pi_m(A)$. By consequence, we may
suppose that the integrand is a single term of the form like in
Eq.~(\ref{int}) and that this term is absolutely integrable. It
suffices to show that the integral
\begin{equation}\label{intint}
\ddd(\lambda)\int\limits_{t\in A_\lambda}
|(t-\gamma(\lambda))^a\mu^{-a}|^\frac{1}{n}
v(t-\gamma(\lambda))^l|dt|
\end{equation}
is in $\OOL$. Write $u=t-\gamma(\lambda)$; since $A$ is a cell
with center $\gamma$ and coset $\mu P_n$, the set $A$ is  of the
form
\[
A= \{(\lambda,u)\in K^{m+1}\mid   \lambda\in D, \
|\aaa(\lambda)|\sq_1 |u|\sq_2 |\bbb(\lambda)|,\
  u\in \mu P_n\},
\]
with $\sq_i$ either $<$ or no condition, $D$ a cell, and
$\aaa,\bbb:K^m\to K^\times$ subanalytic functions. Taking into
account that the integral (\ref{intint}) is, by supposition,
integrable, only a few possibilities can occur (with respect to
the sign of the integers $a$ and $l$, the conditions $\sq_i$ being
$<$ or no condition, and $\mu$ being zero or nonzero). If $\mu=0$
the set $A_\lambda$ is a point for each $\lambda\in D$, thus the
statement is clear. Suppose $\mu\not=0$. In case that both $\sq_1$
and $\sq_2$ represent no condition, the integrand has to be zero
by the supposition of integrability, and the above integral
trivially is in $\OOL$. We suppose from now on that $\sq_1$ is
$<$; the other cases can be treated similarly. The integral of
Eq.~(\ref{intint}) can be written as
\[
\ddd(\lambda)\cdot\int\limits_{t\in A_\lambda}
|(t-\gamma(\lambda))^a\mu^{-a}|^\frac{1}{n}
v(t-\gamma(\lambda))^l|dt|
\]
\begin{eqnarray*}
& = &\ddd(\lambda)\sum_k (q^{-ak}|\mu^{-a}|)^\frac{1}{n}
k^{l}\cdot\mathrm{Measure}\{t\in
A_\lambda\mid v(t)=k\}\\
&= &\ \epsilon\ddd(\lambda) \sum_k
(q^{-ak}|\mu^{-a}|)^\frac{1}{n}q^{-k} k^{l}
\end{eqnarray*}
for some $\epsilon\in\Q$, where the summation is over those
integers $k\equiv v(\mu) \bmod{n}$ satisfying
 \[
 |\aaa(\lambda)|<
q^{-k}\sq_2|\bbb(\lambda)|.
 \]
 For fixed residue classes
\[
v(\aaa(\lambda))\pmod{n}\quad\mbox{ and  }\quad
v(\bbb(\lambda))\pmod{n},
 \]
this sum is equal to a $\Q$-linear combination of products of the
functions $\ddd$, $|\aaa|$, $v(\aaa)$, $|\bbb|$ and $v(\bbb)$.
Since the characteristic function of sets of the form
\[
\{\lambda\mid v(\aaa(\lambda))\equiv a_1\bmod{n}\mbox{ and }
v(\bbb(\lambda))\equiv a_2\bmod{n}\}
\]
is in $\OOL$ for any integers $a_i$, the integral \ref{intint} is
in $\OOL$ as was to be shown.
\end{proof}
As a corollary we will formulate another version of the basic
integration theorem, conjectured in \cite{Denef1} in the remark
following theorem~2.6.
\begin{definition}
A set $A\subset \Z^n\times \Q_p^m$ is called \emph{simple} if
\[\{(\lambda,x)\in K^{n+m}\mid (v(\lambda_1),\ldots,v(\lambda_n),x)\in
A,\quad \lambda_i\not=0\}\]
 is a subanalytic set. A function $h:\Z^n\times \Q_p^m\to \Z$ is called \emph{simple} if
its graph is simple. Let $\OOLp$ be the $\Q$-algebra generated by
all simple functions and all functions of the form $q^h$ where $h$
is a simple function.
\par
For a function $f:\Z^{k+l}\times K^{m+n}\to\Q$ in $\OOLp$ we
define $I_{k,m}(f):\Z^{k}\times K^{m}\to\Q$ as
\[
\sum_{z'\in \Z^{l}}\int\limits_{K^n}f(z,z',\lambda,x)|dx|
\]
if the function $(z',x)\mapsto f(z,z',\lambda,x)$ is absolutely
integrable with respect to the Haar measure on $K^n$ and the
discrete measure on $\Z^l$ and we define $I_{k,m}(f)(z,\lambda)=0$
otherwise.
\end{definition}
\begin{theorem}\label{Thm:conjecture}
The algebra $\OOLp$ is closed under integration, to be precise,
for each $f:\Z^{k+l}\times K^{m+n}\to \Q$ in $\OOLp$, the function
$I_{k,m}(f):\Z^k\times K^l\to\Q$ is in $\OOLp$.
\end{theorem}
\begin{proof}
It is enough to prove that for a function $f:\Z^k\times K^{m}\to
\Q$ in $\OOLp$ in the variables $(z_1,\ldots,z_k,x_1,\ldots,x_m)$
the function obtained by integrating $x_m$ out, resp. by
integrating $z_k$ out, is in $\OOLp$. To $f:\Z^k\times K^{m}\to\Q$
we can associate a function $g:K^{k+m}\to \Q$ by replacing the
variables $z$ running over $\Z^k$ by variables $\lambda$ running
over $K^k$ in such a way that
$g(\lambda,x)=f(v(\lambda_1),\ldots,v(\lambda_k),x)$ for each
$\lambda\in (K^\times)^k$ and $g(\lambda,x)=0$ if one of the
$\lambda_i$ is zero. By the definitions, $g$ is a  simple
$q$-exponential function and the given discrete sum corresponds to
a $p$-adic integral of the function $g$. If we integrate the
$x_m$-variable out, then we get a function
$I(\lambda,x_1,\ldots,x_{m-1})$ which is in $\OOL$ by theorem
\ref{thm:basic}. This function only depends on
$(v(\lambda_1),\ldots,v(\lambda_k),x)$ and thus corresponds to a
function in $\OOLp$. The same arguments hold if we integrate
$\lambda_k$ out and thus the prove is complete.
 \end{proof}
\section*{Some Remarks} Many of the
results of \cite{DvdD}, as well of \cite{vdDHM} are formulated for
$\Q_p$ and not for finite field extensions of $\Q_p$,
nevertheless, all results referred to in this paper, also hold for
finite field extensions of $\Q_p$ (see the remark in (3.31) of
\cite{DvdD}). I want to remark that many of the applications in
\cite{Denef1} as well as many results of \cite{Denef3},
\cite{vdDHM} and \cite{DvdD} can be obtained immediately in the
subanalytic context using the cell decomposition and the
integration theorems of this paper, for example: Cor.~1.8.2.~of
\cite{Denef3} on local singular series; the $p$-adic Lojasiewicz
inequalities (3.37), the subanalytic selection theorem (3.6),
Thm.~(3.2) and the stratification theorem.~(3.29) of \cite{DvdD};
Thm.~3.1 of \cite{Denef1}.
\par
Very often, Theorem \ref{Thm:Cell} acts as a $p$-adic analogue of
the preparation theorem \cite{LR} for real subanalytic functions.
Theorem \ref{Thm:Cell} can also be seen as a $p$-adic analogue of
the cell decomposition theorem for o-minimal structures on $\R$
(see e.g.~\cite{vdD}).

\section*{History and acknowledgment} I would like to thank
Denef for pointing out to me his conjecture on $p$-adic
subanalytic integration and for stimulating conversations on the
subject. I also thank D.~Haskell and M.-H.~Mourgues for many
interesting discussions.
\par
Partial results towards a subanalytic cell decomposition theorem
have been obtained in \cite{Liu1},  \cite{Mourgues}, \cite{Mil}
and \cite{Wilc}, where the authors of \cite{Liu1} and \cite{Mil}
were supervised by L.~Lipshitz and the author of \cite{Wilc} by
A.~Macintyre. In \cite{Mourgues} and \cite{Wilc}, a partitioning
of arbitrary subanalytic sets into cells is obtained, but not the
description of the norm of subanalytic functions on these cells.
However, the description of the valuation of the subanalytic
functions in theorem \ref{Thm:Cell} is essential for the
applications.

\clearemptydoublepage

\chapter[Decay Rates]{Multi-variate Igusa theory:\\
 Decay rates of $p$-adic exponential sums}\label{chap:decay}
\begin{abstract}
\footnote{This chapter corresponds to \cite{Cexp}.} We investigate
the asymptotic behaviour of the multi-dimensional analogon of the
one-dimensional  exponential sum
$$
\sum_{x\in (\Z_{p}/p^{m})^{n}}\exp(2\pi i \frac{f(x)}{p^{m}})
$$
as $m$ tends to infinity, where $f$ is a polynomial over $\Z_p$.
To be precise, we replace $f$ by a polynomial map $(f_1,\ldots,
f_r)$ and $1/p^m$ by a tuple of $p$-adic numbers. J.~Igusa used
resolution of singularities to study the asymptotics of the one
dimensional sum. The method of Igusa, based on resolution of
singularities, is not suitable to be generalized to the general
multivariate case. To overcome this difficulty, we use integration
techniques based on $p$-adic cell decomposition. All results of
the paper also hold for analytic maps and for finite field
extensions of $\Q_p$.
\end{abstract}
\section{Introduction}
Let $f=(f_1,\ldots, f_r)$ be a dominant polynomial map in the
variables $x=(x_1,\ldots,x_n)$, given by polynomials
$f_i\in\Z_p[x]$, and let $\psi$ be the standard additive
character\footnote{the definition is given below; we postpone
several definitions to later sections.} on $\Q_p$. For $y\in
\Q_p^r$, we consider the exponential integral
$$
E(y)=\int_{\Z_{p}^{n}}\psi(\langle y,f(x)\rangle)|dx|,
$$
where $\langle,\rangle$ denotes the inner product on $K^r$, and
$|dx|$ stand for the Haar measure on $\Q_p^n$, normalized so that
$\Z_{p}^{n}$ has measure 1. We obtain upper bounds for $|E(y)|$ in
terms of the $p$-adic norm $|y|=\max_i|y_i|$ of $y$ and we
investigate the asymptotic behaviour of $E(y)$ as $|y|$ tends to
$\infty$. This integral (which is actually a generalized Gauss
sum, see below) converges to zero whenever $|y|$ goes to infinity;
this is not unexpected because of more and more cancellation of
the integrand due to the character $\psi$. We obtain the following
qualitative decay rate of $E(y)$ when $|y|$ goes to $\infty$.
\begin{theorem}\label{thm:decayexp:A}
There exist real numbers $\alpha<0$ and $c>0$ such that
$$
|E(y)|<c\min \{|y|^{\alpha},1\},\qquad\mbox{ for all } y\in
\Q_p^r.
$$
\end{theorem}
In his book \cite{Igusa3}, J.~Igusa proves theorem
\ref{thm:decayexp:A} in the case that $r=1$ (the univariate case)
and assuming that $f_1$ is homogeneous. Loxton \cite{Lox}, Denef
and Sperber, and many others also study the case $r=1$. B.~Lichtin
studies $E(y)$ for $r=2$ and in some special cases for $r>2$, see
\cite{Lichtin} and \cite{Lichtin2}. Lichtin and Igusa relate
$\alpha$ to fine geometrical invariants of the locus of $f=0$ like
the numerical data of an embedded resolution of singularities of
$f$.\\
 To study this exponential sum for general $r$
is part of the open problem proposed by Igusa on page 39 of
\cite{Igusa3}; he describes the problem of generalizing the whole
theory of local zeta functions to the multivariate setting.
\par
If $r=1$ and if we put $y=up^{-m}$, with $u\in \Z_{p}^{\times}$,
the above formula reduces to
$$
E(y)=p^{-nm}\sum_{x\in (\Z_{p}/p^{m})^{n}}\exp(2\pi
i\frac{uf(x)}{p^{m}})
$$
which is a classical exponential sum modulo $p^{m}$. It is
surprising that the (normalized) sum converges to zero with decay
at least a positive power of $p^{-m}$.
\par
Let $F(z)$ be the local singular series of $f$, precisely,
\[
F(z) = \int_{f^{-1}(z)\cap \Z_p^n}|\frac{dx}{df}|,
\]
whenever $z$ is a regular value of $f$. The functions $E$ and $F$
were introduced by A.~Weil in \cite{Weil}. Igusa  showed that $E$
is the Fourier transform of $F$.  Denef  proves in \cite{Denef1}
that $F$ is a simple $p$-exponential function (see section
\ref{sec:decay:cell:decomp} for the definition of this class of
functions); this is the multidimensional analogue of the
asymptotic expansions of local singular series given by Igusa in
\cite{Igusa1}, \cite{Igusa2} and \cite{Igusa3} when $r=1$. Igusa
\cite{Igusa3} relates the asymptotic properties of $E$ and $F$ in
a very precise way. We obtain the following upper bounds for the
norm of the Fourier transform of integrable simple $p$-exponential
functions.
\begin{theorem}\label{thm:decaysimple:A}
Let $G:\Z_p^r\to \Q$ be a simple $p$-exponential function in the
variables $z=(z_1,\ldots,z_r)$, which is absolutely integrable
over $\Z_p^r$. For $y\in \Q_p^r$, put
\[
G^*(y)=\int_{z\in \Z_p^r}G(z)\psi(\langle y,z\rangle)|dz|,
\]
then there exist real numbers $\alpha<0$ and $c>0$ such that
\[
|G^*(y)|< c \min\{\abs{(y)}^{\alpha},1\},\qquad\mbox{ for all
}y\in \Q_p^r.
\]
\end{theorem}
We will first prove theorem \ref{thm:decaysimple:A} and derive
theorem \ref{thm:decayexp:A} from it.
\par
All results of the paper also hold for finite field extensions of
$\Q_p$, and also when $f$ is replaced by a dominant restricted
analytic map.
\par
The proofs are based on $p$-adic cell decomposition; this is
proven by Denef \cite{Denef1} for polynomial maps, and for
analytic maps by the author \cite{Ccell}, see chapter
\ref{chap:cell} and below. Cell decomposition yields a very
precise qualitative, piecewise description of simple
$p$-exponential functions in a way which is very suitable for
integration, see section \ref{sec:decay:cell:decomp}.
\par
A second key result is the fact that local singular series are
simple $p$-exponential functions, proven by Denef \cite{Denef1} in
the semialgebraic case and the author in the analytic case, see
theorem \ref{prop:local:pol} and \ref{prop:local:an} below.
 \par
In the case $r=1$ and assuming $f_1\in\Z[x]$ to be homogeneous, it
is known that the decay rate of $E$ can be bounded by means of the
numerical data of an embedded resolution of $f$. Namely, if
$(N_{i},\nu_{i})$ are the data of such a resolution, and
$\sigma<0$ is a real number satisfying $\sigma>-\nu_{i}/N_{i}$ for
each couple $(N_{i},\nu_{i})$ different from $(1,1)$, then we can
find for each $p$-adic completion $\Q_p$ of $\Q$ a real constant
$c(\Q_p)$ such that $|E(y,\Q_p)|\leq c(\Q_p)|z|^{\sigma}$. Igusa
\cite{Igusa3} conjectured that this constant $c(\Q_p)$ can be
chosen independently from $p$. Although partial results towards
this conjecture have been proven by Igusa, Sperber and Denef, and
others, the general case remains open.
\par
I thank J.~Denef and B.~Lichtin for pointing my attention to these
open problems on exponential sums, and for interesting
conversations on this problem.
 \section{Notation and terminology}
Let $K$ be a $p$-adic field (i.e. $[K:\Q_p]$ is finite). Write $R$
for the valuation ring of $K$, $M$ for its maximal ideal, $\pi_0$
for a generator of $M$, $k$ for the residue field $R/M$ and $q$
for the cardinality of $k$. For $x\in K$ let
$v(x)\in\Z\cup\{\infty\}$ denote the $p$-adic valuation of $x$ and
$|x|=q^{-v(x)}$ the $p$-adic norm. We write $P_m$ for the
collection of $m$-th powers in $K^\times=K\setminus\{0\}$, $m>0$.
A Schwartz-Bruhat functions on $K^n$ or $R^n$ is a locally
constant function to $\Q$ with compact support. We denote by
$\psi$ the additive character $K\to \C^\times:x\mapsto \exp( 2\pi
i \mathrm{Tr}_{K/\Q_p}(x))$, where $\mathrm{Tr}_{K/\Q_p}$ is the
trace map. If we put $\langle x,y\rangle=x_1y_1+\ldots+x_ny_n$ for
$x,y\in K^n$, then the dual group of $K^n$ can be identified with
$K^n$ by $y\in K^n\mapsto \psi_y$ with
\[
\psi_y:K^n\to \C^\times:x\mapsto \psi(\langle x,y\rangle).
\]
\par
A restricted analytic function $R^n\to K$ is an analytic function,
given by a single Tate series over $K$ in $n$ variables (by
definition, these are the series which converge on $R^n$). A
subset of $K^n$ is called \emph{semialgebraic} if it can be
obtained by taking finite unions, finite intersections,
complements and linear projections of the zero locus of
polynomials in $K^{n+e}$, $e\geq0$. A subset of $K^n$ is called
\emph{subanalytic} if it can be obtained by taking finite unions,
finite intersections, complements and linear projections of
semialgebraic subsets of $K^{n+e}$ and zerosets of restricted
analytic functions inside $R^{n+e'}\subset K^{n+e'}$ for some
numbers $e$ and $e'$ (see chapter \ref{chap:cell} for equivalent
definitions). A function $f:X\subset K^{m}\rightarrow K^{n}$ is
called semialgebraic (resp.~subanalytic) if its graph is a
semialgebraic (resp.~subanalytic) set. Remark that each
semialgebraic and subanalytic function is piecewise analytic (in
the sense of \cite{Bour}), see \cite{DvdD}. We refer to
\cite{DvdD}, \cite{Ccell}, \cite{Denef1}, and chapter
\ref{chap:cell} for the basic theory of subanalytic sets.
\section{Cell decomposition and $p$-adic integration}\label{sec:decay:cell:decomp}
Cell decomposition is a powerful tool to describe piecewise
several kinds of $p$-adic maps, for example:
\begin{itemize}
\item
polynomial maps $K^n\rightarrow K^{r}$;
\item restricted analytic maps $R^n\rightarrow K^{r}$;
\item local singular series $K^n\to\Q$ (both for polynomial mappings as for restricted analytic mappings);
\item simple $q$-exponential functions $K^n\to\Q$ and also simple
analytic $q$-exponential functions.
\end{itemize}
It allows us to partition the space $K^n$ into $p$-adic manifolds
of a simple form, called cells, and to obtain on each of these
pieces a nice description of the way the function depends on a
specific special variable $t$ (see proposition
\ref{prop:descrip:simple}).
\par
Cells are defined by induction on the number of variables.
 \begin{definition}\label{def::cell}
A cell $A\subset K$ is either a point or a set of the form
 \begin{equation}
\{t\in K\mid |\aaa|\sq_1 |t-c|\sq_2 |\bbb|,\
  t-c\in \lP\},
\end{equation}
with constants $n>0$, $\lambda,c\in K$, $\aaa,\bbb\in K^\times$,
and $\square_i$ either $<$ or no condition.
 A semialgebraic (resp.~subanalytic) cell $A\subset K^{m+1}$, $m\geq0$, is a set  of the
form
 \begin{equation}\label{Eq:cell:decay}
 \begin{array}{ll}
\{(x,t)\in K^{m+1}\mid
 &
 x\in D, \  |\aaa(x)|\sq_1 |t-c(
 x)|\sq_2 |\bbb(x)|,\\
 &
  t-c(x)\in \lP\},
  \end{array}
\end{equation}
 with $(x,t)=(x_1,\ldots,
x_m,t)$, $n>0$, $\lambda\in K$, $D=\pi_m(A)$ a cell, semialgebraic
(resp.~subanalytic) functions $\aaa,\bbb:K^m\to K^\times$, and
$c:K^m\to K$, and $\square_i$ either $<$ or no condition, such
that the functions $\aaa,\bbb$, and $c$ are analytic on $D$. We
call $c$ the center of the cell $A$ and $\lambda P_n$ the coset of
the cell.
 \end{definition}
Remark that a cell  $A\subset K^{m+1}$ is a $p$-adic analytic
manifold, and it is either the graph of an analytic function
defined on the projection $\pi_m(A)\subset K^{m}$ (namely if
$\lambda=0$), or, for each $x\in \pi_m(A)$, the fiber $A_x\subset
K$ contains a nonempty open (if $\lambda\not=0$).
\begin{theorem}[$p$-adic cell decomposition]\label{thm:CellDecomp}
Let $X\subset K^{m+1}$ be a semi\-algebraic (resp.~subanalytic)
set and $f_j:X\to K$ semialgebraic (resp.~subanalytic) functions
for $j=1,\ldots,r$. Then there exists a finite partition of $X$
into semialgebraic (resp.~subanalytic) cells $A$ with analytic
center $c:K^m\to K$ and coset $\lambda P_n$ such that
 \begin{equation}
 |f_j(x,t)|=
 |\ddd_j(x)|\cdot|(t-c(x))^{a_j}\lambda^{-a_j}|^\frac{1}{n},\quad
 \mbox{ for each }(x,t)\in A,
 \end{equation}
with $(x,t)=(x_1,\ldots, x_m,t)$, integers $a_j$, and
$\ddd_j:K^m\to K$  semialgebraic (resp.~subanalytic) functions,
analytic on $\pi_m(A)$ for $j=1,\ldots,r$. If $\lambda=0$, we
understand $a_j=0$, and we use the convention $0^0=1$.
 \end{theorem}
One might think of the numbers $a_{j},n$ for each cell as in
equation (\ref{Eq:cell:decay}) as the 'numerical data' of a
specific decomposition.
\par
Denef (\cite{Denef}, \cite{Denef2}) used Macintyre's elimination
of quantifiers and techniques of P.~Cohen to prove the
semialgebraic part of Thm.~\ref{thm:CellDecomp} and the author
(\cite{Ccell}, chapter \ref{chap:cell}) used rigid $p$-adic
analysis and the theory of subanalytic $p$-adic sets to prove the
subanalytic part.
\par
We explain how cell decomposition can be used to calculate
$p$-adic integrals and to prove a basic theorem on $p$-adic
integrals.
\begin{definition}\label{def:basic:algebra}
We define $\OOs$, the algebra of simple $q$-exponential functions,
as the $\Q$-algebra generated by the functions $K^m\to\Q$ of the
form $x\mapsto v(f(x))$ and $x\mapsto |g(x)|$ where $f:K^m\to
K^\times$ and $g:K^m\to K$ are semialgebraic. Similarly, we define
the algebra $\OOK$ of simple analytic $q$-exponential functions as
the $\Q$-algebra generated by the functions $K^m\to\Q$ of the form
$v(f(x))$ and $|g(x)|$ where $f:K^m\to K^\times$ and $g:K^m\to K$
are subanalytic. To any function $f:K^{m+n}\to \Q$ in $\OOs$
(resp. in $\OOK$) we associate a function $I_m(f):K^m\to \Q$ by
putting
 \begin{equation}\label{I_l:decay}
I_m(f)(\lambda)= \int\limits_{K^n}f(\lambda,x)|dx|
 \end{equation}
if the function $x\mapsto f(\lambda,x)$ is absolutely integrable
and $I_m(f)(\lambda)=0$ otherwise.
\end{definition}
Cell decomposition yields the following piecewise description of
simple $q$-exponential functions.
\begin{prop}\label{prop:descrip:simple}
Let $G:K^n\to \Q$ be a simple $q$-exponential function in the
variables $(x_1,\ldots, x_{n-1},t)$. Then there exists a finite
partition of $K^n$ into cells $A$ with center $c$ and coset
$\lambda P_m$, such that the restriction $G|_A$ is a finite sum of
functions of the form
\[
|(t-c(x))^a\lambda^{-a}|^\frac{1}{m}v(t-c(x))^{s}h(x),
\]
where $h:K^{n-1}\to\Q$ is a simple $q$-exponential function which
is nowhere zero, and $a$ and $s\geq 0$ integers. If $\lambda=0$,
we understand $a=0$, and we use the convention $0^0=1$.
\par
The analogous statement holds also for simple analytic
$q$-exponential functions.
\end{prop}
\begin{proof}
By the definition of simple $q$-exponential functions, $G$ is a
finite sum of  products of functions of the form $v(g(x,t))$ and
$|g'(x,t)|$, where $g:K^n\to K^\times$ and $g':K^n\to K$ are
semialgebraic functions. Now apply the cell decomposition theorem
to all such functions $g$ and $g'$ occurring in such a description
of $F$. This way, one finds a finite partition of $K^n$ into cells
$A$ with center $c$ and coset $\lambda P_m$, such that the
restriction $G|_A$ is a finite sum of functions of the form
\[
|(t-c(x))^{a}\lambda^{-a}|^\frac{1}{m}v(t-c(x))^{s}h(x),
\]
where $h:K^{n-1}\to\Q$ is a simple $q$-exponential function, and
$a$ and $s\geq 0$ are integers (where we assume $a=0$ if
$\lambda=0$). By partitioning further with respect to the
$x$-space, one can obtain that the functions $h$ are nowhere zero.
(For this, one uses also the easy fact that the characteristic
function of a semialgebraic set is simple $q$-exponential.)
\end{proof}
The description of proposition \ref{prop:descrip:simple} can be
used to integrate one variable at a time out and obtain the
following theorem on $p$-adic integration.
 \begin{theorem}[Basic theorem on $p$-adic integrals]\label{thm:basic:int}
The algebra's $\OOs$ and $\OOK$ are closed under $p$-adic
integration. To be precise, for a function $f:K^{m+n}\to\Q$,
$m,n\geq0$ in $\OOs$, the function $I_m(f):K^m\to \Q$ is in $\OOs$
and similarly for $f\in\OOK$.
 \end{theorem}
The semialgebraic part of theorem \ref{thm:basic:int} is proven by
Denef \cite{Denef1}, and the subanalytic part by the author
\cite{Ccell} (see chapter \ref{chap:cell}).
\par
\section{Local singular series}\label{sec:local}
From now on, we fix a dominant\footnote{The map $f$ is dominant if
its image is not contained in the zero locus of a polynomial}
polynomial map $f=(f_1,\ldots, f_r)$, with $f_i\in K[x]$, in the
variables $x=(x_1,\ldots,x_n)$. Let $\phi:K^n\to \Q$ be a
Schwartz-Bruhat function. If $z$ is a regular value of the
polynomial mapping $f$, i.e. if there is an $x$ with $f(x)=z$ and
if the rank of the Jacobian matrix of $f$ is maximal in each point
$x$ satisfying $f(x)=z$, we can define a value $F_\phi(z)$
as\footnote{For dx/df we can take any $n-r$-form $\omega$ such
that $\omega\wedge df=dx$; the integral is independent of the
choice of $\omega$.}
$$
F_\phi(z)=\int_{f(x)=z}\phi(x)|\frac{dx}{df}|.
$$
This formula defines a function $F_\phi$ on the set of regular
values $z\in K^{r}$ of $f$, and we extend $F_\phi$ by putting
$F_\phi(z)$ equal to zero whenever $z$ is a critical value of $f$.
\par
The function $F=F_\phi$ with $\phi$ the characteristic function of
$R^n$ is called the local singular series of $f$ and plays an
important role in number theory, for example in the circle method.
It is easy to see that $F_\phi$ is constant in the neighbourhood
of any regular value $z$; however, when $z$ tends to a singular
value of $f$, $F_\phi$ shows nontrivial singular behaviour, which
is closely related to the behaviour of the exponential sum $E(y)$.
\par
In \cite{Denef1}, the following is proven.
\begin{theorem}\label{prop:local:pol}
For each Schwartz-Bruhat function $\phi$, $F_\phi$ is a simple
$q$-exponential function.
\end{theorem}
We prove an analytic analogue of Thm.~\ref{prop:local:pol}. Fix a
dominant\footnote{By dominant we mean here that $f'(R^n)$ contains
a nonempty open subset of $K^r$.} analytic map
$f'=(f'_1,\ldots,f'_r):R^n\to K^r$, given by restricted analytic
functions $f'_i$. Let $\phi:R^n\to \Q$ be a Schwartz-Bruhat
function, and let $z\in K^r$ be a regular value of $f'$, then we
can define $F'_\phi(z)$ as
$$
F'_\phi(z)=\int_{f'(x)=z}\phi(x)|\frac{dx}{df'}|,
$$
and extend it by zero on the nonregular values $z$. We can now
state the analytic analogon of Thm.~\ref{prop:local:pol}; its
proof uses the same ideas as in the proof of
Thm.~\ref{prop:local:pol} together with some specific results on
subanalytic sets.
\begin{theorem}\label{prop:local:an}
For each Schwartz-Bruhat function $\phi:R^n\to \Q$, the function
$F'_\phi:K^r\to\Q$ is a simple $q$-exponential analytic  function.
\end{theorem}
\begin{proof}
We can cover $R^n$ by finitely many balls on which all partial
derivatives of the functions $f'_i$ are given by convergent power
series. Therefore, the set of regular values is subanalytic. Let
$i$ be the tuple $(1,2,\ldots,n-r)$ and put
$x_i=(x_1,\ldots,x_{n-r})$ and $dx_i=dx_1\wedge\ldots\wedge
dx_{n-r}$. We define $U_i\subset R^n$ as the open set given by
$h(x)\not=0$ where $h$ is the analytic function determined by
\[h(x)dx=
df'_1\wedge\ldots\wedge df'_r\wedge dx_i.
\]
The set $U_i$ is subanalytic. For a regular value $z$ and $x_i\in
R^{n-r}$, let $\Lambda_{(z,x_i)}$ be the set
\[
\{(y_1,\ldots,y_r)\in R^r \mid f'(x_i,y)=z,\ (x_i,y)\in U_i\}.
\]
When $z$ is not a regular value, we let $\Lambda_{(z,x_i)}$ be the
empty set. It follows by the implicit function theorem that each
set $\Lambda_{(z,x_i)}$ is discrete. Since $\Lambda_{(z,x_i)}$ is
also subanalytic, it is a finite set. Moreover, by Thm.~(3.2) of
\cite{DvdD} on the existence of bounds, the number of elements in
$\Lambda_{(z,x_i)}$ is uniformly bounded by a number $t>0$ when
$z$ and $x_i$ vary (here we use that the dependence on the
parameters $(z,x_i)$ is subanalytic). By the subanalytic selection
theorem (3.6) of \cite{DvdD}, there are subanalytic functions
$g_1,\ldots,g_t$ such that $\Lambda_{(z,x_i)}\subset
\{g_1(z,x_i),\ldots,g_t(z,x_i)\}$ for any regular value $z$ of
$f'$.
\par
We will express the residue form $dx/df'$ in local coordinates on
a finite cover of $f^{'-1}(z)$ for regular values $z$, using the
functions $g_j$. For a fixed regular value $z$, let $\omega_1$ be
equal to the $m-r$-form
\[
\frac{1}{h(g_1(z,x_i))}dx_i
\]
on
 \[
V_1=U_i\cap\{x=(x_i,x_{n-r+1},\ldots,x_n)\in R^n\mid
x=g_1(z,x_i)\}.
 \]
The set $V_1$ is subanalytic and has $x_i$ as local coordinates.
Subsequently, let $\omega_2$ be the $m-r$-form
\[
\frac{1}{h(g_2(z,x_i))}dx_i
\]
on
 \[
V_2=U_i\cap \{x\in R^n\mid x=g_2(z,x_i)\}.
 \]
Proceeding this way for $j=1,\ldots,t$ we find a finite cover of
$U_i\cap f^{'-1}(z)$ on each part of which we have a local
description of the residue form $dx/df'$. Now repeat this process
for all tuples $i$ in $\{1,\ldots,n\}^{n-r}$. Working with
appropriate disjoint pieces, the desired parameterized integral
can now be calculated easily, and by theorem \ref{thm:basic:int}
we find that $F_\phi$ is a simple $q$-exponential analytic
function.
\end{proof}
For the fixed maps $f$ and $f'$ and  for any Schwartz-Bruhat
functions $\phi:K^n\to \Q$ and $\phi':R^n\to\Q$ we define
exponential integrals
 \[
E_\phi(y)=\int_{K^{n}}\phi(x)\psi(\langle y,f(x)\rangle)|dx|,
 \]
and
 \[
E'_{\phi'}(y)=\int_{R^{n}}\phi'(x)\psi(\langle
y,f'(x)\rangle)|dx|.
 \]
We will use the following well-known result (for this result it is
necessary that $f$ and $f'$ are dominant).
\begin{prop}\label{prop:fourier} The functions $F_\phi$ and
$F'_{\phi'}$ are integrable over $K^r$, 
and the following Fourier transformation formulas hold
 \[
E_\phi(y)=\int_{z\in K^{r}}F_\phi(z)\psi(\langle y,z\rangle )|dz|,
 \]
 \[
 E'_{\phi'}(y)=\int_{z\in K^{r}}F'_{\phi'}(z)\psi(\langle
y,z\rangle )|dz|.
 \]
\end{prop}
This can be proven similarly as theorem 8.3.1 of
\cite{Igusa:intro}.
\section{Decay rates of exponential sums} We
use the notation of the previous section. The following theorem is
a generalization of theorem \ref{thm:decayexp:A}.
\begin{theorem}\label{thm:decayexp}
For any Schwartz-Bruhat function $\phi:K^n\to \Q$, the function
$|E_\phi(y)|$ has a nontrivial decay in terms of $|y|$, in the
sense that there are real numbers $\alpha<0$ and  $c>0$ such that
\[
|E_\phi(y)|
<c\min\{ |y|^{\alpha},1\}\quad\mbox{ for all }y\in K^r.
\]
The same statement also holds for $E'_{\phi'}$ and any
Schwartz-Bruhat function $\phi':R^n\to\Q$.
\end{theorem}
\begin{proof}
The theorem for $E_\phi$ follows from theorems
\ref{prop:local:pol} and \ref{thm:decaysimple} and proposition
\ref{prop:fourier}. For the analytic case one uses theorem
\ref{prop:local:an} instead of \ref{prop:local:pol}.
\end{proof}
\par
For a simple $q$-exponential function $G:K^r\to \Q$ which is
integrable over $K^r$, and for a fixed additive character $\psi$
on $K$, one can define the Fourier transform of $G$ as
\[
G^*:K^r\to \C:y\mapsto\int_{K^{r}}G(x)\psi(\langle
y,x\rangle)|dx|.
\]
We prove the following generalization of theorem
\ref{thm:decaysimple:A}; its proof is the technical heart of the
paper.
\begin{theorem}\label{thm:decaysimple}
Let $G:K^r\to \Q$ be a simple $q$-exponential function in the
variables $(x_1,\ldots,x_r)$, which is absolutely integrable over
$K^r$ and which has compact support. Then $|G^*(y)|$ has a
nontrivial decay in terms of $|y|$, in the sense  that there are
real numbers $\alpha<0$ and $c>0$ such that
\[
|G^*(y)|
<c\min\{\abs{(y)}^{\alpha},1\}\qquad\mbox{ for all }y\in K^r.
\]
The same statement holds also for an integrable simple analytic
$q$-exponential function with compact support $G':K^r\to\Q$, and
its Fourier transform $G^{'*}$.
\end{theorem}
\begin{proof} 
We prove the theorem in the algebraic case; the analytic case is
completely similar. It is clear that $G^*$ is uniformly bounded
since
\[
|G^*(y)|\leq \int_{K^{r}}|G(x)||dx|,
\]
thus we only have to find a decay when $|y|$ is big. It is enough
to find a decay in terms of $|y_r|$ when $|y_r|$ is big. We will
write $x=(\hat x,x_r)$ with $\hat x=(x_1,\ldots,x_{r-1})$, and we
focus on what happens if we integrate $x_r$ out.
 \par
By proposition \ref{prop:descrip:simple}, we can split $G$ up as a
sum of finitely many simple terms; we will repeat this process $r$
times. Precisely, we can partition the support of $G$ into cells
$A$ with center $c$ and coset $\lambda P_m$, such that $G|_A$ is a
finite sum of functions of the form
\[
|(x_r-c(\hat x))^a\lambda^{-a}|^\frac{1}{m}v(x_r-c(\hat
x))^{s}h(\hat x),
\]
where $h:K^{r-1}\to \Q$ is a simple $q$-exponential function,
$h(\hat x)\not=0$, $a$ and $s\geq0$ integers. We repeat this
argument to the function $h$ to partition each $\pi_{r-1}(A)$ into
cells $A'$ with center $c'$ and coset $\lambda'P_{m'}$ such that
$h|_{A'}$ is a sum of functions of the form
 \[ h(x)=|(x_{r-1}-c'(\tilde x))^{a'}\lambda^{'-a'}|^\frac{1}{m'}
 v(x_{r-1}-c'(\tilde x))^{s'}h'(\tilde x),
 \]
with $\tilde x=(x_1,\ldots, x_{r-2})$, integers $a',s'$ and a
simple $q$-exponential function $h'$ which is nowhere zero.
Repeating this construction $r$ times, we find a partition of the
support of $G$ into cells $A$ on which the function $G$ is a sum
of terms which split completely on $A$, in the sense that $A$ has
center $c$ and coset $\lambda P_m$ and $G|_A$ is a sum of
functions $H$ of the form
\begin{equation}\label{eq:term:H}
 H(x)=|(x_r-c(\hat
x))^a\lambda^{-a}|^\frac{1}{m}v(x_r-c(\hat x))^{s}h(\hat x),
\end{equation}
where $a$ and $s$ are integers, and $h$ is a nowhere zero function
which splits in a similar way on the cell $\pi_{r-1}(A)$ with
center $c'$ and so on $r$ times. By partitioning further, we may
also suppose that either $v(x_r-c)$ is constant on $A$, either it
takes infinitely many values on $A$, and similarly for
$v(x_{r-1}-c')$ and so on. In the case that $v(x_r-c)$ is constant
on $A$, we may put $a=s=0$ and similarly when $v(x_{r-1}-c')$ is
constant on $A$, and so on. We may suppose that the terms as in
equation (\ref{eq:term:H}) all have at least one different
exponent in the sequence $a,s,a',s',\ldots,s^{(r)}$. Now fix such
a term $H$ as in (\ref{eq:term:H}) and such a cell $A$ with center
$c$ and coset $\lambda P_m$.
\par
We may suppose that $A$ has dimension $r$ (as an analytic $p$-adic
manifold), hence nonzero measure. By the description
(\ref{eq:term:H}) of $H$ and by the definition of cells, the fact
that $H$ is integrable over $A$ only depends on the exponents
$a,s,a',s',\ldots,s^{(r)}$ and the particular form of the cell
$A$. Also, all terms described above have a different asymptotic
behaviour if $x\in A$ tends to a point in the boundary of $A$,
hence they cannot cancel if one integrates their sum over $A$.
Since $G$ is absolutely integrable over $A$, it follows that each
term is integrable over $A$. In particular, $H$ is integrable over
$A$.
\par
By a variant of Fubini's theorem, the function
 \[
H(\hat x,\cdot):A_{\hat x}\to \Q: x_r\mapsto H(\hat x,x_r)
 \]
is integrable over $A_{\hat x}$ for almost all $\hat x\in
\pi_{r-1}(A)$. Since this integrability only depends on $a$ and
$s$ and $A$ has nonzero measure, $H_(\hat x,\cdot)$ is integrable
for each $\hat x\in \pi_{r-1}(A)$.
\par
Suppose that $A$ has the following form
\[
\begin{array}{ll}
A=\{x\mid & \hat x\in A',\,v(\alpha(\hat x))\sq_1 v(x_r-c(\hat
x))\sq_2 v(\beta(\hat x)),
\\
 &
 x_r-c(\hat x)\in\lambda P_m \},
\end{array}
\]
where $A'=\pi_{r-1}(A)$, $\sq_i$ is $<$ or no condition, and
$\alpha,\beta:K^{r-1}\to K^\times$ and $c:K^{r-1}\to K$ are
analytic semialgebraic functions.  Since $A$ has dimension $r$, we
have $\lambda\not=0$. Since the support of $G$ is compact and
since $A$ is contained in the support of $G$ we know that the
closure of $A$ is compact. Thereby we find  that $\sq_1$ is $<$.
When $\sq_2$ stands for no condition, we will put $\beta(\hat
x)=0$, such that
$v(\beta(\hat x))=\infty$. 
\par
Write $\chi_{\lambda P_{m}}$ for the characteristic function of
$\lambda P_{m}$. For $\hat x\in A'$ and $y\in K^r$, we denote by
$I(\hat x,y)$ the value
 \[
\int_{x_r\in A_{\hat x}} H(x)\, \psi(\langle x,y\rangle )\,|dx_r|.
 \]
We put $u=x_r-c(\hat x)$ and $\hat y=(y_{1},\ldots,y_{r-1})$. We
can calculate the integral $I$ and find that $I(\hat x,y)$ equals
\begin{equation}\label{eq:sum0}
\begin{array}{c}
\psi(\langle \hat x,\hat y\rangle+c(\hat x)y_r)\ h(\hat
x)\sum\limits_j (q^{-ja}|\lambda|^{-a})^\frac{1}{m}\, j^{s}
\int\limits_{v(u)=j} \chi_{\lambda P_m}(u)\psi(u y_{r})\,|du|,
\end{array}
\end{equation}
where the summation is over those $j$ satisfying
\[
v(\alpha(\hat x)) < j \sq_2 v(\beta(\hat x)).
\]
\par
By Hensel's lemma, there exists an integer $e$ such that all units
$\alpha$ with $\alpha\equiv 1\bmod \pi_0^e$ are $m$-th powers
(here, $\pi_0$ is such that $v(\pi_0)=1$). Remark that
\[
\int_{v(u)=j}\chi_{\lambda P_m}(u)\psi(u y_{r})\,|du|
\]
is zero whenever $j<-v(y_r)-e$ (since in this case one essentially
sums a nontrivial character over a finite group). Therefore, the
only terms contributing to the sum (\ref{eq:sum0}) are those for
which $-v(y_{r})-e\leq j$.
\par
In the case where $v(\beta)$ is bounded on $A'$, the intervals of
integers $[v(\alpha),v(\beta)]$ and $[-v(y_r)-e,\infty]$ will
eventually become disjoint. Thus in this case, the function $H$ on
$A$ has, eventually, a zero contribution to the whole of the
exponential sum, when $|y_r|$ goes to $\infty$.
\par
Now suppose that $v(\beta)$ is unbounded on $A'$ (possibly
$v(\beta)=\infty$). By the above arguments, we then have
\begin{eqnarray}
|\int_{x\in A}H(x)\psi(\langle x,y\rangle)|dx|\, |
 & = &
 |\int_{\hat x\in A'} I(\hat x,y)|d\hat x|\,|\label{eq:int:I}
 \\
& \leq &
 \int_{(\ref{eq:cond:int})}|H(\hat x,u)||dx|\label{eq:int:bound}
 \\
& \leq &
 \int_{x\in A}|H(\hat x,u)||dx|\label{eq:int:bound:2}
\end{eqnarray}
 where the integral (\ref{eq:int:bound}) is over
\begin{equation}\label{eq:cond:int}
 \{(\hat x,u)\in K^r \mid (\hat x,x_r)\in A\mbox{ and } -v(y_r)-e \leq v(u) 
 \},
\end{equation}
and where still $u=x_r-c(\hat x)$. 
The  integral (\ref{eq:int:bound:2}) is a finite rational number,
since $H$ is integrable over $A$. The integral
(\ref{eq:int:bound}) has the same integrand as
(\ref{eq:int:bound:2}), but its domain of integration gets
arbitrarily small when $|y_r|$ tends to infinity. Together with
the fact that this domain of integration is contained in a compact
set (the support of $G$), this implies that (\ref{eq:int:bound})
goes to zero when $|y_r|$ goes to infinity. Also, the integral
(\ref{eq:int:bound}) only depends on $y_r$ and not on $\hat y$,
and moreover, the function $K\to \Q$ in the variable $y_r$ given
by the integral (\ref{eq:int:bound}) is a simple $q$-exponential
function by theorem \ref{thm:basic:int} (this follows from the
fact that the set where $H$ is positive, resp.~negative, is
semialgebraic). In other words, the left hand side of
(\ref{eq:int:I}) is bounded by (\ref{eq:int:bound}) which is a
simple $q$-exponential function $K\to \Q$ in the variable $y_r$
and which goes to zero whenever $|y_r|$ goes to infinity.
\par We can repeat this argument for each term $H$ of the form
(\ref{eq:term:H}) of $G|_A$ and for each cell $A$ in the partition
of the support of $G$, and find for each term a simple
$q$-exponential functions in $y_r$ yielding a respective bound as
in the calculations above.
\par
Using these bounding functions, is not difficult to construct a
semialgebraic function in the variable $y_r$ which bounds
$|G^*(y)|$ for all $y$ and which goes to zero when $|y_r|$ goes to
infinity. Thus we have found that $|G^*(y)|$ is bounded by a
single simple $q$-exponential function $K\to \Q$ in the variable
$y_r$ which goes to zero whenever $y_r$ goes to infinity. If we
combine this with lemma \ref{lemma:decay} below, the theorem is
proven.
\end{proof}
\begin{lemma}
\label{lemma:decay} Let $f:K\to\Q$ be a simple $q$-exponential
function. Suppose that if $|y|$ tends to $\infty$ then $f(y)$
converges to zero. Then there are real numbers $\alpha<0$ and
$c>0$ such that
\[
|f(y)|<c\ |y|^{\alpha}\ \mbox{ for all $|y|$ close enough to
$\infty$}.
\]
The same statement holds for simple analytic $q$-exponential
functions.
\end{lemma}
\begin{proof} Let $f:K\to\Q$ be a simple $q$-exponential function.
By proposition \ref{prop:descrip:simple}, we can partition $K$
into cells, such that on each cell $A$, with center $c$ and coset
$\lambda P_m$, the restriction $f|_A$ is a finite sum of functions
of the form
\begin{equation}\label{eq:formsimple}
y\mapsto r v(y-c)^s\ |(y-c)^a\lambda^{-a}|^\frac{1}{m},
\end{equation}
with $r$ a rational number and $a,s$ integers. Regroup these terms
in such a way that for each  $s$ and $a$, there is only one term
of the form (\ref{eq:formsimple}). Remark that in the case that
$\infty$ is in the boundary of the cell $A$, all terms have a
different asymptotic behaviour for $y$ tending to $\infty$. In
general, $\infty$ is in the boundary of finitely many cells of the
form
\[
A=\{y\in K\mid  |a| \sq |y|,\ y\in\lambda P_m\}
\]
with $\lambda\not=0$, $m>0$ an integer and $\square$ either $<$ or
no condition. (We could suppose that the center $c$ of $A$ is
zero, otherwise we could have refined the partition to obtain
this.) On such a cell $A$, suppose that $f|_A$ is a sum of terms
of the form (\ref{eq:formsimple}), regrouped like above. Since
$f|_A$ converges to zero if $y\to\infty$, and because all terms
have a different asymptotic behaviour, each of the terms must
converge to zero. Such a term clearly has a nontrivial decay  in
terms of $|y|$ like in the statement of the lemma, and we can find
a decay for $f(y)$ if $y\to \infty$ by taking the weakest decay
for all such terms and all cells having $\infty$ in its boundary.
\end{proof}

\clearemptydoublepage

\backmatter
\clearemptydoublepage

\bibliographystyle{amsplain}
\bibliography{anbib}
\clearemptydoublepage


\setcounter{page}{1}
\renewcommand{\thepage}{N\arabic{page}}

\renewcommand{\chaptermark}[1]{\markboth{#1}{}}
\renewcommand{\sectionmark}[1]{\markright{Nederlandse samenvatting}}

\setcounter{secnumdepth}{0}

\setcounter{chapter}{14}
\renewcommand{\thechapter}{\Alph{chapter}}

\chapter{Nederlandse samenvatting}\label{chap:sam:nl}
\section{$p$-adische getallen}
Zij $\Q_p$ de verzameling der $p$-adische getallen en $\Z_p$ de
verzameling der $p$-adische gehelen. We noteren met $P_n$ de
verzameling van de $n$-de machten in
$\Q_p^\times=\Q_p\setminus\{0\}$. We noemen een verzameling
$X\subset \Q_p^n$ semialgebra\"isch als $X$ kan bekomen worden
door het nemen van eindige unies, eindige doorsnedes, complementen
en lineaire projecties van Zariski-gesloten\footnote{Een
Zariski-gesloten verzameling is een nulpuntsverzameling van een
collectie veeltermen.} deelverzamelingen van $\Q_p^{n+e}$ voor
zekere $e$. Een functie $X\to Y$ tussen twee semialgebra\"ische
verzamelingen heet ook semialgebra\"isch als zijn grafiek
semialgebra\"isch is. We kunnen meteen het (chronologisch) eerste
resultaat van het gedane onderzoek formuleren, bekomen in
samenwerking met D.~Haskell.
\begin{theorema}[\cite{CH}, zie hoofdstuk \ref{chap:CH}]\label{thm:nl:bij}
Er bestaat een semialgebra\"ische bijectie tussen $\Z_p$ en
$\Z_p\setminus\{0\}$.
\end{theorema}
Dit beantwoordt een vraag gesteld door L.~B\'elair en is
belangrijk in verband met de zoektocht naar een niet triviale
Euler karakteristiek op de semialgebra\"ische verzamelingen.
Bovenstaand resultaat leert ons dat zo'n Euler karakteristiek niet
bestaat. Vele resultaten, analoog aan theorema \ref{thm:nl:bij}
worden bekomen voor andere velden, meestal velden die een Henselse
deelring bevatten of die een valuatie naar de gehele getallen
dragen. Gegeven deze $p$-adische bijectie rees onmiddellijk de
vraag naar een criterium voor het bestaan van semialgebra\"ische
bijecties tussen gegeven semialgebra\"ische verzamelingen, of naar
het bestaan van invarianten onder zulke bijecties. Deze vragen
worden volledig beantwoord met de volgende classificatie.
\begin{theorema}[\cite{C}, zie hoofdstuk
\ref{chap:classification}]\label{thm:nl:class}
\begin{sloppypar}
Gegeven twee se\-mi\-al\-ge\-bra\"ische verzamelingen  $X$ en $Y$
met oneindig veel punten, dan bestaat er een semialgebra\"ische
bijectie $X\to Y$ als en slechts als $X$ en $Y$ dezelfde dimensie
hebben.
\end{sloppypar}
\end{theorema}
De dimensie van semialgebra\"ische verzamelingen is in de jaren
tachtig ingevoerd \cite{SvdD} en volgt de meetkundige intu\"itie.
Men kan namelijk een semialgebra\"ische verzameling opsplitsen in
gladde oppervlakken ($p$-adische gladde vari\"eteiten) en dan de
grootste dimensie van deze stukken bekijken. Theorema
\ref{thm:nl:class} is dus een erg meetkundig resultaat en ook erg
verrassend, omdat er zeer veel bijecties blijken te zijn.
 \par
Een verzameling $Y\subset \Q_p^n$ wordt subanalytisch genoemd als
$Y$ kan bekomen worden door het nemen van eindig veel unies,
doorsnedes, complementen en lineaire projecties van
semialgebra\"ische deelverzamelingen van $\Q_p^{n+e}$ en van
nulpuntsverzamelingen van convergerende machtreeksen op
$\Z_p^{n+e'}\subset \Q_p^{n+e'}$. Een functie $X\to Y$ tussen twee
subanalytische verzamelingen wordt subanalytisch genoemd als de
grafiek subanalytisch is.
\par
De theorie van $p$-adische subanalytische verzamelingen is
ontwikkeld door Denef en van den Dries \cite{DvdD} in navolging
van Hironaka's theorie van subanalytische re\"ele verzamelingen en
is dus ook uiterst meetkundig van aard. In de bewijstechnieken
wordt veelvuldig gebruik gemaakt van logica, meerbepaald van
modeltheorie. Er was echter een belangrijk open probleem in deze
theorie in ver\-ge\-lij\-king met de theorie van
semialgebra\"ische verzamelingen: het ontbreken van een
celdecompositie. Deze tekortkoming wordt ingevuld in de thesis.
\begin{theorema}[\cite{Ccell}, zie hoofdstuk \ref{chap:cell}]\label{th:nl:cel}
Zij $f_1,\ldots,f_r$ subanalytische functies van een
subanalytische verzameling $X\subset \Q_p^n$ naar $\Q_p$. Dan
bestaat er een eindige partitie van $X$ in stukken van de vorm
\begin{equation}\label{eq:nl:cell}
\begin{array}{ll}
A=\{x\in\Q_p^n\mid
 &
  \hat x \in C,\ |\alpha(\hat x)|\sq_1
|x_n-c(\hat x)|\sq_2|\beta(\hat x)|,\\
 &  x_n-c(\hat x)\in\lambda P_n\},
\end{array}
\end{equation}
waarbij $\hat x=(x_1,\ldots,x_{n-1})$, $\sq_i$ is ofwel $<$ ofwel
geen voorwaarde, $\lambda$ is in $\Q_p$, $n>0$, $\alpha$ en
$\beta$ zijn subanalytische functies van $\Q_p^{n-1}$ naar
$\Q_p^\times$ en $c$ is een subanalytische functie van
$\Q_p^{n-1}$ naar $\Q_p$ zodat  $\alpha,\beta$ en $c$ analytische
zijn op $C$ en zodat $C$ gelijk is aan de projectie van $A$ op de
eerste $n-1$ co\"ordinaten en opnieuw een verzameling van dezelfde
vorm. Deze stukken kunnen zo gekozen worden dat voor een
verzameling $A$ zoals in (\ref{eq:nl:cell}) voor alle $x\in A$
geldt
\[
|f_i(x)|=|h_i(\hat x)|\, |(x_n-c(\hat
x))^{a_i}\lambda^{-a_i}|^\frac{1}{n},\quad i=1,\ldots,r,
\]
waarbij $a_i$ gehele getallen zijn en $h_i$ subanalytische
functies van $\Q_p^n$ naar $\Q_p$, analytisch op $C$ voor elke
$i$. Hierbij stellen we $a_i=0$ als $\lambda=0$ en $0^0=1$.
\end{theorema}
Verzamelingen als in (\ref{eq:nl:cell}) zijn geometrisch
eenvoudige objecten, ze worden cellen genoemd. Het zijn eerst en
vooral $p$-adische gladde vari\"eteiten (in de zin van
\cite{Bour}) en ten tweede is er een speciale variabele (namelijk
$x_n$) die op uiterst eenvoudige wijze in verband staat met de
andere variabelen. Ook de dimensie van zo'n cel kan eenvoudig
afgelezen worden uit zijn vorm. Deze celdecompositie geeft als
onmiddellijk gevolg ook een classificatie van subanalytische
verzamelingen; het bewijs is namelijk een reductie tot het
semialgebra\"ische analogon.
\begin{theorema}[\cite{Ccell}, zie hoofdstuk \ref{chap:cell}]
Gegeven twee subanalytische verzamelingen  $X$ en $Y$ met oneindig
veel punten, dan bestaat er een subanalytische bijectie $X\to Y$
als en slechts als $X$ en $Y$ dezelfde dimensie hebben.
\end{theorema}
\par
Celdecomposities zijn uiterst geschikt voor integratie van
allerlei functies. We zullen nu uitleggen hoe theorema
\ref{th:nl:cel} in de thesis gebruikt wordt om een conjectuur van
Denef te bewijzen.
\begin{definitie}
Een simpele $p$-exponenti\"ele functie $f:\Q_p^n\to\Q$ is een
functie die bekomen kan worden als $\Q$-lineaire combinatie van
producten van  functies van de vorm $x\mapsto |f(x)|$ en $x\mapsto
v(g(x))$ waarbij $f$ en $g$ semialgebra\"ische functies zijn en
$g$ nooit nul. Laat $\OOs$ de $\Q$-algebra zijn van simpele
$p$-exponenti\"ele functies.
\par
Wanneer in bovenstaande $f$ en $g$ subanalytische functies mogen
zijn spreken we analoog van simpele analytische $p$-exponenti\"ele
functies en we laten $\OOK$ de $\Q$-algebra zijn van simpele
analytische $p$-exponenti\"ele functies.
\end{definitie}
Op deze algebra's $\OOs$ en $\OOK$ kunnen we een
integratieoperator defini\"eren: voor $f:\Q_p^{n+m}\to \Q$ in
$\OOs$ defini\"eren we een functie $I_n(f):\Q_p^n\to\Q$ als
\[
I_n(f)(x)=\int_{y\in\Q_p^m}f(x,y)|dy|
\]
(waarbij $x=(x_1,\ldots,x_n)$ en $y=(y_1,\ldots, y_m)$) als deze
integrand integreerbaar is en we zetten  $I_n(f)(x)$ gelijk aan
nul als de integrand niet integreerbaar is. We kunnen natuurlijk
hetzelfde doen voor $f\in\OOL$. Denef bewijst in \cite{Denef1} dat
$\OOs$ gesloten is onder deze integratieoperator en conjectureerde
hetzelfde voor $\OOL$. Dit wordt positief beantwoord met het
volgende theorema.
\begin{theorema}[\cite{Ccell}, zie hoofdstuk \ref{chap:cell}]
De algebra $\OOL$ is gesloten onder integratie, in de zin dat voor
elke $f\in\OOL$ en $n\geq0$ de functie $I_n(f)$ ook in $\OOL$ zit.
\end{theorema}
Het bewijs van dit resultaat verloopt analoog met het
se\-mi\-al\-ge\-bra\-\"\i\-sche bewijs door Denef en maakt gebruik
van specifieke eigenschappen van subanalytische verzamelingen; het
steunt op celdecompositie.
\section{Presburgerarithmetiek}
Getaltheorie is een gecompliceerde branche van wiskunde. Precies
daarom heeft M.~Presburger in de jaren dertig een afgeslankte
versie voorgesteld, vooral voor het bestuderen van
basiseigenschappen. Deze afgeslankte arithmetiek wordt nu
Presburgerarithmetiek genoemd en wordt in de logica beschouwd als
een eenvoudige basis-structuur. Toch zijn over deze arithmetiek
nog niet alle basisresultaten gekend, onder meer omdat er
natuurlijk steeds nieuwe basisvragen ontstaan naargelang wiskunde
als wetenschap verder evolueert. Enkele basisvragen worden in deze
thesis beantwoord (zie ook \cite{CPres} en hoofdstuk
\ref{chap:Pres}).
\par
Een Presburgerverzameling $A\subset \Z^n$ is een verzameling die
bekomen kan worden door het nemen van eindige unies, intersecties,
complementen, lineaire projecties en Cartesische producten van de
verzameling $\N\subset\Z$ en van nevenklassen van deelmodules van
$\Z^{n+e}$ voor zekere $e$. Een Presburgerafbeelding is een
functie tussen twee Presburgerverzamelingen zodat zijn grafiek een
Presburgerverzameling is. In de thesis wordt een celdecompositie
geformuleerd voor Presburgerverzamelingen: elke
Presburgerverzameling kan in eindig veel stukken gepartitioneerd
worden waarbij deze stukken een eenvoudige vorm hebben, genoemd
cellen. Dit is op zich een beetje een triviaal resultaat, maar het
kan gebruikt worden om de volgende classificatie aan te tonen:
\begin{theorema}[\cite{CPres}, zie hoofdstuk \ref{chap:Pres}]
Gegeven twee oneindige Presburgerverzamelingen $X$ en $Y$, dan
bestaat er een Presburgerbijectie $X\to Y$ als en slechts als $X$
en $Y$ dezelfde dimensie hebben.
\end{theorema}
In de thesis  wordt de dimensie van Presburgerverzamelingen
ingevoerd gebruik makende van technieken uit de logica. Gelukkig
sluit de definitie van de dimensie mooi aan bij de meetkundige
intu\"itie, evenals de basiseigenschappen van de dimensie. Ook
wordt er een fundamenteel logisch resultaat bewezen over
equivalentierelaties.
\begin{theorema}[\cite{CPres}, zie hoofdstuk \ref{chap:Pres}]
Zij $R$ een equivalentierelatie op $\Z^n$, zodanig dat de
verzameling van  koppels van equivalente $x$ en $y$ een
Presburgerverzameling is. Dan bestaat er een Presburgerfunctie
\[
f:\Z^n\to\Z^m
\]
voor een zekere $m$ met de eigenschap dat $f(x)=f(y)$ als en
slechts als $x$ equivalent is met $y$.
\end{theorema}
\section{Exponenti\"ele sommen}
Als kroon op het werk worden exponenti\"ele sommen bestudeerd,
zowel in \cite{CHerr} dat niet werd opgenomen in de thesis, als in
\cite{Cexp}, zie hoofdstuk \ref{chap:decay}. Het betreft een groot
project van Igusa voor het veralgemenen van de volledige
Igusa-theorie naar het multidimensionele geval. Igusa zelf heeft
een uitgebreide theorie ontworpen, in de lijn van A.~Weil, over
functies van $n$-dimensionale ruimtes naar \'e\'endimensionale
ruimtes. Hij werkte vooral met lokaal compacte velden, dit zijn
$\R$, $\C$, $\Q_p$, $\F_p((t))$ en eindige velduitbreidingen
daarvan. Wij werken vooral met $\Q_p$ en eindige velduitbreidingen
en veralgemenen een deel van de Igusa-theorie nu ook naar functies
van $n$-dimensionale ruimtes naar $r$-dimensionale ruimtes. Bekijk
de volgende genormaliseerde exponenti\"ele som
\[
\frac{1}{p^{mn}}\sum_{x\in (\Z_{p}/p^{m})^{n}}\exp(2\pi i
\frac{f(x)}{p^{m}})
\]
met $f$ een polynoom over $\Z$, dan kunnen we dit herschrijven als
\[
\int_{x\in \Z_p^n}\psi(yf(x))|dx|,
\]
waarbij $|dx|$ de genormaliseerde Haar maat op $\Z_p^n$ voorstelt,
$\psi$ het standaard additieve karakter op $\Q_p$ is en $y=1/p^m$.
Deze exponenti\"ele som is in zekere zin \'e\'endimensionaal omdat
$y$ slechts \'e\'en $p$-adische variabele voorstelt.
\par
Laten we nu starten met een polynoomafbeelding
$f=(f_1,\ldots,f_r)$, dus $f_i$ is een polynoom over $\Z_p$ in de
variabelen $x=(x_1,\ldots,x_n)$. Noteer met
$\langle\cdot,\cdot\rangle$ het inproduct op $\Q_p^r$. We
veronderstellen van af nu altijd dat $f$ een dominante
polynoomafbeelding is, dus het beeld is niet bevat in de
nulpuntsverzameling van een polynoom. De meerdimensionale
exponenti\"ele som kunnen we nu defini\"eren als
\[
E(y)=\int_{\Z_{p}^{n}}\psi(\langle y,f(x)\rangle)|dx|
\]
voor $y\in\Q_p^r$. We bestuderen vooral het asymptotische gedrag
van $E(y)$ als $|y|$ naar oneindig gaat. Zoals gezegd werd dit ook
door Igusa bestudeerd voor $r=1$ en later ook door Lichtin voor
$r=2$ \cite{Lichtin}, \cite{Lichtin2}. In de thesis  bewijzen we
het volgende fundamentele resultaat voor willekeurige $r$.
\begin{theorema}\label{thm:nl:decayexp:A}
Er bestaat een $\alpha<0$ in $\Q$ en een constante $c$ zodanig dat
\[
|E(y)|<c\min \{|y|^{\alpha},1\}\qquad\mbox{ voor alle } y\in
\Q_p^r.
\]
\end{theorema}
Dit geeft dus een bovengrens en een convergentiegedrag voor een
heel brede klasse van exponenti\"ele sommen opgebouwd met
dominante polynoomafbeeldingen. Het bewijs is weer volledig
gebaseerd op $p$-adische celdecompositie. Een groot pluspunt van
zulk bewijs is dat het eigenlijk in zekere zin onafhankelijk wordt
van de algebra\"ische setting: alles gaat meteen ook op voor de
subanalytische setting. Meerbepaald als $f'=(f_1',\ldots,f_r')$
een analytische afbeelding is, gegeven door convergerende
machtreeksen op $\Z_p^n$ en als deze afbeelding dominant is,
kunnen we analoog als boven de meerdimensionale exponenti\"ele som
\[
E'(y)=\int_{\Z_{p}^{n}}\psi(\langle y,f'(x)\rangle)|dx|
\]
invoeren en we bekomen hetzelfde asymptotische resultaat als
voor\-heen. \par Het blijft een opdracht voor de toekomst om meer
te zeggen hoe deze $\alpha$ afhangt van $f$. Eveneens is het een
belangrijke uitdaging om uniformiteit met betrekking tot het
priemgetal $p$ na te gaan. In deze zin conjectureerde Igusa dat er
een erg scherpe $\alpha$ bestaat en een constante $c$ die
bovengrenzen bepalen zoals in de theorema's hierboven voor bijna
alle $p$ tegelijk (natuurlijk heeft dit slechts zin als alle
informatie over $\Z$ gegeven is).

\end{document}